\documentclass[11pt,oneside]{amsart}

\usepackage{tikz-cd}

\usepackage{fix-cm}

\DeclareMathAlphabet{\mathcal}{OMS}{cmsy}{m}{n}

\usepackage[arrow,curve,matrix,tips,2cell]{xy}

\usepackage[hmargin=2cm,vmargin=1.5cm]{geometry}

\usepackage{bm}

\usepackage{fancyhdr}

\usepackage[obeyspaces]{url}

\usepackage{tensor}

\usepackage{mathtools}

\usepackage{accents}

\usepackage{amsmath}
\usepackage{amscd}
\usepackage{amssymb}
\usepackage{latexsym}
\usepackage{url}

\usepackage{graphicx}

\newtheorem{theorem}{Theorem}[section]
\newtheorem*{theorem*}{Theorem}
\newtheorem{lemma}[theorem]{Lemma}
\newtheorem*{lemma*}{Lemma}
\newtheorem{corollary}[theorem]{Corollary}
\newtheorem{proposition}[theorem]{Proposition}

\newtheorem{remark}[theorem]{Remark}
\newtheorem*{remark*}{Remark}
\newtheorem{definition}[theorem]{Definition}
\newtheorem*{definition*}{Definition}
\newtheorem{question}[theorem]{Question}
\newtheorem*{question*}{Question}

\newtheorem{example}[theorem]{Example}
\newtheorem{examples}[theorem]{Examples}

\makeatletter
\def\revddots{\mathinner{\mkern1mu\raise\p@
\vbox{\kern7\p@\hbox{.}}\mkern2mu
\raise4\p@\hbox{.}\mkern2mu\raise7\p@\hbox{.}\mkern1mu}}
\makeatother % ddots gespiegelt
\newcommand{\bgl}{\begin{equation}} %eine Gleichung mit Ziffer
\newcommand{\egl}{\end{equation}}
\newcommand{\bgloz}{\begin{equation*}} %eine Gleichung ohne Ziffer
\newcommand{\egloz}{\end{equation*}}
\newcommand{\bgln}{\begin{eqnarray}} %mehrere Gleichungen mit Ziffer
\newcommand{\egln}{\end{eqnarray}}
\newcommand{\bglnoz}{\begin{eqnarray*}} %mehrere Gleichungen ohne Ziffer
\newcommand{\eglnoz}{\end{eqnarray*}}
\newcommand{\btheo}{\begin{theorem}}
\newcommand{\etheo}{\end{theorem}}
\newcommand{\btheooz}{\begin{theorem*}}
\newcommand{\etheooz}{\end{theorem*}}

\newcommand{\blemma}{\begin{lemma}}
\newcommand{\elemma}{\end{lemma}}
\newcommand{\blemmaoz}{\begin{lemma*}}
\newcommand{\elemmaoz}{\end{lemma*}}
\newcommand{\bproof}{\begin{proof}}
\newcommand{\eproof}{\end{proof}}
\newcommand{\bbew}{\begin{beweis}}
\newcommand{\ebew}{\end{beweis}}
\newcommand{\bremark}{\begin{remark}}
\newcommand{\eremark}{\end{remark}}
\newcommand{\bremarkoz}{\begin{remark*}}
\newcommand{\eremarkoz}{\end{remark*}}
\newcommand{\bdefin}{\begin{definition}}
\newcommand{\edefin}{\end{definition}}
\newcommand{\bdefinoz}{\begin{definition*}}
\newcommand{\edefinoz}{\end{definition*}}
\newcommand{\bex}{\begin{example}}
\newcommand{\eex}{\end{example}}
\newcommand{\bexs}{\begin{examples}}
\newcommand{\eexs}{\end{examples}}

\newcommand{\bprop}{\begin{proposition}}
\newcommand{\eprop}{\end{proposition}}
\newcommand{\bcor}{\begin{corollary}}
\newcommand{\ecor}{\end{corollary}}
\newcommand{\bfa}{\begin{cases}} %Fallunterscheidung
\newcommand{\efa}{\end{cases}}

\newcommand{\bquestion}{\begin{question}}
\newcommand{\equestion}{\end{question}}
\newcommand{\bquestionoz}{\begin{question*}}
\newcommand{\equestionoz}{\end{question*}}
\newcommand{\cG}{\mathcal G}

\newcommand{\cJ}{\mathcal J}
\newcommand{\cK}{\mathcal K}

\newcommand{\cO}{\mathcal O}

\def\Cz{\mathbb{C}}

\def\Nz{\mathbb{N}}

\def\Rz{\mathbb{R}}

\def\Zz{\mathbb{Z}}

\def\1z{\mathbb{1}}
\newcommand{\fE}{\mathfrak E}

\newcommand{\mfd}{\mathfrak d}
\newcommand{\mfe}{\mathfrak e}
\newcommand{\mff}{\mathfrak f}
\newcommand{\mfl}{\mathfrak l}

\newcommand{\an}[1]{``#1''} % Anfuehrungsstriche
\newcommand{\ti}{\tilde}

\newcommand{\ma}{\mapsto} % wird abgebildet auf
\newcommand{\into}{\hookrightarrow} % injektiv
\newcommand{\Rarr}{\Rightarrow} % Folgerung
\newcommand{\sgn}{\rm sgn}

\def\SEMI{\mbox{$\times\kern-2pt\vrule height5pt width.6pt \kern3pt $}}

\newcommand{\Spec}{{\rm Spec\,}} % Spektrum

\newcommand{\defeq}{\mathrel{:=}} % per Definition
\newcommand{\eqdef}{\mathrel{=:}} % auch per Definition

\newcommand{\dop}{\text{: }} % in Mengen
\newcommand{\lge}{\left\{} % links geschweift
\newcommand{\rge}{\right\}} % rechts geschweift
\newcommand{\lsp}{\left\langle} % links spitz
\newcommand{\rsp}{\right\rangle} % links spitz
\newcommand{\gekl}[1]{\lge #1 \rge} % geschweifte Klammer
\newcommand{\spkl}[1]{\lsp #1 \rsp} % spitze Klammer
\newcommand{\menge}[2]{\gekl{ #1 \dop #2 }} % Menge
\def\bf1{\mathbf{1}}
\newcommand{\bo}{\bm{o}}
\newcommand{\beps}{\bm{\epsilon}}
\newcommand{\chT}{\check{T}}
\newcommand{\Zzg}{{\Zz_{\geq 0}}}
\newcommand{\Stab}{\rm Stab}

\newcommand{\Conv}{\mathop{\scalebox{2}{\raisebox{-0.2ex}{$\ast$}}}}

\newcommand{\pars}{\setlength{\parindent}{0cm} \setlength{\parskip}{0.5cm}}
\newcommand{\pari}{\setlength{\parindent}{0.5cm} \setlength{\parskip}{0cm}}
\newcommand{\nopar}{\setlength{\parindent}{0cm} \setlength{\parskip}{0cm}}

\usepackage{newtxtext}
\usepackage{newtxmath}

\begin{document}

\title{Semigroup C*-algebras arising from graphs of monoids}

\thispagestyle{fancy}

\author{Cheng Chen}

\author{Xin Li}

\address{Xin Li, School of Mathematics and Statistics, University of Glasgow, University Place, Glasgow G12 8QQ, United Kingdom}
\email{Xin.Li@glasgow.ac.uk}

\address{Cheng Chen, School of Mathematics and Statistics, University of Glasgow, University Place, Glasgow G12 8QQ, United Kingdom}
\email{2527724c@student.gla.ac.uk}

\subjclass[2010]{Primary 46L05, 20E08; Secondary 46L80, 46L35}

\thanks{This project has received funding from the European Research Council (ERC) under the European Union's Horizon 2020 research
and innovation programme (grant agreement No. 817597). The first author has been partially supported by a China Scholarship Council (CSC) PhD Scholarship.}

\begin{abstract}
We study groupoids and semigroup C*-algebras arising from graphs of monoids, in the setting of right LCM monoids. First, we establish a general criterion when a graph of monoids gives rise to a submonoid of the fundamental group which is right LCM. Moreover, we carry out a detailed analysis of structural properties of semigroup C*-algebras arising from graphs of monoids, including closed invariant subspaces and topological freeness of the groupoids as well as ideal structure, nuclearity and K-theory of the semigroup C*-algebras. As an application, we construct families of pairwise non-conjugate Cartan subalgebras in every UCT Kirchberg algebra.
\end{abstract}

\maketitle

%\tableofcontents

\setlength{\parindent}{0cm} \setlength{\parskip}{0.5cm}

\section{Introduction}

Bass-Serre theory \cite{Serre,Bass} plays an important role in group theory and geometric topology. Roughly speaking, it builds a precise dictionary between group actions on trees and decompositions of groups as fundamental groups of graphs of groups, which are constructed from vertex stabilizers and edge stabilizers of corresponding actions on trees. This dictionary has become a standard tool, with many applications, and it has been generalized in several directions.

The goal of the present paper is to study right LCM submonoids of fundamental groups, which arise from graphs of monoids, and the groupoids and C*-algebras generated by left regular representations of these submonoids. In general, an arbitrary submonoid does not necessarily admit nice presentations, but we identify a setting when our submonoids can be described by essentially the same presentations as the fundamental groups attached to the corresponding graph of groups. The crucial idea is that if all the groups in our graph of groups are totally ordered, and if the structure maps are order preserving or order reversing, then the positive cones form a graph of monoids which naturally leads to a submonoid of the fundamental group. Such monoids and the C*-algebras generated by their left regular representations have been studied in special situations, for instance for free products \cite{DP}, the particular case of Baumslag-Solitar monoids \cite{Spi}, or more general HNN extensions \cite{HNSY}. Our more general framework allows for a rich supply of examples which are natural, yet tractable, and exhibit interesting phenomena.

For more information about semigroup C*-algebras in general, we refer the reader to \cite{CELY} and the references therein. Generally speaking, it is interesting that, unlike in the case of groups, the structure of a semigroup automatically produces a topological dynamical system (in terms of a topological groupoid) which on the one hand serves as a model for the C*-algebra generated by the left regular representation of our semigroup and on the other hand interacts with the algebraic structure of the semigroup.

In our case of graphs of monoids, dynamical properties of the associated groupoid are intimately connected to the behaviour of normal forms of finite and infinite words in the generators under left multiplication by group elements, which in turn is closely related --- very much in the spirit of Bass-Serre theory --- to the corresponding group action on the Bass-Serre tree and its boundary. Following this philosophy, we succeed in analyzing the dynamical properties of the groupoids arising from graphs of monoids and deduce results about structural properties of the corresponding semigroup C*-algebras, both in a general context and for specific example classes. Our results on structural properties and classification for groupoids and C*-algebras attached to our graphs of monoids show a fascinating interplay between algebraic (in terms of normal forms), geometric (via action on trees), dynamic (through groupoids, their orbit structures and related properties) and analytic (C*-algebraic) aspects. Let us now summarize our main achievements in the general setting:
\nopar

\begin{itemize}
\item We identify a general criterion when our graph of monoids gives rise to a submonoid of the fundamental group which is right LCM, i.e., non-empty intersections of principal right ideals are again principal right ideals (Proposition~\ref{prop:LCM}). This right LCM property is crucial for a detailed analysis of the groupoids and semigroup C*-algebras attached to graphs of monoids.
\item We establish a general criterion, in the right LCM case, when the semigroup C*-algebra attached to a graph of monoids is purely infinite simple (Corollary~\ref{cor:I,pis}).
\item We show, again in the right LCM case, that the semigroup C*-algebra of a graph of monoids is nuclear if and only if this is the case for the graph of monoids obtained by restricting to a maximal subtree (Theorem~\ref{thm:NucP-PT}).
\end{itemize}
\pars

In the special case where all the individual monoids of our graph of monoids are positive cones of subgroups of $(\Rz,+)$, we obtain, under some extra assumptions, the following stronger results:
\nopar

\begin{itemize}
\item We determine all closed invariant subspaces of the unit space of the groupoid attached to our graph of monoids (Theorems~\ref{thm:clinvsub_gen} and \ref{thm:clinvsub_GBS}).
\item For each closed invariant subspace, we characterize when the restricted groupoid is topologically free (Theorem~\ref{thm:topfree}). This leads to a criterion when ideals of the semigroup C*-algebra attached to a graph of monoids are in one-to-one correspondence to closed invariant subspaces of the corresponding groupoid (Corollary~\ref{cor:ideals=clinvsub}).
\item We completely characterize which graphs of monoids give rise to amenable groupoids, or equivalently, nuclear semigroup C*-algebras (Theorem~\ref{thm:CPTnuc}).
\item We compute K-theory for all groupoid C*-algebras induced from closed invariant subspaces (Theorem~\ref{thm:K}).
\item We establish a criterion when boundary quotients of semigroup C*-algebras arising from graphs of monoids are UCT Kirchberg algebras and compute their K-groups, and thus classify them completely (Theorem~\ref{thm:bC}). The boundary quotient is a distinguished quotient of the semigroup C*-algebra corresponding to a minimal closed invariant subspace of the underlying groupoid. In typical examples, boundary quotients play a special role and have very interesting properties. For instance, they are more likely to be simple and hence fall into the scope of classification results.
\end{itemize}
\pars

Our results cover many concrete examples, for instance free products or amalgamated free products. At the same time, it is surprising that even though our framework allows for a very rich supply of examples, we are still able to obtain general results on important structural properties such as closed invariant subspaces, topological freeness, nuclearity, K-theory and classifiability.

As an application of our results, we construct families of pairwise non-conjugate Cartan subalgebras in all UCT Kirchberg algebras.
\btheooz[see Theorem~\ref{thm:CartanInKirchberg}]
Let $A$ be a UCT Kirchberg algebra. For every abelian, torsion-free, finite rank group $\Gamma$ which is not free abelian, there exists a Cartan subalgebra $B_{\Gamma}$ of $A$ such that $\Spec B_{\Gamma}$ is homeomorphic to the Cantor space if $A$ is unital and to the non-compact locally compact Cantor space if $A$ is not unital, and, for all such groups $\Gamma$ and $\Lambda$, $(A,B_{\Gamma}) \cong (A,B_{\Lambda})$ implies $\Gamma \cong \Lambda$.
\etheooz
\nopar

Here \an{non-compact locally compact Cantor space} refers to the up to homeomorphism unique totally disconnected, second countable, locally compact non-compact Hausdorff space without isolated points. Cartan subalgebras of C*-algebras have been introduced in \cite{Ren}, based on \cite{Kum}. They provide a general framework, in C*-algebraic terms, for producing groupoid models for C*-algebras. It was shown in \cite{LR} that Cartan subalgebras in UCT Kirchberg algebras are not unique. Our theorem above is a strengthening of this result. It not only implies that every UCT Kirchberg algebra has uncountably many pairwise non-conjugate Cartan subalgebras, but actually produces the more precise statement that given a UCT Kirchberg algebra $A$, the classification of Cartan subalgebras of $A$ is at least as complex as the classification of all abelian, torsion-free, finite rank groups. The complexity of the latter classification problem has been discussed in \cite{Tho}.
\pars

\bremarkoz
We point out that another construction of C*-algebras from graphs of groups has been studied in \cite{BMPST}. However, monoids do not feature in \cite{BMPST}, and our results show that the two constructions have different properties.
\eremarkoz

We would like to thank the anonymous referees for very helpful comments which helped to improve this paper.

This paper is partly based on contents of the PhD thesis of the first author, completed at Queen Mary University of London and the University of Glasgow.

\section{Graphs of groups and monoids}

\subsection{Preliminaries}

We collect some basics about graphs of groups and their fundamental groups. Our exposition follows \cite{Serre} (see also \cite{Bass}). Here and in the sequel, we write $\beps$ for the identity of a group.

A graph of groups consists of a graph $(V,E)$, where $V$ is the set of vertices and $E$ is the set of edges. Given $e \in E$, $o(e) \in V$ denotes its origin and $t(e)$ denotes its target. Edges come in pairs $e, \bar{e} \in E$ (as in \cite{Serre}). We assume throughout that $(V,E)$ is connected. Part of the definition of a graph of groups is the following data: For every $v \in V$, we are given a group $G_v$, and for every edge $e \in E$, we are given a group $G_e$ such that $G_e = G_{\bar{e}}$. Moreover, we are given group embeddings $G_e \into G_{o(e)}, \, x \ma x^{\bar{e}}$ and $G_e \into G_{t(e)}, \, x \ma x^{e}$. Let $G_e^{\bar{e}}$ denote the image of $G_e \into G_{o(e)}, \, x \ma x^{\bar{e}}$ and $G_e^e$ denote the image of $G_e \into G_{t(e)}, \, x \ma x^{e}$.

Let $\bo$ be a fixed vertex which we think of as the base vertex. 
\bdefin
\label{def:o-word}
An $\bo$-word is a word $X$ of the form $X = h_0 e_1 h_1 \dotso e_n h_n$, where $h_l \in G_{v_l}$, $e_l \in E$ satisfies $o(e_l) = v_{l-1}$, $t(e_l) = v_l$, and $v_0 = \bm{o} = v_n$.  We set $\ell(X) \defeq n$.
\edefin
\bdefin
The fundamental group $\pi_{1,\bo}$ attached to our graph of groups is the subgroup of
$$
 F \defeq \spkl{\menge{G_v}{v \in V} \cup E \ \vert \ \bar{e} = e^{-1} \ \forall \ e \in E, \, g^{\bar{e}} e = e g^e \ \forall \ g \in G_e, \, e \in E}
$$
generated by all elements of $F$ which can be represented by $\bo$-words.
\edefin
Here and in the sequel, given two $\bo$-words $X$ and $X'$, we write $X=X'$ if these words represent the same element of $\pi_{1,\bo}$, while we write $X \equiv X'$ if the words are identical.

Let us now turn to normal forms. 
\bdefin
An $\bo$-word $X = h_0 e_1 h_1 \dotso e_n h_n$ is called reduced if $n=0$ (and we allow $h_0 = \beps$), or $n \geq 1$ and $e_l = \bar{e}_{l-1}$ implies $h_l \notin G_{e_{l-1}}^{e_{l-1}}$ for all $1 \leq l \leq n$.
\edefin

The following is proved in \cite[Chapter~1, \S~5.2]{Serre}.
\btheo
\label{thm:normal}
Every element of $\pi_{1,\bo}$ can be represented by a reduced $\bo$-word.
\pari

Given two reduced $\bo$-words $X = h_0 e_1 h_1 \dotso e_n h_n$ and $X' = h'_0 e'_1 h'_1 \dotso e'_{n'} h'_{n'}$, we have $X=X'$ if and only if $n=n'$, $e_l = e_{l'}$ for all $1 \leq l \leq n$, and there exist $a_l \in G_{e_l}$, $1 \leq l \leq n$ such that 
$$
 h'_0 = h_0 a_1^{\bar{e}_1}, \, a_l^{e_l} h'_l = h_l a_{l+1}^{\bar{e}_{l+1}} \ \forall \ 1 \leq l \leq n-1, \, a_n^{e_n} h'_n = h_n.
$$
\etheo
\pars

Now assume that there is a decomposition $E = A \amalg \bar{A} \amalg T$, where $T$ is a maximal tree of $(V,E)$. Give two vertices $v, w$, we write $[v,w]$ for the geodesic in $T$ from $v$ to $w$. The following gives a presentation for $\pi_{1,\bo}$.
\bdefin
Define
$$
 \pi_{1,T} \defeq \spkl{\menge{G_v}{v \in V} \cup A \ \vert \ g^{\bar{e}} e = e g^e \ \forall \ g \in G_e, \, e \in A, \, g^{\bar{e}} = g^e \ \forall \ g \in G_e, \, e \in T}.
$$
\edefin

\bprop[{\cite[Chapter~1, \S~5.1, Proposition~20]{Serre}}]
The composite $\pi_{1,\bo} \to F \to \pi_{1,T}$ is an isomorphism, where the first map is the canonical inclusion and the second map is the canonical projection.
\eprop

As explained in \cite[\S~5]{Serre}, there is a one-to-one correspondence between group actions on trees and group representations as fundamental groups of graphs of groups such that the groups $G_v$, $v \in V$, are identified with vertex stabilizers, and the groups $G_e$, $e \in E$, are identified with edge stabilizers of the corresponding action on a tree. This explains why there are embeddings $G_e \into G_{o(e)}$ and $G_e \into G_{t(e)}$, because edge stabilizers naturally embed into vertex stabilizers.

\bexs
\label{ex:pi}
\begin{enumerate}
\item[(i)] If $(V,E)$ is a tree and $G_e=\{\beps\}$ for all edges $e \in E$, then the fundamental group $\pi_{1,\bo} \cong \pi_{1,T}$ is the free product of all the groups $G_{v},\ v \in V$. More generally, without the assumption that $G_e=\{\beps\}$ for all edges $e \in E$, we obtain amalgamated free products.
\item[(ii)] If $(V,E)$ is a bouquet of circles, i.e., $V = \gekl{v}$, and $G_{v} \cong \mathbb{Z}$, $G_{e}\cong \mathbb{Z}$ for all edges $e\in E$, then the fundamental group $\pi_{1,\bo} \cong \pi_{1,T}$ is called a one vertex generalised Baumslag-Solitar (one vertex GBS) group. In that case, the fundamental group admits the presentation
$$
 \langle \gekl{b} \cup A \ \vert \ b^{n_e} e = e b^{\sgn(e) m_e} \ \forall \ e \in A \rangle.
$$
Here $A \subseteq E$ is such that $E = A \amalg \bar{A}$ (and $T = \emptyset$), $n_e, m_e \in \Zz_{\geq 1}$ for all $e \in A$, and $\sgn(e) \in \gekl{\pm 1}$ for all $e \in A$. In particular, we obtain the classical Baumslag-Solitar groups if $\# A=1$.
\end{enumerate}
\eexs
Note that HNN extensions, as studied for instance in \cite{HNSY}, also fit naturally into our framework.

In the following, it will be more convenient to work in $\pi_{1,T}$. Therefore, we set $G \defeq \pi_{1,T}$ and set up terminology in order to make use of Theorem~\ref{thm:normal}.
\bdefin
\label{def:cT-word}
A $\check{T}$-word is a word $W$ of the form 
$$
 W = g_{0,1} \dotsm g_{0,\mu_0} d_1 g_{1,1} \dotsm g_{1,\mu_1} d_2 g_{2,1} \dotsm g_{m-1,\mu_{m-1}} d_m g_{m,1} \dotsm g_{m,\mu_m},
$$
where $g_{k,\lambda} \in G_{u_{k,\lambda}}$ for all $0 \leq k \leq m$, $1 \leq \lambda \leq \mu_k$, and $d_k \in A \amalg \bar{A}$ for all $1 \leq k \leq m$.
\edefin
\nopar

Note that $\mu_k = 0$ is allowed, in which case no $g_{k,\lambda}$ appears between $d_k$ and $d_{k+1}$.
\pars

As above, given two $\check{T}$-words $W$ and $W'$, we write $W=W'$ if these words represent the same element of $G$, while we write $W \equiv W'$ if the words are identical. 

\bremark
Every word in $\gekl{G_v}_{v \in V} \cup E$ represents an element of $G$: Given such a word, first delete all letters which lie in $T$ to obtain a word in $\gekl{G_v}_{v \in V} \cup (A \amalg \bar{A})$, which in turn represents an element of $G$.
\eremark

The following map allows us to pass from $\check{T}$-words to $\bo$-words.
\bdefin
Given a $\chT$-word $W$ as in Definition~\ref{def:cT-word}, define
$$
 \fE(W) \defeq e_{0,0} g_{0,1} e_{0,1} \dotsm g_{0,\mu_0} e_{1,0} g_{1,1} \dotsm g_{1,\mu_1} e_{2,0} g_{2,1} \dotsm g_{m-1,\mu_{m-1}} e_{m,0} g_{m,1} \dotsm g_{m,\mu_m} e_{m,\mu_m},
$$
where $e_{k,\lambda} \defeq [u_{k,\lambda},u_{k,\lambda+1}]$ if $0 \leq k \leq m$, $1 \leq \lambda \leq \mu_k$ with $(k,\lambda) \neq (m,\mu_m)$, $e_{0,0} \defeq [\bo,u_{0,1}]$, $e_{m,\mu_m} \defeq [u_{m,\mu_m},\bo]$, $e_{k,0} \defeq [u_{k-1,\mu_{k-1}}, o(d_k)] d_k [t(d_k), u_{k,1}]$ if $1 \leq k \leq m$, and if $\mu_k = 0$, then we set $u_{k,0} \defeq u_{k,1} \defeq t(d_k)$ if $1 \leq k \leq m-1$ and $u_{k,0} \defeq u_{k,1} \defeq \bo$ if $k \in \gekl{0,m}$. 
\pari

We view $\fE(W)$ as an $\bo$-word by filling up with $\beps$ if necessary.

We set $\ell(W) \defeq \ell(\fE(W))$. 

A $\check{T}$-word $W$ is called reduced if $\fE(W)$ is a reduced $\bm{o}$-word.
\edefin
\pars

Note that we obtain the following by analyzing normal forms.
\bcor[{\cite[Corollary~1.14]{Bass}}]
For any connected subgraph $(V',E')$ of $(V,E)$ with maximal subtree $T' \subseteq T$, if $\pi'_{1,T'}$ denotes the fundamental group of the graph of groups given by $G'_{v'} \defeq G_{v'}$, $G'_{e'} \defeq G_{e'}$ and the same group embeddings as for the original graph of groups, then the canonical map $\pi'_{1,T'} \to \pi_{1,T}$ is a group embedding.
\ecor
\nopar

This allows for the following conventions, which will be convenient.
\bdefin
Let $G_T \defeq \spkl{ \menge{G_v}{v \in V} } \subseteq G$. 
\pari

A $\check{T}$-word in compact form is a word of the form $g_0 d_1 g_1 \dotso d_m g_m$, where $g_k \in G_T$ or $g_k = \emptyset$, and $d_k \in A \amalg \bar{A}$.
\edefin
\nopar
By definition, we can pass from a $\chT$-word as in Definition~\ref{def:cT-word} to one in compact form by setting $g_k \defeq g_{k,1} \dotsm g_{k,\mu_k}$.
\pars

The following is an immediate consequence of Theorem~\ref{thm:normal}.
\blemma
\label{lem:W=W'}
\begin{enumerate}
\item[(i)] Let $X = h_0 e_1 h_1 \dotso e_n h_n$ and $X' = h'_0 e'_1 h'_1 \dotso e'_{n'} h'_{n'}$ be two reduced $\bm{o}$-words with $h_l \in G_{v_l}$ and $h'_{l'} \in G_{v_{l'}}$ for all $l, l'$. If $X = X'$, then $n = n'$, $e_l = e'_l$ for all $1 \leq l \leq n$, and for all $1 \leq l \leq n$, there exists $a \in G_{e_l}^{e_l}$ such that $h_0 e_1 h_1 \dotso h_{l-1} e_l = h'_0 e_1 h'_1 \dotso h'_{l-1} e_l a$, and for all $1 \leq l < n$, there exists $\ti{a} \in G_{e_{l+1}}^{\bar{e}_{l+1}}$ such that $p_0 e_1 p_1 \dotso e_l p_l = p'_0 e_1 p'_1 \dotso e_l p'_l \ti{a}$.
\item[(ii)] Let $W = p_0 d_1 p_1 \dotso d_m p_m$ and $W' = p'_0 d'_1 p'_1 \dotso d'_{m'} p'_{m'}$ be two $\chT$-words in compact form. If $W = W'$, then $m = m'$, $d_k = d'_k$ for all $1 \leq k \leq m$, and for all $1 \leq k \leq m$, there exists $a \in G_{d_k}^{d_k}$ such that $q_0 d_1 q_1 \dotso q_k d_k = q'_0 d_1 q'_1 \dotso q'_k d_k a$.
\end{enumerate}
\elemma

\subsection{Presentation and normal forms for submonoids}
\label{ss:presentation}

In the following, we assume that for all $v$, $G_v$ is totally ordered with positive cone $P_v$, i.e., $G_v = P_v \cup P_v^{-1}$ and $P_v \cap P_v^{-1} = \gekl{\bm{\epsilon}}$. Define, for $e \in E$, $P_e \defeq \menge{g \in G_e}{g^e \in P_{t(e)}}$ and $P_{\bar{e}} \defeq \menge{g \in G_{\bar{e}} = G_e}{g^{\bar{e}} \in P_{o(e)}}$. Note that $P_e^e = G_e^e \cap P_{t(e)}$ and $P_{\bar{e}}^{\bar{e}} = G_{\bar{e}}^{\bar{e}} \cap P_{o(e)}$. We furthermore assume that $P_e = P_{\bar{e}}$ for all $e \in T$, and that there is a decomposition $A = A_+ \amalg A_-$ such that, for all $e \in A_+$, we have $P_{\bar{e}}^e \subseteq P_{t(e)}$, and that for all $e \in A_-$, $P_{\bar{e}}^e \subseteq P_{t(e)}^{-1}$. Note that this implies that $P_{\bar{e}} = P_e$ for all $e \in A_+$ and $P_{\bar{e}} = P_e^{-1}$ for all $e \in A_-$. 

\bdefin
Define the following submonoids of $G$:
\begin{eqnarray*}
 P &\defeq& \spkl{ \menge{P_v}{v \in V} \cup A }^+ \subseteq G\\
 P_T &\defeq& \spkl{ \menge{P_v}{v \in V} }^+ \subseteq P.
\end{eqnarray*}
\edefin
Here $\spkl{\bm{X}}^+$ denotes the submonoid generated by a subset $\bm{X}$. 

We focus on positive cones in totally ordered groups and assume that the embeddings $(\cdot)^e$, $(\cdot)^{\bar{e}}$, for $e \in T \amalg A$, are order preserving or order reversing, because this will guarantee a presentation for $P$ which is analogous to the one for $G$. This in turn leads to normal forms for elements of $P$.

\bdefin
\label{def:posword}
A positive word is a $\chT$-word as in Definition~\ref{def:cT-word} with $g_{k,\lambda} \in P_{u_{k,\lambda}}$ for all $0 \leq k \leq m$, $1 \leq \lambda \leq \mu_k$, and $d_k \in A$ for all $1 \leq k \leq m$.
\edefin

In the sequel, we denote by $\spkl{\Sigma \ \vert \ R}^+$ the universal monoid given by generators $\Sigma$ subject to relations $R$. For $P$ arising from positive cones as above, we have the following:
\blemma
\label{lem:normal}
\begin{enumerate}
\item[(i)] Every $p \in P$ is represented by a reduced positive word.
\item[(ii)] Assume that if we have two reduced positive words $W$ and $W'$ with $\fE(W) = q_0 e_1 q_1 \dotso e_n q_n$ and $\fE(W') = q'_0 e'_1 q'_1 \dotso e'_{n'} q'_{n'}$ such that $W = W'$. Then $n = n'$, $e_l = e'_l$ for all $1 \leq l \leq n$, and either $q'_0 = q_0 a_0$, where $a_0 = a^{\bar{d}_1}$, and
\begin{itemize}
\item if $e_1 \in T \cup A_+$, then we must have $a_0 e_1 = e_1 a_1$, where $a_1 = a^{e_1}$, and $a_1 q'_1 e_2 q'_2 \dotso = q_1 e_2 q_2 \dotso$,
\item if $e_1 \in A_-$, then we must have $a_0 e_1 a_1 = e_1$, where $a_1 = (a^{e_1})^{-1}$, and $a_1 q_1 e_2 q_2 \dotso = q'_1 e_2 q'_2 \dotso$;
\end{itemize}
or $q_0 = q'_0 a_0$, where $a_0 = a^{\bar{e}_1}$, and
\begin{itemize}
\item if $e_1 \in T \cup A_+$, then we must have $a_0 e_1 = e_1 a_1$, where $a_1 = a^{e_1}$, and $a_1 q_1 e_2 q_2 \dotso = q'_1 e_2 q'_2 \dotso$,
\item if $e_1 \in A_-$, then we must have $a_0 e_1 a_1 = e_1$, where $a_1 = (a^{e_1})^{-1}$, and $a_1 q'_1 e_2 q'_2 \dotso = q_1 e_2 q_2 \dotso$.
\end{itemize}
\item[(iii)] We have canonical isomorphisms
\begin{eqnarray}
\label{e:presentation}
P 
&\cong& 
  \spkl{\menge{P_v}{v \in V} \cup A \ \vline 
  \begin{array}{cl}
  a^e = a^{\bar{e}} & \forall \ a \in P_e = P_{\bar{e}}, \, e \in T;
  \\
  e a^e = a^{\bar{e}} e & \forall \ a \in P_e, \, e \in A_+;
  \\
  e = a^{\bar{e}} e (a^e)^{-1} & \forall \ a \in P_{\bar{e}}, \, e \in A_-
  \end{array}}^+.
\\
P_T
&\cong& 
  \spkl{\menge{P_v}{v \in V} \ \vert \ 
  a^e = a^{\bar{e}} \ \forall \ a \in P_e = P_{\bar{e}}, \, e \in T
  }^+.
\nonumber
\end{eqnarray}
\end{enumerate}
\elemma
\bproof
(i) Let $W$ be a positive $\check{T}$-word such that $\fE(W) = q_0 e_1 q_1 \dotso e_n q_n$. We proceed inductively on $n = \ell(\fE(W))$. If $\fE(W)$ is not reduced, then we must have $n \geq 1$ and there exists $l$ such that $e_l = \bar{e}_{l-1}$, $q_{l-1} \in P_{e_{l-1}}^{e_{l-1}}$, say $q_{l-1} = a^{e_{l-1}}$ for some $a \in P_{e_{l-1}}$. $e_l = \bar{e}_{l-1}$ implies that $e_l, e_{l-1} \in T$ as all $e \in A$ which appear in $\fE(W)$ must have already appeared in $W$. So
\begin{equation*}
 W = q_0 \dotso e_{l-1} q_{l-1} e_l \dotso q_n = q_0 \dotso e_{l-1} a^{e_{l-1}} \bar{e}_{l-1} \dotso q_n = q_0 \dotso a^{\bar{e}_{l-1}} \dotso q_n.
\end{equation*}
In this way, we arrive at a word whose length has decreased, so that it can be represented by a reduced positive $\check{T}$-word. 

(ii) Assume that $e_1 \in T \cup A_+$. If $q'_0 = q_0 a^{\bar{e}_1}$ for some $a \in G_{e_1}$, then
$$
 q_0 e_1 q_1 \dotso = q'_0 e_1 q'_1 \dotso = q_0 a^{\bar{e}_1} q'_1 \dotso = q_0 e_1 a^{e_1} q'_1
$$
implies that $q_1 e_2 q_2 \dotso = a^{e_1} q'_1 e_2 q'_2 \dotso$. If $a \in P_{e_1}$, set $a_0 \defeq a^{\bar{e}_1}$. If $a \in P_{e_1}^{-1}$, then $a^{-1} \in P_{e_1}$, and we set $a_0 \defeq (a^{-1})^{e_1}$. Since $q'_0 = q_0 a^{\bar{e}_1}$, we have $q_0 = q'_0 (a^{-1})^{\bar{e}_1} = q'_0 a_0$ and $(a^{-1})^{e_1} q_1 e_2 q_2 \dotso = q'_1 e_2 q'_2 \dotso$.

(iii) Let $Q$ be the semigroup on the right-hand side of \eqref{e:presentation}. We need to show that for two words $W$ and $W'$ in $\gekl{P_v}_{v \in V} \cup A$, $W = W'$ in $P$ implies $W = W'$ in $Q$. So assume that $W = W'$ in $P$. Without loss of generality we may assume that $\fE(W)$ and $\fE(W')$ are reduced and of the same form as in (ii). Now we proceed inductively on $m = \ell(\fE(W)) = \ell(\fE(W'))$. We treat the first case in (ii), the other cases are similar. We have the following equation in $Q$: 
$$
 W' = q_0 a_0 e_1 q'_1 \dotso = q_0 e_1 a_1 q'_1 \dotso = q_0 e_1 q_1 \dotso \equiv W,
$$
where we used that $a_1 q'_1 \dotso = q_1 \dotso$ in $P$ implies $a_1 q'_1 \dotso = q_1 \dotso$ in $Q$ by induction hypothesis.

The presentation for $P_T$ follows.
\eproof

We derive the following consequence.
\bcor
\label{cor:NoUnits}
We have $P \cap P^{-1} = \gekl{\bm{\epsilon}}$.
\ecor
\bproof
Assume $p, q \in P$ with $p, q \neq \bm{\epsilon}$ satisfy $pq = \bm{\epsilon}$. Then Lemma~\ref{lem:normal} implies that $p, q \in P_T$. Let $p_0 \dotso p_m$ be a reduced positive word representing $p$, with $p_k \in P_{v_k}$. Similarly, let $q_0 \dotso q_n$ be a reduced positive word representing $q$, with $q_l \in P_{w_l}$. Since $p, q \neq \bm{\epsilon}$, there exists $k$ with $p_k \neq \bm{\epsilon}$ and $l$ with $q_l \neq \bm{\epsilon}$. If now $r_0 \dotso r_N$ is a reduced positive word representing $pq$, then $P_v \cap P_v^{-1} = \gekl{\bm{\epsilon}}$ for all $v \in V$ and Lemma~\ref{lem:normal} imply that there exists $1 \leq M \leq N$ such that $r_M \neq \bm{\epsilon}$. Now it follows from Lemma~\ref{lem:normal} that $pq \neq \bm{\epsilon}$, as desired.
\eproof

The following notion will be useful when dealing with products of words.
\bdefin
Given a $\check{T}$-word $W = g_{0,1} \dotsm g_{0,\mu_0} d_1 g_{1,1} \dotsm g_{1,\mu_1} d_2 g_{2,1} \dotsm g_{m-1,\mu_{m-1}} d_m g_{m,1} \dotsm g_{m,\mu_m}$ as in Definition~\ref{def:cT-word}, let 
$
 \fE(W) = e_{0,0} g_{0,1} e_{0,1} \dotsm g_{0,\mu_0} e_{1,0} g_{1,1} \dotsm g_{1,\mu_1} e_{2,0} g_{2,1} \dotsm g_{m-1,\mu_{m-1}} e_{m,0} g_{m,1} \dotsm g_{m,\mu_m} e_{m,\mu_m}
$
be as above. 
\pari

We define $\mfl(W) \defeq \ell(e_{0,1}) + \dotso \ell(e_{0,\mu_0}) + \dotso + \ell(e_{m,0}) + \dotso + \ell(e_{m,\mu_m - 1})$.

We call $W$ properly reduced if the following are satisfied:
\begin{itemize}
\item $W$ is reduced,
\item if $\mu_0 \neq 0$ and $e_{0,1}$ starts with $d \in T$, then $g_{0,1} \notin G_d^{\bar{d}}$,
\item if $\mu_m \neq 0$ and $e_{m,\mu_m-1}$ ends with $e \in T$, then $g_{m,\mu_m} \notin G_e^e$.
\end{itemize}
\edefin
\nopar

Note that if $\fE(W)$ is reduced, then if $\mu_0 \neq 0$, $e_{0,1}$ starts with $d \in T$ and $g_{0,1} \in G_d^{\bar{d}}$, or if $\mu_m \neq 0$, $e_{m,\mu_m-1}$ ends with $e \in T$ and $g_{m,\mu_m} \in G_e^e$, then we must have $m \geq 1$ and $\mfl(W) \geq 1$.
\pars

\blemma
\label{lem:properlyreduced}
Every element of $G$ is represented by a properly reduced $\check{T}$-word.
\elemma
\nopar

\bproof
Every element of $G$ is represented by a $\check{T}$-word 
$$
 W = g_{0,1} \dotsm g_{0,\mu_0} d_1 g_{1,1} \dotsm g_{1,\mu_1} d_2 g_{2,1} \dotsm g_{m-1,\mu_{m-1}} d_m g_{m,1} \dotsm g_{m,\mu_m}
$$
such that  
$$
 \fE(W) = e_{0,0} g_{0,1} e_{0,1} \dotsm g_{0,\mu_0} e_{1,0} g_{1,1} \dotsm g_{1,\mu_1} e_{2,0} g_{2,1} \dotsm g_{m-1,\mu_{m-1}} e_{m,0} g_{m,1} \dotsm g_{m,\mu_m} e_{m,\mu_m}
$$
is a reduced $\bm{o}$-word. Now proceed inductively on $\mfl(W)$. In case $\mfl(W) = 0$ there is nothing to do. Now assume that $e_{0,1}$ starts with $d \in T$, say $e_{0,1} = d d'$, and that $g_0 \in G_d^{\bar{d}}$, i.e., $g_0 = a^{\bar{d}}$. Then 
$$
 W \equiv g_{0,1} g_{0,2} \dotsm = a^{\bar{d}} g_{0,2} \dotsm = a^d g_{0,2} \dotsm \eqdef W'.
$$ 
Then we have $\fE(W') = e_{0,0} d a^d d' g_{0,2} \dotsm$ with $\mfl(W') < \mfl(W)$, and we can apply induction hypothesis.
\eproof
\pars

\bcor
\label{cor:GvP=Pv}
In the setting of Lemma~\ref{lem:normal}, we have $G_v \cap P = P_v$ for all $v \in V$.
\ecor
\nopar

\bproof
Take $g \in G_v$ such that $g \in P$. By Lemma~\ref{lem:normal} and Lemma~\ref{lem:properlyreduced}, there exists a properly reduced positive word 
$$
 W = p_{0,1} \dotsm p_{0,\mu_0} d_1 p_{1,1} \dotsm p_{1,\mu_1} d_2 p_{2,1} \dotsm p_{m-1,\mu_{m-1}} d_m p_{m,1} \dotsm p_{m,\mu_m}
$$
representing $g$. Since $g$ lies in $P$, we must have $m=0$. Thus we obtain that $\fE(W) = \mfd_0 p_1 \mfd_1 p_2 \dotso \mfd_{\mu - 1} p_{\mu} \mfd_{\mu}$ is a reduced $\bm{o}$-word, where $\mfd_{\lambda}$ are paths in $T$ and $p_{\lambda} \in P_{v_{\lambda}}$. At the same time, $g$ lies in $G_v$, so we obtain another reduced $\bm{o}$-word representing $g$ of the form $\mff h \bar{\mff}$, for some $u \in V$, $h \in G_u$, $\mff = [\bm{o},u]$. 
\pars

Now we proceed inductively on $\ell([u,v])$. If $\ell([u,v]) = 0$, then we have $h = g$ in $G_v$. Comparing reduced forms, we must have $u = v_1$, and since $P_{v_1} = P_u$ is the positive cone in a totally ordered group, we must have either $h = p_1 a$ or $h a = p_1$ for some $a \in P_{v_1}$. In the first case, we are done. In the second case, we conclude that $\bar{\mff} = a \mfd_1 p_2 \dotso \mfd_{\mu-1} p_{\mu} \mfd_{\mu}$. However, induction on $\ell(\bar{\mff}) = \ell(\mfd_1) + \dotso + \ell(\mfd_{\mu})$ shows that this is impossible unless $a = \beps$. Thus $g = h \in P_v$, as desired. 

If $\ell([u,v]) \geq 1$, assume that $[u,v]$ ends with $e \in E$. Then $g = a^e = a^{\bar{e}}$, and $a^{\bar{e}} \in G_{o(e)}$. Since $\ell([u,o(e)]) < \ell([u,v])$, induction hypothesis implies that $a^{\bar{e}} \in P_{o(e)}$, and thus $a^{\bar{e}} \in P_{\bar{e}}^{\bar{e}} = P_e^e \subseteq P_v$ (using that $e \in T$). This implies that $g \in P_v$, as desired.
\eproof
\pars

\blemma
\label{lem:l>1notinPv}
Let $W$ be a properly reduced positive word with $\mfl(W) \geq 1$. Then $W \notin P_v$ for any $v \in V$.
\elemma
\nopar

\bproof
Let $W = p_{0,1} \dotsm p_{0,\mu_0} d_1 p_{1,1} \dotsm$ be a properly reduced word and $W' \equiv q \in P_v$ for some $v \in V$. If $W = W'$, then we have $m=0$ and $W \equiv p_{0,1} \dotsm p_{0,\mu_0} \in P_T$ by Lemma~\ref{lem:W=W'}~(ii). Write $p_{\lambda} \defeq p_{0,\lambda}$. It follows that $\fE(W) = e_0 p_1 e_1 \dotsm e_{\mu-1} p_{\mu} e_{\mu}$ and $\fE(W') = d_0 q d_1$. Without loss of generality, we may assume that $\fE(W')$ is reduced. Since $\fE(W) = \fE(W')$, Lemma~\ref{lem:W=W'}~(i) implies $\ell(d_0) + \ell(d_1) = \ell(e_0) + \ell(e_{\mu}) + \mfl(W)$. So either $\ell(d_0) > \ell(e_0)$ or $\ell(d_1) > \ell(e_{\mu})$. Suppose that $\ell(d_1) > \ell(e_{\mu})$, the other case is similar. Assume that $e_{\mu-1}$ ends with $e \in T$, say $e_{\mu-1} = e'e$. Then Lemma~\ref{lem:W=W'}~(i) yields $e_0 p_1 e_1 \dotsm e' a = d_0 q d'_1$, for some $a \in G_e^e$, where $d'_1$ is the subpath of $d_1$ of length $\ell(d_1) - \ell(e_{\mu}) - 1$ starting from $v$. Going back to $\chT$-words, we deduce that $p_1 \dotsm p_{\mu-1} a = q$ which implies $p_{\mu} = a \in G_e^e$. But this contradicts $p_{\mu} \notin G_e^e$.
\eproof
\pars

The following is a straightforward consequence of Lemma~\ref{lem:l>1notinPv}.
\blemma
\label{lem:W=W'implies...}
Let $W$ and $W'$ be properly reduced positive words with $W, W' \in P_T$. If $W = W'$, then $\mfl(W) = \mfl(W')$, and if $W = p_1 \dotso p_{\mu}$ with $p_{\lambda} \in P_{u_{\lambda}}$ and $W' = p'_1 \dotso p'_{\mu'}$ with $p'_{\lambda'} \in P_{u'_{\lambda'}}$, then $u_1 = u'_1$ and $u_{\mu} = u'_{\mu'}$.
\elemma

\blemma
\label{lem:WW'}
Suppose that $W = p_{0,1} \dotsm p_{0,\mu_0} d_1 \dotsm d_m p_{m,1} \dotsm p_{m,\mu_m}$ and $W' = q_{0,1} \dotsm q_{0,\nu_0} e_1 \dotsm e_n q_{n,1} \dotsm q_{n,\nu_n}$ are properly reduced positive words with 
\begin{eqnarray*}
 \fE(W) &=& \bm{d}_{0,0} p_{0,1} \dotsm p_{0,\mu_0} \bm{d}_{1,0} \dotsm \bm{d}_{m,0} p_{m,1} \dotsm p_{m,\mu_m} \bm{d}_{m,\mu_m} \\
 \fE(W') &=& \bm{e}_{0,0} q_{0,1} \dotsm q_{0,\nu_0} \bm{e}_{1,0} \dotsm \bm{e}_{n,0} q_{n,1} \dotsm q_{n,\nu_n} \bm{e}_{n,\nu_n}.
\end{eqnarray*}
Then $WW'$ is a reduced positive word unless $p_{m,\mu_m} \in P_{u_{m,\mu_m}}$, $q_{0,1} \in P_{v_{0,1}}$, $u_{m,\mu_m} = v_{0,1}$, $\bm{d}_{m,\mu_m-1}$ ends with $\mfd \in T$, $\bm{e}_{0,1}$ starts with $\mfe \in T$, $\mfe = \bar{\mfd}$ and $p_{m,\mu_m} q_{0,1} \in P_{\mfd}^{\mfd}$.
\elemma
\nopar

\bproof
If $u_{m,\mu_m} \neq v_{0,1}$, then 
$$
 \fE(WW') = \bm{d}_{0,0} p_{0,1} \dotsm p_{0,\mu_0} \bm{d}_{1,0} \dotsm \bm{d}_{m,0} p_{m,1} \dotsm p_{m,\mu_m} [u_{m,\mu_m},v_{0,1}] q_{0,1} \dotsm q_{0,\nu_0} \bm{e}_{1,0} \dotsm \bm{e}_{n,0} q_{n,1} \dotsm q_{n,\nu_n} \bm{e}_{n,\nu_n}
$$
is reduced because $\ell([u_{m,\mu_m},v_{0,1}]) \geq 1$. 
\pars

If $u_{m,\mu_m} = v_{0,1}$, then
$$
 \fE(WW') = \bm{d}_{0,0} p_{0,1} \dotsm p_{0,\mu_0} \bm{d}_{1,0} \dotsm \bm{d}_{m,0} p_{m,1} \dotsm \bm{d}_{m,\mu_m-1} (p_{m,\mu_m} q_{0,1}) \bm{e}_{0,1} \dotsm q_{0,\nu_0} \bm{e}_{1,0} \dotsm \bm{e}_{n,0} q_{n,1} \dotsm q_{n,\nu_n} \bm{e}_{n,\nu_n}
$$
If $\bm{d}_{m,\mu_m-1}$ ends with an edge in $A$ or $\bm{e}_{0,1}$ starts with an edge in $A$, then $WW'$ is reduced. Now assume that $\bm{d}_{m,\mu_m-1}$ ends with $\mfd \in T$ and $\bm{e}_{0,1}$ starts with $\mfe \in T$. If $\mfe \neq \bar{\mfd}$, then $WW'$ is again reduced. If $\mfe = \bar{\mfd}$ and $p_{m,\mu_m} q_{0,1} \notin P_{\mfd}^{\mfd}$, then $WW'$ is again reduced.
\eproof
\pars

\bremark
Note that in the situation of Lemma~\ref{lem:WW'}, if $\mfl(W), \mfl(W') \geq 1$, then we actually obtain that $WW'$ is properly reduced, not only reduced. This however does not need to be the case if $\mfl(W) = 0$ or $\mfl(W') = 0$.
\eremark

\blemma
\label{lem:XY_m}
\begin{enumerate}
\item[(i)] Let $X = x_1 \dotsm x_{\mu_0} d_1 x_{\mu_0+1} \dotsm x_{\mu_0 + \mu_1} d_2 \dotsm d_m x_{\mu_0 + \dotso + \mu_{m-1} + 1} \dotsm x_{\mu_0 + \dotso + \mu_m}$ be a positive word and $Y$ another positive word. Let $\nu \defeq \mu_0 + \dotso + \mu_m$. Then there exists a properly reduced positive word $W = p_1 \dotsm p_{\nu_0} f_1 p_{\nu_0+1} \dotsm p_{\nu_0 + \nu_1} f_2 \dotsm f_n p_{\nu_0 + \dotso + \nu_{n-1} + 1} \dotsm p_{\nu_0 + \dotso + \nu_n}$ representing $XY$ such that, if $W_{\nu} \defeq p_1 \dotsm p_{\nu_0} f_1 p_{\nu_0+1} \dotsm p_{\nu_0 + \nu_1} f_2 \dotsm f_l p_{\nu_0 + \dotso + \nu_{l-1} + 1} \dotsm p_{\nu}$ if $\nu_0 + \dotso + \nu_{l-1} < \nu \leq \nu_0 + \dotso + \nu_l$ for some $l$ and $W_{\nu} \defeq W$ otherwise, then $W_{\nu} \in XP$. 
\item[(ii)] Let $x, y \in P$ be represented by positive reduced positive words $X, Y$. Set $\ell \defeq \ell(X)$. Then there exists a reduced positive $\bm{o}$-word $q_0 f_1 q_1 \dotso f_n q_n$ such that $q_0 f_1 q_1 \dotso q_{\ell-1} f_{\ell} q_{\ell} \in xP$.
\end{enumerate}
\elemma
\nopar

\bproof
(i) We proceed inductively on $\mfl(Y)$. We may assume that 
$$
 X = x_1 \dotsm x_{\mu_0} d_1 x_{\mu_0+1} \dotsm x_{\mu_0 + \mu_1} d_2 \dotsm d_m x_{\mu_0 + \dotso + \mu_{m-1} + 1} \dotsm x_{\mu_0 + \dotso + \mu_m}
$$
and $Y = y_1 y_2 \dotsm$ are properly reduced. If $XY$ is reduced and $\mfl(X), \mfl(Y) \geq 1$, then $XY$ is properly reduced and we can just take $W = XY$. Now consider the case that $XY$ is reduced but not properly reduced. Then we must have $\mfl(Y) = 0$ or $\mfl(X) = 0$. In the first case, either $X y'_1$ is properly reduced for some $y'_1 \in P_{w'_1}$ with $y'_1 = y_1$ or there exists a properly reduced positive word $W = p_1 \dotsm p_{\nu_0} f_1 p_{\nu_0+1} \dotsm p_{\nu_0 + \nu_1} f_2 \dotsm f_n p_{\nu_0 + \dotso + \nu_{n-1} + 1} \dotsm p_{\nu_0 + \dotso + \nu_n}$ representing $XY$ with $\nu_0 + \dotso + \nu_n \leq \nu$. In both cases, our claim follows. If $\mfl(X) = 0$, either $x'_1 Y$ is properly reduced for some $x'_1 \in P_{v'_1}$ with $x'_1 = x_1 = X$ or we can write $XY = (x'_1 y_1) Y'$ for some $x'_1 \in P_{v'_1}$ with $x'_1 = x_1 = X$ and $y_1 \in P_{w_1}$ with $w_1 = v'_1$, where $Y' = y_2 \dotso$ satisfies $\mfl(Y') < \mfl(Y)$, and then we can proceed inductively on $\mfl(Y)$. Again, in both cases, our claim follows.
\pari

If $XY$ is not reduced, then Lemma~\ref{lem:WW'} implies that $x_{\nu} y_1 \in P_{\mfd}^{\mfd}$, where $\mfd$ is as in Lemma~\ref{lem:WW'}. Then define $X' = X y_1$ and $Y' = y_2 \dotsm$. Then $\mfl(Y') < \mfl(Y)$, so induction hypothesis produces a properly reduced positive word $W$ representing $X'Y'$ such that $W_{\nu} \in X'P \subseteq XP$.
\pars

(ii) We may assume that $X$ and $Y$ are properly reduced. If $XY$ is reduced, then we can simply take $Z = \fE(XY)$. We can write $\fE(X) = Z'_x Z''_x$, $\fE(Y) = Z'_y Z''_y$ such that $Z = Z'_x Z''_y$ and $\ell(Z'_x) \leq \ell$.
\pari

If $XY$ is not reduced, then we proceed inductively on $\mfl(Y)$. Lemma~\ref{lem:WW'} implies that $xy = (X y_1) Y'$ for some properly reduced positive word $Y'$ with $\mfl(Y') < \mfl(Y)$. Hence it suffices to treat the case where $\mfl(Y) = 0$, say $Y = y_1$ and $X y_1$ is not reduced. But then $X = X' x_{\nu}$ and $XY = X' (x_{\nu} y_1)$ for some positive word $X'$, and we see that $\ell(XY) \leq \ell$. Thus we can take any reduced positive $\bm{o}$-word representing $xy$.
\eproof
\pars

\section{The right LCM property}

Assume that $P$ is a monoid as in \S~\ref{ss:presentation}. Our goal is to establish a criterion when $P$ is right LCM, i.e., for all $p, q \in P$, either $pP \cap qP = \emptyset$ or $pP \cap qP = rP$ for some $r \in P$. It is convenient to introduce the following notation: For $p$, $q \in P$ we write $p \prec q$ if $q \in pP$. Given $p, q \in P$, we denote by $p \vee q$ the (necessarily unique by Corollary~\ref{cor:NoUnits}) minimal element $p \vee q \in P$ satisfying $p, q \prec p \vee q$ if such a minimal element exists. Here minimality refers to $\prec$. In this language, $P$ is right LCM if for all $p, q \in P$, either $pP \cap qP = \emptyset$ or $p \vee q$ exists.

Given $e \in T \amalg A$ and $p \in P$, we set
$$
 p^{-1}(P_{\bar{e}}^{\bar{e}}) \defeq \menge{x \in P_{o(e)}}{px \in P_{\bar{e}}^{\bar{e}}}.
$$

\bdefin
\label{def:LCM}
We say that condition (LCM) is satisfied if for all $e \in T \amalg A$ and $p \in P_{o(e)}$, either $p^{-1}(P_{\bar{e}}^{\bar{e}}) = \emptyset$ or $p^{-1}(P_{\bar{e}}^{\bar{e}}) = q P_{\bar{e}}^{\bar{e}}$ for some $q \in P_{o(e)}$. In the latter case, we set $p^{-1,e} \defeq q$.
\edefin

The main result of this section reads as follows.
\bprop
\label{prop:LCM}
If $P$ is as in \S~\ref{ss:presentation}, then $P$ is right LCM if condition (LCM) is satisfied.
\eprop

In the following, we assume that we are always in the setting of \S~\ref{ss:presentation}. 

We start with the following:
\blemma
\label{lem:p-1Pee}
If condition (LCM) is satisfied, then for all $p \in P$ and $e \in T \amalg A$, either $p^{-1}(P_{\bar{e}}^{\bar{e}}) = \emptyset$ or $p^{-1}(P_{\bar{e}}^{\bar{e}}) = q P_{\bar{e}}^{\bar{e}}$ for some $q \in P$.
\elemma
\nopar

\bproof
Note that if $p \notin P_T$, then $p^{-1}(P_{\bar{e}}^{\bar{e}}) = \emptyset$. Moreover, for all $p \in P$, $e \in T \amalg A$ and $x \in P$, $px \in P_{\bar{e}}^{\bar{e}}$ implies that $x \in P_T$. So we may work in $P_T$.
\pars

We first consider the case that $p \in P_v$ for some $v \in V$. Let $d_1 \dotso d_l = [v,o(e)]$ and set $d_{l+1} \defeq e$. Define $p_0 \defeq p$, $q_1 \defeq p^{-1,d_1}$ if $p^{-1} (P_{d_1}^{\bar{d}_1}) \neq \emptyset$, and, for all $1 \leq i \leq l$, $p_i \defeq p q_1 \dotso q_i$, $q_{i+1} \defeq p_i^{-1,d_{i+1}}$ if $p_i^{-1} (P_{d_{i+1}}^{\bar{d}_{i+1}}) \neq \emptyset$. We claim that $p^{-1}(P_{\bar{e}}^{\bar{e}}) \neq \emptyset$ if and only if $p_i^{-1} (P_{d_{i+1}}^{\bar{d}_{i+1}}) \neq \emptyset$ for all $0 \leq i \leq l$, and that $p^{-1}(P_{\bar{e}}^{\bar{e}}) = q_1 \dotso q_{l+1} P_{\bar{e}}^{\bar{e}}$ in that case, i.e., $p^{-1,e} = q_1 \dotso q_{l+1}$.

To see that, observe that it is easy to see that $p^{-1}(P_{\bar{e}}^{\bar{e}}) \neq \emptyset$ if $p_i^{-1} (P_{d_{i+1}}^{\bar{d}_{i+1}}) \neq \emptyset$ for all $0 \leq i \leq l$. We now prove the converse and $p^{-1,e} = q_1 \dotso q_{l+1}$ inductively on $\ell \defeq \ell([v,o(e)])$. The case $\ell = 0$ follows from condition (LCM) together with Corollary~\ref{cor:GvP=Pv}. Now assume that $\ell \geq 1$. 
%We may assume that $p$ is properly reduced. 
Suppose that $p^{-1}(P_{\bar{e}}^{\bar{e}}) \neq \emptyset$ and take $x \in P$ with $px \in P_{\bar{e}}^{\bar{e}}$.

Using Lemma~\ref{lem:WW'}, we can find positive words $W_m$, $X_m$, $Y_m$ and $f_m \in T$ for $1 \leq m \leq n$ such that $W_1 \equiv p$, $x = X_1 Y_1$, $W_m = W_{m-1} X_{m-1}$, $Y_{m-1} = X_m Y_m$, $W_m \in P_{v_m}$, $X_m \in P_{v_m}$, $W_m X_m \in P_{f_m}^{f_m} \subseteq P_{v_{m+1}}$, $o(f_m) = v_m$, $t(f_m) = v_{m+1}$, and $W_n \equiv px \in P_{o(e)}$. Note that we allow the possibility that $X_m = \emptyset$ or $Y_m = \emptyset$.

Now let $M \in \gekl{1, \dotsc, n}$ be maximal such that $v_M = v$. Then we must have $f_M = d_1$ as $[v,o(e)]$ starts with $d_1$. By construction, we have $W_m = W_1 X_1 \dotso X_{m-1}$, $x = X_1 \dotso X_m Y_m$ and $W_m X_m Y_m = px$. Therefore, if we set $x' \defeq X_1 \dotso X_M$, $x'' \defeq Y_M$, then we have $x = x' x''$ and $px' = W_M X_M$. $W_M \in P_v$ and $X_M \in P_v$ imply that $px' \in P_v$. As $p \in P_v$, this implies that $x' \in G_v \cap P = P_v$ by Corollary~\ref{cor:GvP=Pv}. As $px' \in P_{d_1}^{d_1} = P_{\bar{d}_1}^{d_1}$, this implies that $p^{-1} (P_{d_1}^{\bar{d}_1}) \neq \emptyset$, so that condition (LCM) implies that $x' = p^{-1,d_1} y = q_1 y$ for some $y \in P_{d_1}^{\bar{d}_1}$. Thus $px = p x' x'' = (p q_1) y x'' \in P_{\bar{e}}^{\bar{e}}$, and $p q_1 \in P_{d_1}^{\bar{d}_1} = P_{d_1}^{d_1} \subseteq P_{t(d_1)}$. Now $\ell([t(d_1),o(e)]) < \ell([v,o(e)])$. Hence induction hypothesis implies that $p_i^{-1} (P_{d_{i+1}}^{\bar{d}_{i+1}}) \neq \emptyset$ for all $1 \leq i \leq l$ and that $yx'' \in q_2 \dotso q_{l+1} P_{\bar{e}}^{\bar{e}}$. This shows that $x = x'x'' = q_1 y x'' \in q_1 q_2 \dotso q_{l+1} P_{\bar{e}}^{\bar{e}}$, as desired.

Now let $p \in P_T$ be arbitrary. Let $W = p_0 \dotso p_m$ be a properly reduced positive word representing $p$ with $p_k \in P_{v_k}$. Assume that $\fE(W) = \mfd_0 p_0 \mfd_1 p_1 \dotso p_m \mfd_{m+1}$. We proceed inductively on $m$. We have just dealt with the case $m=0$. Take $x \in P$ such that $px \in P_{\bar{e}}^{\bar{e}}$. Let $X = x_0 \dotso x_n$ be a properly reduced positive word representing $x$ with $x_l \in P_{w_l}$. Write $\fE(X) = \mfe_0 x_0 \mfe_1 x_1 \dotso x_n \mfe_{n+1}$. By Lemma~\ref{lem:WW'} and Lemma~\ref{lem:l>1notinPv}, if $\mfl(X) \geq 1$, then $\mfd_m$ must end with $\bar{d} \in T$, $\mfe_1$ must start with $d \in T$, and we must have $w_0 = v_m$ and $p_m x_0 \in P_d^d$. If $\mfl(X) = 0$, then we can still arrange that $\mfd_m$ ends with $\bar{d} \in T$, and that $w_0 = v_m$ and $p_m x_0 \in P_d^d$. In either case, we obtain that $x_0 \in p_m^{-1} P_d^d$, i.e., $x_0 = p_m^{-1,d} x'_0$ for some $x'_0 \in P_d^d$. Then
$$
 px = p_0 \dotso p_m x_0 \dotso x_n = p_0 \dotso p_m p_m^{-1,d} x'_0 \dotso x_n \in P_{\bar{e}}^{\bar{e}}.
$$
Now $\ti{p} \defeq p_0 \dotso p_m p_m^{-1,d}$, and induction hypothesis produces an element $\ti{q} \in P$ such that $\ti{p}^{-1} P_{\bar{e}}^{\bar{e}} = \ti{q} P_{\bar{e}}^{\bar{e}}$. Hence it follows that $x_0' \dotso x_n \in \ti{p}^{-1} P_{\bar{e}}^{\bar{e}} = \ti{q} P_{\bar{e}}^{\bar{e}}$, and therefore $x = p_m^{-1,d} x'_0 \dotso x_n \in p_m^{-1,d} \ti{q} P_{\bar{e}}^{\bar{e}}$. It is now easy to see that $q \defeq p_m^{-1,d} \ti{q}$ has the desired property.
\eproof

We extend the notation introduced in Definition~\ref{def:LCM} as follows:
\bdefin
We write $p^{-1,e}$ for the element $q$ in Lemma~\ref{lem:p-1Pee} if $p^{-1}(P_{\bar{e}}^{\bar{e}}) \neq \emptyset$.
\edefin
In other words, we define $p^{-1,e}$ as the unique element such that $p^{-1}(P_{\bar{e}}^{\bar{e}}) = p^{-1,e} P_{\bar{e}}^{\bar{e}}$ whenever $p^{-1}(P_{\bar{e}}^{\bar{e}}) \neq \emptyset$.

\blemma
\label{lem:px=q}
Let $p \in P_v$, $x \in P$ such that $px$ is represented by a properly reduced positive word of the form $q_0 q_1 \dotso$ with $q_0 \in P_w$. Let $\ell([v,w]) \geq 1$ such that $[v,w]$ ends with $f \in T$. Then $x \in p^{-1,f} P$.
\elemma
\nopar

\bproof
As in the proof of Lemma~\ref{lem:p-1Pee}, we use Lemma~\ref{lem:WW'} to find positive words $W_m$, $X_m$, $Y_m$ and $f_m \in T$ for $1 \leq m \leq n$ such that $W_1 \equiv p$, $x = X_1 Y_1$, $W_m = W_{m-1} X_{m-1}$, $Y_{m-1} = X_m Y_m$, $W_m \in P_{v_m}$, $X_m \in P_{v_m}$, $W_m X_m \in P_{f_m}^{f_m} \subseteq P_{v_{m+1}}$, $o(f_m) = v_m$, $t(f_m) = v_{m+1}$, and $W_n Y_n$ is a properly reduced positive word representing $px$. Note that we allow the possibility that $X_m = \emptyset$ or $Y_m = \emptyset$. We must have $v_n = w$ by Lemma~\ref{lem:W=W'implies...}.
\pars

Let $M$ be minimal such that $v_M = w$. Then we must have $f_{M-1} = f$. A similar argument as in the proof of Lemma~\ref{lem:p-1Pee} shows that $p X_1 \dotso X_{M-1} = W_{M-1} X_{M-1} \in P_f^{\bar{f}} = P_f^f$, so that $X_1 \dotso X_{M-1} \in p^{-1,f} P_f^{\bar{f}}$. This implies that $x = X_1 \dotso X_{M-1} Y_{M-1} \in p^{-1,f} P$, as desired.
\eproof
\pars

Looking at the way $p^{-1,f}$ has been constructed in the proof of Lemma~\ref{lem:p-1Pee}, the following is an immediate consequence.
\bcor
\label{cor:px=q}
In the situation of Lemma~\ref{lem:px=q}, assume that $[v,w]$ starts with $d \in T$. Then $x \in p^{-1,d} P$.
\ecor

Let $\prec_T$ and $\vee_T$ be the analogues of $\prec$ and $\vee$ with $P_T$ in place of $P$.
\bprop
\label{prop:LCMP,PT}
Assume that condition (LCM) is satisfied.
\pari

Given $p, q \in P_T$, $p P_T \cap q P_T = \emptyset$ if and only if $p P \cap q P = \emptyset$, and $p \vee_T q$ exists if and only if $p \vee q$ exists, and in the latter case, we have $p \vee_T q = p \vee q$. 

Moreover, $P$ is right LCM if and only if $P_T$ is right LCM.
\eprop
\nopar

\bproof
Given $p, q \in P_T$, it is clear that $p P_T \cap q P_T \neq \emptyset$ implies that $p P \cap q P \neq \emptyset$. Now assume that $p P \cap q P \neq \emptyset$, i.e., we can find $x, y \in P$ with $px = qy$. Let $p_0 d_1 p_1 \dotso$ and $q_0 e_1 q_1 \dotso$ be positive words in compact form representing $x$ and $y$. Then we obtain $p p_0 d_1 p_1 \dotso = q q_0 e_1 q_1 \dotso$. Lemma~\ref{lem:W=W'}~(ii) implies that $p p_0 a = q q_0$ or $p p_0 = q q_0 a$ for some $a \in P_T$. This shows that $p P_T \cap q P_T \neq \emptyset$. This also shows that if $p \vee q$ exists, then $p \vee q \prec p p_0 a$ or $p \vee q \prec p p_0$, and in both cases, we obtain $p \vee q \in P_T$. Furthermore, if $px = qy$ holds in $P_T$, then $px \in (p \vee q) P$ implies that $px \in (p \vee q) P_T$ since $px \in P_T$. This shows that $p \vee q = p \vee_T q$. Now if $p \vee_T q$ exists in $P_T$, and if we have $px = qy$ for some $x, y \in P$, then we have $p \vee_T q \prec q q_0$ or $p \vee_T q \prec p p_0$. In the first case, we obtain $p \vee_T q \prec qy = px$ and in the second case, we obtain $p \vee_T q \prec px = qy$. This shows that $p \vee_T q = p \vee q$.
\pars

Now we turn to the second statement. We have already shown \an{$\Rightarrow$}, so let us show \an{$\Leftarrow$}. Let $p, q \in P$ and $p_0 d_1 p_1 \dotso d_m p_m$, $q_0 e_1 q_1 \dotso e_n q_n$ be positive words in compact form representing $p$, $q$. Without loss of generality, assume that $n \geq m$. If $pP \cap qP \neq \emptyset$, i.e., there exist $x, y \in P$ with $px = qy$, then we obtain a positive word in compact form representing $px$ of the form $p_0 d_1 p_1 \dotso p_{m-1} d_m \dotso$, and a positive word in compact form representing $qy$ of the form $q_0 e_1 q_1 \dotso q_{m-1} e_m \dotso$. $px = qy$ implies that $d_l = e_l$ for all $1 \leq l \leq m$. Moreover, Lemma~\ref{lem:W=W'}~(ii) implies that $p_0 d_1 p_1 \dotso p_{m-1} d_m = q_0 d_1 q_1 \dotso q_{m-1} d_m a$ or $p_0 d_1 p_1 \dotso p_{m-1} d_m a = q_0 d_1 q_1 \dotso q_{m-1} d_m$ for some $a \in P_T$. In the first case, we have that
\begin{align*}
 p \vee q &= (q_0 d_1 q_1 \dotso q_{m-1} d_m a p_m) \vee (q_0 d_1 q_1 \dotso q_{m-1} d_m q_m e_{m+1} q_{m+1} \dotso)\\
 &= q_0 d_1 q_1 \dotso q_{m-1} d_m ((a p_m) \vee (q_m e_{m+1} q_{m+1} \dotso))
\end{align*}
exists if and only if $(a p_m) \vee (q_m e_{m+1} q_{m+1} \dotso)$ exists. In the second case, we deduce that 
\begin{align*}
 p \vee q &= (p_0 d_1 p_1 \dotso p_{m-1} d_m p_m) \vee (p_0 d_1 p_1 \dotso p_{m-1} d_m a q_m e_{m+1} q_{m+1} \dotso)\\
 &= p_0 d_1 p_1 \dotso p_{m-1} d_m (p_m \vee (a q_m e_{m+1} q_{m+1} \dotso))
\end{align*}
exists if and only if $p_m \vee (a q_m e_{m+1} q_{m+1} \dotso)$ exists. In both cases, we see that we may assume that $p \in P_T$.

Now let $p \in P_T$, and recall that $q_0 e_1 q_1 \dotso e_n q_n$ is a positive word in compact form representing $q$. We proceed inductively on $n$ to show that $p \vee q$ exists. The case $n=0$ is our assumption that $P_T$ is right LCM. Now assume $n \geq 1$. Since $P_T$ is right LCM, either $p P \cap q_0 P = \emptyset$ or $p \vee q_0$ exists. In the first case, we deduce $p P \cap q P = \emptyset$. In the second case, write $p \vee q_0 = q_0 r$ for some $r \in P_T$. Then $p \vee q = p \vee q_0 \vee q = q_0 r \vee q = q_0 (r \vee e_1 q_1 \dotso e_n q_n)$ exists if and only if $r \vee e_1 q_1 \dotso e_n q_n$ exists. To show the latter, take $x \in P$ such that $r x \in e_1 P$. Then similar arguments involving Lemma~\ref{lem:WW'} as in the proof of Lemma~\ref{lem:p-1Pee} imply that we must have $x = x' x''$ such that $r x' \in P_{e_1}^{\bar{e}_1}$ and $x'' \in e_1 P$. Now Lemma~\ref{lem:p-1Pee} implies that $x' \in r^{-1,e_1} P_{e_1}^{\bar{e}_1}$, say $x' = r^{-1,e_1} y^{\bar{e}_1}$. Let $r r^{-1,e_1} = a^{\bar{e}_1}$.

If $e_1 \in A_-$, then $e_1 = a^{\bar{e}_1} e_1 (a^{e_1})^{-1} = r r^{-1,e_1} e_1 (a^{e_1})^{-1} \in rP$. This implies that $e_1 q_1 \dotso e_n q_n \subseteq e_1 P \subseteq rP$ and hence $r \vee (e_1 q_1 \dotso e_n q_n) = e_1 q_1 \dotso e_n q_n$.

If $e_1 \in A_+$, then $a^{\bar{e}_1} e_1 = e_1 a^{e_1}$, $y^{\bar{e}_1} e_1 = e_1 y^{e_1}$, so that $rx = r r^{-1,e_1} y^{\bar{e}_1} e_1 \dotso = e_1 a^{e_1} y^{e_1} \dotso \in e_1 a^{e_1} P$ implies $r \vee e_1 = e_1 a^{e_1} P$. Thus $r \vee e_1 q_1 \dotso e_n q_n = (r \vee e_1) \vee e_1 q_1 \dotso e_n q_n = e_1 a^{e_1} \vee e_1 q_1 \dotso e_n q_n = e_1 (a^{e_1} \vee q_1 \dotso e_n q_n)$ and $a^{e_1} \vee q_1 \dotso e_n q_n$ exists by induction hypothesis.
\eproof
\pars

We now prove the following proposition, which together with Proposition~\ref{prop:LCMP,PT} implies Proposition~\ref{prop:LCM}.
\bprop
\label{prop:LCMPT}
If condition (LCM) is satisfied, then $P_T$ is right LCM.
\eprop
\nopar

\bproof
We start with $p \in P_v$, $q \in P_w$ and show inductively on $\ell([v,w])$ that $p \vee q$ exists. Without loss of generality we can assume that $p$ and $q$ are reduced when viewed as positive words. The case $\ell([v,w]) = 0$ is clear, so let us consider the case when $v \neq w$.
\pars

Suppose that $x, y \in P_T$ satisfy $px = qy$. Using Lemma~\ref{lem:WW'} as before, we can find positive words $W_m$, $X_m$, $Y_m$ and $f_m \in T$ for $1 \leq m \leq n$ such that $W_1 \equiv p$, $x = X_1 Y_1$, $W_m = W_{m-1} X_{m-1}$, $Y_{m-1} = X_m Y_m$, $W_m \in P_{v_m}$, $X_m \in P_{v_m}$, $W_m X_m \in P_{f_m}^{f_m} \subseteq P_{v_{m+1}}$, $o(f_m) = v_m$, $t(f_m) = v_{m+1}$, and $W_n Y_n$ is a properly reduced positive word representing $px$. Note that we allow the possibility that $X_m = \emptyset$ or $Y_m = \emptyset$. Again using Lemma~\ref{lem:WW'} as before, we can find positive words $W'_{m'}$, $X'_{m'}$, $Y'_{m'}$ and $f'_{m'} \in T$ for $1 \leq m' \leq n'$ such that $W'_1 \equiv q$, $y = X'_1 Y'_1$, $W'_{m'} = W'_{m'-1} X'_{m'-1}$, $Y'_{m'-1} = X'_{m'} Y'_{m'}$, $W'_{m'} \in P_{v'_{m'}}$, $X'_{m'} \in P_{v'_{m'}}$, $W'_{m'} X'_{m'} \in P_{f'_{m'}}^{f'_{m'}} \subseteq P_{v'_{m'+1}}$, $o(f'_{m'}) = v'_{m'}$, $t(f'_{m'}) = v'_{m'+1}$, and $W'_{n'} Y'_{n'}$ is a properly reduced positive word representing $qy$. As before, we allow the possibility that $X'_{m'} = \emptyset$ or $Y'_{m'} = \emptyset$. We must have $v_n = v'_{n'}$ by Lemma~\ref{lem:W=W'implies...}. Assume that the paths $v_1, v_2, \dotsc$ and $v'_1, v'_2, \dotsc$ meet for the first time at $u \in V$. We must have $u \in [v,w]$. So we have $x = x' x''$ and $y = y' y''$ such that $px', qy' \in P_u$. Now use that $P_u$ is the positive cone in the totally ordered group $G_u$. We obtain that $px'z = qy'$ or $px' = qy'z$ for some $z \in P_u$. In the first case, observe that $q y' y'' = p x' z y'' = p x' x''$ implies $z y'' = x''$ and thus $x = (x'z)y''$, $y = y' y''$ and $p (x'z) = q y'$. In the second case, observe that $p x' x'' = q y' z x'' = q y' y''$ implies $y'' = z x''$ and thus $y = (y'z)x''$, $x = x' x''$ and $p x' = q (y'z)$. So we may assume that $px' = qy'$.

a) Suppose that for all $x, y \in P$ with $px = qy$, the vertex $u$ as above satisfies $u \in [v,w] \setminus \gekl{v,w}$. Let $[v,w]$ start with $d$ and end with $f$. Then it follows as in the proof of Lemma~\ref{lem:p-1Pee} that $x \in p^{-1,d}P$ and $y \in q^{-1,\bar{f}}P$. Therefore $p \vee q = p p^{-1,d} \vee q q^{-1,\bar{f}}$. Now $p p^{-1,d} \in P_d^{\bar{d}} = P_d^d \subseteq P_{t(d)}$ and $q q^{-1,\bar{f}} \in P_{\bar{f}}^f = P_{\bar{f}}^{\bar{f}} \subseteq P_{o(f)}$, and $\ell([t(d),o(f)]) < \ell([v,w])$. Hence induction hypothesis implies that $p p^{-1,d} \vee q q^{-1,\bar{f}}$ exists.

b) Suppose that there exist $x, y \in P$ with $px = qy$ such that the vertex $u$ as above satisfies $u = v$. Then $px' = qy' \in P_u$, where $x'$ and $y'$ are as before. We must have $y' \in q^{-1,\bar{d}} P$, i.e., $y' = q^{-1,\bar{d}} \ti{y}$, by the same argument as above (where $d$ is as in a)). 

b.1) If $q q^{-1,\bar{d}} \prec p$ in $P_u$, then $p \in q q^{-1,\bar{d}} P \subseteq qP$, so that $pP \subseteq qP$ and hence $p \vee q = p$.

b.2) If $p \prec q q^{-1,\bar{d}}$ in $P_u$, i.e., $q q^{-1,\bar{d}} = p z$, then $px' = qy' = q q^{-1,\bar{d}} \ti{y} = p z \ti{y}$. This implies $x' = z \ti{y}$ and thus $x' \in zP$. As $pz = q q^{-1,\bar{d}} \in P_{\bar{d}}^d = P_d^{\bar{d}}$, we obtain that $z \in p^{-1,d}P$. Hence it follows that $px \in px'P \subseteq p p^{-1,d}P$. Moreover, we obtain as in a) that $qy \in q q^{-1,\bar{f}} P$ (where $f$ is as in a)). Thus $p \vee q = p p^{-1,d} \vee q q^{-1,\bar{f}}$ exists by induction hypothesis.

The case $u = w$ is analogous to b).

Now consider $p, q \in P_T$ arbitrary. Let $W_p = p_0 p_1 \dotso p_m$ and $W_q = q_0 q_1 \dotso q_n$ be properly reduced positive words representing $p$ and $q$, with $p_k \in P_{v_k}$ and $q_l \in P_{w_l}$. We proceed inductively on $\mfl(W_p) + \mfl(W_q)$. 

Suppose that $px = qy$ for some $x, y \in P$. First consider the case that $m, n \geq 1$. If there are properly reduced positive words $W_x$ and $W_y$ representing $x$ and $y$ such that $W_p W_x$ and $W_q W_y$ are properly reduced, then Lemma~\ref{lem:W=W'implies...} implies that $p_0 = q_0 a$ or $p_0 a = q_0$ because our semigroups are positive cones in totally ordered groups. If $p_0 = q_0 a$, then $p = p_0 p_1 \dotso p_m = q_0 a p_1 \dotso p_m$, and $p \vee q = (q_0 a p_1 \dotso p_m) \vee (q_0 q_1 \dotso q_n) = q_0 ((a p_1 \dotso p_m) \vee (q_1 \dotso q_n))$ exists if and only if $(a p_1 \dotso p_m) \vee (q_1 \dotso q_n)$ exists. The latter now follows from induction hypothesis as $a p_1 \dotso p_m$ and $q_1 \dotso q_n$ can be represented by properly reduced positive words with smaller $\mfl$. The case $p_0 a = q_0$ is analogous.

It remains to consider the case that for all properly reduced positive words $W_x$ and $W_y$ representing $x$ and $y$, $W_p W_x$ or $W_q W_y$ is not properly reduced. As we proceed inductively on $\mfl(W_p) + \mfl(W_q)$, we may assume that $v_0 \neq w_0$. Write $\fE(W_p) = \mfd_0 p_0 \mfd_1 \dotso \mfd_m p_m \mfd_{m+1}$ and $\fE(W_q) = \mfe_0 q_0 \mfe_1 \dotso \mfe_n q_n \mfe_{n+1}$. Suppose that $\mfd_m$ ends with $d$ and $\mfe_n$ ends with $e$. Given $x, y \in P$ with $px = qy$, we claim that $x \in p_m^{-1,\bar{d}} P$ and $y \in q_n^{-1,\bar{e}} P$ unless $p \in qP$ or $q \in pP$. By Lemma~\ref{lem:WW'}, this is clear if both $W_p W_x$ and $W_q W_y$ are not properly reduced. Let us now consider the case that $W_p W_x$ is not properly reduced, while $W_q W_y$ is properly reduced. (The other case is similar.) Lemma~\ref{lem:WW'} implies that $x \in p_m^{-1,d}P$, say $x = p_m^{-1,d} x'$. We proceed inductively on $\mfl(W_p)$ to show that $y \in q_n^{-1,\bar{e}} P$ unless $q \in pP$. We have $qy = px = p_0 \dotso p_{m-1} (p_m p_m^{-1,d}) x'$, so that induction hypothesis implies $y \in q_n^{-1,e} P$ unless $q \in pP$. Thus it is enough to treat the case when $p \in P_v$. Assume that $[v,w_0]$ ends with $e_0 \in E$. Since $px = qy = W_q W_y$, because $W_q W_y$ is properly reduced and starts with $q_0 \in P_w$, Lemma~\ref{lem:px=q} implies that $x \in p^{-1,e_0}P$. So $x = p^{-1,e_0} x_1$. We have $p p^{-1,e_0} \in P_{w_0}$. If $p p^{-1,e_0} \prec q_0$, then $q_0 \in pP$ and thus $q \in pP$. Otherwise, we have $q_0 \prec p p^{-1,e_0}$, say $p p^{-1,e_0} = q_0 p_1$ for some $p_1 \in P_{e_0}^{\bar{e}_0} = P_{e_0}^{e_0} \subseteq P_{w_0}$. Then $q y = q_0 q_1 \dotso q_n y = px = q_0 p_1 x_1$ and thus $q_1 \dotso q_n y = p_1 x_1$. Let $q_0^{(1)} \dotso q_{n_1}^{(1)} y$ be the properly reduced positive word representing $q_1 \dotso q_n y$ obtained via the algorithm from the proof of Lemma~\ref{lem:properlyreduced}. We have $q_{n_1}^{(1)} \in P_{\bar{e}}^{\bar{e}} q_n$. Again, Lemma~\ref{lem:px=q} implies that $x_1 \in p_1^{-1,e_1} P$, say $x_1 = p_1^{-1,e_1} x_2$, where $e_1 \in E$ lies in $[w_1,w_2] [w_2,w_3] \dotso [w_{n-1},w_n]$. Continuing in this way, we obtain elements $x_{\lambda} \in P$ and $p_{\lambda} \in P_{e_{\lambda-1}}^{\bar{e}_{\lambda-1}}$ such that $x_{\lambda} = p^{-1,e_{\lambda}} x_{\lambda+1}$ and $p_{\lambda} p_{\lambda}^{-1,e_{\lambda}} = q_0^{(\lambda)} p_{\lambda+1}$, where $e_{\lambda} \in E$ lies in $[w_1,w_2] [w_2,w_3] \dotso [w_{n-1},w_n]$ and $q_0^{(\lambda)} \dotso q_{n_{\lambda}}^{(\lambda)} y$ is a properly reduced positive word representing $q_0^{(\lambda-1)} \dotso q_{n_{\lambda - 1}}^{(\lambda - 1)} y$ with $q_{n_{\lambda}}^{(\lambda)} \in P_{\bar{e}}^{\bar{e}} q_n$. We end up with $q_0^{(\nu)} y = p_{\nu} x_{\nu}$. Again, Lemma~\ref{lem:px=q} implies that $x_{\nu} \in p_{\nu}^{-1,e} P$. If $p_{\nu} p_{\nu}^{-1,e} \prec q_0^{(\nu)}$ then $q_0^{(\nu)} \in p_{\nu} P$ and hence $q_0 \dotso q_n = q_0^{(0)} \dotso q_0^{(\nu-1)} q_0^{(\nu)} P \subseteq q_0^{(0)} \dotso q_0^{(\nu-1)} p_{\nu} P \subseteq \dotso \subseteq pP$. Otherwise, we have $q_0^{(\nu)} \prec p_{\nu} p_{\nu}^{-1,e}$, say $p_{\nu} p_{\nu}^{-1,e} = q_0^{(\nu)} p_{\nu+1}$. 

As $p_{\nu} p_{\nu}^{-1,e} \in P_{e}^{\bar{e}} = P_e^e$, we conclude that $p_{\nu+1} \in (q_0^{(\nu)})^{-1,\bar{e}}P$. This yields $q_0^{(\nu)} y = p_{\nu} x_{\nu} \in p_{\nu} p_{\nu}^{-1,e} P = q_0^{(\nu)} p_{\nu+1} P$ and thus $y \in p_{\nu+1} P \subseteq (q_0^{(\nu)})^{-1,\bar{e}}P \subseteq q_n^{-1,\bar{e}}P$, as desired.

It remains to treat the case when $m=0$, say $p \in P_v$, while $n \geq 1$. (The other case is similar.) As we proceed inductively on $\mfl(W_q)$, we may assume that $v \neq w_0$. Consider properly reduced positive words $W_x$ and $W_y$ representing $x$ and $y$. If both $p W_x$ and $W_q W_y$ are properly reduced, then Lemma~\ref{lem:W=W'implies...} implies $v = w_0$. So it suffices to treat the case that $p W_x$ or $W_q W_y$ is not properly reduced. As before, write $\fE(W_q) = \mfe_0 q_0 \mfe_1 \dotso \mfe_n q_n \mfe_{n+1}$. Suppose that $[v,w_0]$ starts with $d$ and $\mfe_n$ ends with $e$. Given $x, y \in P$ with $px = qy$, we claim that $x \in p^{-1,d} P$ and $y \in q_n^{-1,\bar{e}} P$ unless $p \in qP$ or $q \in pP$. We proceed inductively on $\mfl(W_q)$. If $W_q W_y$ is not properly reduced, then Lemma~\ref{lem:WW'} implies that $y \in q_n^{-1,\bar{e}} P$, say $y = q_n^{-1,\bar{e}} y'$. Then $px = qy = (q q_n^{-1,\bar{e}}) y'$, and induction hypothesis implies that $x \in p^{-1,d}P$ because $q q_n^{-1,\bar{e}}$ can be represented by a properly reduced positive word with $\mfl$ strictly less than $\mfl(W_q)$. And the case when $q = q_0$ has already been treated above. If $W_q W_y$ is properly reduced but $p W_x$ is not properly reduced, then Corollary~\ref{cor:px=q} implies that $x \in p^{-1,d} P$. And a similar argument as in the previous case ($m, n \geq 1$) shows that $y \in q_n^{-1,\bar{e}}P$ unless $q \in pP$.

The case $m=0=n$ has already been dealt with at the beginning of the proof.
\eproof
\pars

\bproof[Proof of Proposition~\ref{prop:LCM}]
Proposition~\ref{prop:LCM} follows from Propositions~\ref{prop:LCMP,PT} and \ref{prop:LCMPT}.
\eproof

\section{Semigroup C*-algebras and their groupoid models}
\label{s:SemigroupC}

Let us briefly recall a few facts about semigroup C*-algebras. Given a left cancellative semigroup $P$, its regular representation assigns to every $p \in P$ the isometry $\lambda_p: \: \ell^2 P \to \ell^2 P$ determined by $\lambda_p(\delta_x) = \delta_{px}$ for all $x \in P$, where $\menge{\delta_x}{x \in P}$ is the canonical orthonormal basis of $\ell^2 P$. The (reduced) semigroup C*-algebra $C^*_{\lambda}(P)$ is the C*-algebra generated by $\menge{\lambda_p}{p \in P}$. We refer the reader to \cite{Li12,Li13} as well as \cite[\S~5]{CELY} for more details about general semigroup C*-algebras.

In this paper, we will only consider monoids $P$ with the following properties: $P$ is countable, embeds into a group $G$, has no non-trivial invertible elements, i.e., $P^* = \gekl{\beps}$, and $P$ is right LCM, i.e., for all $p, q \in P$, either $pP \cap qP = \emptyset$ or $pP \cap qP = rP$ for some $r \in P$. In this case, the set of constructible right ideals of $P$ is given by $\cJ = \menge{pP}{p \in P} \cup \gekl{\emptyset}$. $\cJ$ is a semilattice under intersection because $P$ is right LCM. 

Now let $\Omega$ be the set of multiplicative, non-zero maps $\chi: \: \cJ \to \gekl{0,1}$ sending $\emptyset$ to $0$. Equip $\Omega$ with the topology of point-wise convergence. For brevity, we set $\chi(p) \defeq \chi(pP)$ for all $\chi \in \Omega$. Consider the partial action $G \curvearrowright \Omega$ determined by the partial homeomorphisms $U_{g^{-1}} \to U_g, \, \chi \ma g.\chi$, where $U_{g^{-1}}$ is the subspace of all $\chi \in \Omega$ which satisfy $\chi(q) = 1$ for some $q \in P$ such that $g = pq^{-1}$ for some $p \in P$. For such $\chi$, $g.\chi$ is determined by $(g.\chi)(x) = \chi(qy)$ if $xP \cap pP = pyP$ and $(g.\chi)(x) = 0$ if $xP \cap pP = \emptyset$. \cite[Theorem~5.6.41]{CELY} implies that $C^*_{\lambda}(P)$ is canonically isomorphic to the reduced groupoid C*-algebra $C^*_r(G \ltimes \Omega)$ of the transformation groupoid $G \ltimes \Omega$ attached to the partial action $G \curvearrowright \Omega$. 

Let us describe distinguished elements and subspaces of $\Omega$. First of all, every $pP \in \cJ$ determines a point $\chi_p \in \Omega$ given by $\chi_p(x) = 1$ if $pP \subseteq xP$ and $\chi_p(x) = 0$ if $pP \nsubseteq xP$. This allows us to identify $P$ with a subset of $\Omega$ because $P^* = \gekl{e}$. We define $\Omega_{\infty} \defeq \Omega \setminus P$. Among the points in $\Omega_{\infty}$, we single out those $\chi \in \Omega$ for which $\chi^{-1}(1)$ is maximal, i.e., whenever $\omega \in \Omega$ satisfies $\omega(x) = 1$ for all $x \in P$ with $\chi(x) = 1$, then we must have $\omega = \chi$. We set $\Omega_{\max} \defeq \menge{\chi \in \Omega}{\chi^{-1}(1) \ {\rm is} \ {\rm maximal}}$. Note that $\Omega_{\max} \subseteq \Omega_{\infty}$. Moreover, we define $\partial \Omega \defeq \overline{\Omega_{\max}}$. Let us now collect a few facts about $\partial \Omega$, which are obtained in \cite[\S~5.7]{CELY} in greater generality than needed here. $\partial \Omega$ is the minimal non-empty closed $G$-invariant subspace of $\Omega$. Moreover, $\partial \Omega$ reduces to a single point (namely $\chi \in \Omega$ given by $\chi(x) = 1$ for all $x \in P$; we usually denote this $\chi$ by $\infty$) if and only if $P$ is left reversible. If $\partial \Omega$ is not a point, then $G \curvearrowright \partial \Omega$ is purely infinite. Since $\partial \Omega$ is always a closed $G$-invariant subspace of $\Omega$, we can define a quotient of $C^*_{\lambda}(P)$ by setting $\partial C^*_{\lambda}(P) \defeq C^*_r(G \ltimes \partial \Omega)$. $\partial C^*_{\lambda}(P)$ is called the boundary quotient of $C^*_{\lambda}(P)$. We need the following characterization of elements in $\Omega_{\max}$.
\blemma[see {\cite[Lemma~5.7.4]{CELY}}]
\label{lem:Omegamax_genchar}
Let $\chi \in \Omega$. $\chi$ lies in $\Omega_{\max}$ if and only if for any $p \in P$ with $\chi (p)=0$, there exists $q \in P$ such that $\chi(q)=1$ and $pP \cap qP = \emptyset$. 
\elemma
By construction, for every closed subspace $X \subseteq \Omega$, the following sets form a collection of basic open sets in $X$: $X(p, \gekl{q_i}) \defeq \menge{\chi \in X}{\chi(p) = 1, \, \chi(q_i) = 0 \ \forall \ i}$, where $p \in P$ and $\gekl{q_i} \subseteq P$ is a finite set.

Because $P$ is right LCM, elements in $\Omega$ can be described by words in $P$. Let $w$ be a word in $P$, i.e., $w = x_1 x_2 \dotsm$ with $x_{\bullet} \in P$. Set $w_{\text{---}i} \defeq x_1 \dotsm x_i$. Define $\chi_w \in \Omega$ by $\chi_w(p) = 1$ if and only if there exists $i$ such that $w_{\text{---}i} \in pP$. A similar argument as in \cite[\S~2.2]{LOS} shows that every element in $\Omega$ is of the form $\chi_w$ for some word $w$. Moreover, the partial action $G \curvearrowright \Omega$ can be described as follows: Given $g \in G$, and with $w$ as before, $g. \chi_w$ is defined if and only if $g = p w_{\text{---}i}^{-1}$ for some $i$, and then $g.\chi_w = \chi_{w'}$, where $w' = p x_{i+1} x_{i+2} \dotsm$.

Finally, we need the following characterization of topological freeness for $G \curvearrowright \partial \Omega$.
\btheo[{\cite[Proposition~6.18]{LS}} and {\cite[Theorem~A]{Li21a}}]
\label{thm:Gc}
Define $G^c \defeq \menge{g \in G}{(pP) \cap (gpP) \neq \emptyset \quad \forall \ p \in P}$. Then $G \curvearrowright \partial \Omega$ is topologically free if and only if $G^c = \gekl{\beps}$.
\etheo

\section{Closed invariant subspaces}

\label{s:clinvsub}

Throughout this section, assume that $P$ is as in \S~\ref{ss:presentation} and that condition (LCM) is satisfied. In addition, suppose that all our groups $G_v$, $v \in V$, are countable, and that $V$ and $E$ are countable, too. Our goal is to study closed invariant subspaces of $\Omega$, where $\Omega$ is as in \S~\ref{s:SemigroupC}. In our discussion, we first exclude the case of one vertex generalized Baumslag-Solitar groups and their monoids, as in Example~\ref{ex:pi}~(ii); that special case will be discussed separately later on.

\subsection{The general case}

By an infinite positive word, we mean an infinite word in $\gekl{P_v}_{v \in V} \cup A$ of the form $W = x_1 x_2 x_3 \dotso$, where $x_{\bullet} \in \gekl{P_{v} \setminus \gekl{\bm{\epsilon}}}_{v \in V} \cup A$. Let $W_{\text{---}j} \defeq x_1 x_2 \dotso x_j$. Recall that given an infinite positive word $W$, the associated character $\chi_W \in \Omega$ is determined by $\chi_W(p) = 1$ $\Leftrightarrow$ there exists $j$ such that $W_{\text{---}j} \in pP$ (see \S~\ref{s:SemigroupC}). 

We will be interested in the following two situations.
\nopar
\begin{enumerate}
\item[I.] For all $v \in V$, $x \in P_v \setminus \gekl{\bm{\epsilon}}$ or $x \in A$ and $\chi \in \Omega$ there exists a finite or infinite positive word $W$ with $\chi = \chi_W$, a strictly increasing sequence $(j_N)_N$ of positive integers, and a finite positive word $Y$ whose first letter does not lie in $P_v$ in case $x \in P_v$, such that, with $W_N \defeq W_{\text{---}j_N}$, we have
\begin{itemize}
\item $x Y W_N$ is a reduced positive word for all $N$,
\item Whenever $p_1 \dotsm p_{\mu} d_1 \dotsm$ ($\mu \neq 0$) is a properly reduced positive word representing $x Y W_N$, then we must have $x \in p_1 P_T$ if $x \in P_v$ and $x \in p_1 P$ if $x \in A$.
\end{itemize}   
\item[II.] There exists $u \in V$ and $\bm{b} \in P_u$ such that the following holds: For all $v \in V$, $x \in P_v \setminus \gekl{\bm{\epsilon}}$  or $x \in A$ and $\chi \in \Omega$ there exists a finite or infinite positive word $W$ with $\chi = \chi_W$, a strictly increasing sequence $(j_N)_N$ of positive integers, and a positive word $Y$ whose first letter does not lie in $P_v$ in case $x \in P_v$, such that, with $W_N \defeq W_{\text{---}j_N}$, we have
\begin{itemize}
\item $x Y W_N$ is a reduced positive word for all $N$,
\item Whenever $p_1 \dotsm p_{\mu} d_1 \dotsm$ ($\mu \neq 0$) is a properly reduced positive word representing $x Y W_N$, then one of the following holds:
\begin{enumerate}
\item[A)] $x \in p_1 P_T$ if $x \in P_v$ and $x \in p_1 P$ if $x \in A$,
\item[B)] $W_N \in \bm{b}P$ and $x \bm{b}^i \in p_1 P_T$ if $x \in P_v$ and $x \bm{b}^i \in p_1 P$ if $x \in A$, where $i$ is some positive integer.
\end{enumerate}
\end{itemize}
\end{enumerate}
\pars

\blemma
\label{lem:I.gen}
Suppose that condition I. holds. Let $\chi \in \Omega$ be arbitrary. For every $\eta \in \Omega$ such that $\eta = \chi_X$ for some infinite positive word $X$ with $\lim_{l \to \infty} \ell(X_{\text{---}l}) = \infty$, we have $\eta \in \overline{G.\chi}$.
\elemma
\nopar

\bproof
Let $\dotsm f_n x_1 \dotsm x_{\nu}$ be a properly reduced positive word representing $X_{\text{---}l}$. We distinguish between two cases:
If $\nu \neq 0$ and $x_1 \dotsm x_{\nu} \neq \beps$, say $x_{\nu} \in P_v \setminus \gekl{\bm{\epsilon}}$, then let $x = x_{\nu}$, and if $\nu = 0$ or $x_1 \dotsm x_{\nu} = \beps$ and $f_n \in A$, then let $x = f_n$.
Condition I. applied to $\chi$ and $x$ provides $W$, $W_N$ and $Y$ as above. Note that these depend on $l$. We now claim that $\lim_{l \to \infty} \chi_{X_{\text{---}l} Y W} = \eta$.
\pars

If $\eta(p) = 1$, then $X_{\text{---}l} \in pP$ for all sufficiently big $l$, so that $X_{\text{---}l} Y W \in pP$ for all sufficiently big $l$. Thus $\chi_{X_{\text{---}l} Y W}(p) = 1$ for all sufficiently big $l$.

Conversely, suppose that $\chi_{X_{\text{---}l} Y W}(p) = 1$ for all sufficiently big $l$. Then $X_{\text{---}l} Y W \in pP$ for all sufficiently big $l$, say $X_{\text{---}l} Y W = pz$. Let $q_0 e_1 q_1 \dotso q_{M-1} e_M q_M$ be a reduced $\bm{o}$-word representing $pz$. 

For sufficiently big $l$, $X_{\text{---}l}$ can be represented by a reduced $\bm{o}$-word of the form $X_n x \bm{\epsilon} \dotso \bm{\epsilon}$ with $\ell(X_n) > \ell(W_p)$, where $W_p$ is a reduced $\bm{o}$-word representing $p$. Moreover, Lemma~\ref{lem:XY_m}~(ii) applied to $m = \ell(X_n)$ implies that $q_0 e_1 q_1 \dotso q_{m-1} e_m \in pP$, say $q_0 e_1 q_1 \dotso q_{m-1} e_m = p z'$ and $z = z' z''$. Since $Y$ and $W_N$ are as in I., there is a reduced $\bm{o}$-word representing $X_{\text{---}l} Y W$ which starts with $X_n x$. Hence Lemma~\ref{lem:W=W'}~(i) yields that $X_n = q_0 e_1 q_1 \dotso q_{m-1} e_m a$ or $X_n a = q_0 e_1 q_1 \dotso q_{m-1} e_m$. In the first case, we obtain $X_{\text{---}l} \in pP$ and thus $\chi_X(p) = 1$. In the second case, we obtain $pz = pz'z'' = X_n a z''$ and thus $X_n x Y W_N = X_n a z''$, which in turn implies $x Y W_N = a z''$. Lemma~\ref{lem:XY_m}~(i) provides a properly reduced positive word representing $a z''$ starting with $a a' \in P_w$ for some $w \in V$. Now condition I. implies that $x \in a a' P_T \subseteq a P_T$. This in turn yields $X_{\text{---}l} = X_n x \in X_n a P \subseteq pP$ and thus $\chi_X(p) = 1$, as desired.
\eproof
\pars

\blemma
\label{lem:II.gen}
Suppose that condition II. holds. 
\nopar

\begin{enumerate}
\item[(i)] Let $\chi \in \Omega$ satisfy $\chi(\bm{b}) = 0$. For every $\eta \in \Omega$ such that $\eta = \chi_X$ for some infinite positive word $X$ with $\lim_{l \to \infty} \ell(X_{\text{---}l}) = \infty$, we have $\eta \in \overline{G.\chi}$.
\item[(ii)] Let $\chi \in \Omega$ be arbitrary. For every $\eta \in \Omega$ such that $\eta = \chi_X$ for some infinite positive word $X$ with $\lim_{l \to \infty} \ell(X_{\text{---}l}) = \infty$ and $g.\eta(\bm{b}^i) = 1$ for all $g \in G$ for which $g.\eta$ is defined and all positive integers $i$, we have $\eta \in \overline{G.\chi}$.
\end{enumerate}
\elemma

\bproof
Let $\dotsm f_n x_1 \dotsm x_{\nu}$ be a properly reduced positive word representing $X_{\text{---}l}$ as in the proof of Lemma~\ref{lem:I.gen}. Condition II. applied to $\chi$ and $x = x_{\nu}$ (if $\nu \neq 0$ and $x_1 \dotsm x_{\nu} \neq \beps$, say $x_{\nu} \in P_v \setminus \gekl{\bm{\epsilon}}$) and $x = f_n$ (if $\nu = 0$ or $x_1 \dotsm x_{\nu} = \beps$ and $f_n \in A$) provides $W$, $W_N$ and $Y$ as above. Note that these depend on $l$. We now claim that $\lim_{l \to \infty} \chi_{X_{\text{---}l} Y W} = \eta$.
\pars

In (i), B) in II. leads to a contradiction to the assumption that $\chi(\bm{b}) = 0$ because $W_N \in \bm{b}P$ implies $\chi(\bm{b}) = 1$. Hence we must have statement A) in II., and $\lim_{l \to \infty} \chi_{X_{\text{---}l} Y W} = \eta$ follows by the same argument as in the proof of Lemma~\ref{lem:I.gen}.

In (ii), suppose that $\chi_{X_{\text{---}l} Y W}(p) = 1$ for all sufficiently big $l$. We can then use A) in II. and the same argument as in the proof of Lemma~\ref{lem:I.gen} to show $\chi_X(p) = 1$, or we can use B) in II. and the same argument as in the proof of Lemma~\ref{lem:I.gen} to shows that $X_{\text{---}l} \bm{b}^i \in pP$ for some positive integer $i$. Now our assumption that $g.\eta(\bm{b}^i) = 1$ for all $g$ implies for $g = X_{\text{---}l}^{-1}$ that $X_{\text{---}l}^{-1}.\eta(\bm{b}^i) = 1$ and thus $\eta(X_{\text{---}l} \bm{b}^i) = 1$. This, together with $X_{\text{---}l} \bm{b}^i \in pP$, implies that $\eta(p) = 1$.
\eproof
\pars

Suppose that II. holds. Define
$$
 \Omega_{\bm{b},\infty} \defeq \menge{\chi \in \Omega}{(g.\chi)(\bm{b}^i) = 1 \ \forall \ g, \, i},
$$
where $i$ runs through all natural numbers and we only consider those $g \in G$ such that $g.\chi$ is defined. Note that we always have $\Omega_{\bm{b},\infty} \subseteq \Omega_{\infty}$.

To summarize, here is the conclusion.
\blemma
\label{lem:linftyvsPv}
Suppose that I. holds. Then given $\chi \in \Omega$ arbitrary and $\eta \in \Omega$, we have $\eta \in G. \Omega_{P_w}$ for some $w \in V$ or $\eta \in \overline{G. \chi}$.

Suppose that II. holds. 
\nopar

\begin{enumerate}
\item[(i)] Given $\chi \in \Omega$ with $\chi(\bm{b}) = 0$ and $\eta \in \Omega$, we have $\eta \in G. \Omega_{P_w}$ for some $w \in V$ or $\eta \in \overline{G. \chi}$.
\item[(ii)] Given $\chi \in \Omega$ arbitrary and $\eta \in \Omega_{\bm{b},\infty}$, we have $\eta \in G. \Omega_{P_w}$ for some $w \in V$ or $\eta \in \overline{G. \chi}$.
\end{enumerate}
\elemma
Note that $\eta \in G. \Omega_{P_w}$ for some $w \in V$ means that $\eta = \chi_X$ for some infinite positive word $X = x_1 x_2 x_3 \dotso$ such that $x_j \in P_w$ for all sufficiently big $j$.
\bproof
Write $\eta = \chi_X$ for some infinite positive word $X = x_1 x_2 x_3 \dotso$. All we have to show is that if $\sup_j \ell(X_{\text{---}j}) < \infty$, then $\eta \in G. \Omega_{P_w}$ for some $w \in V$. Indeed, $\sup_j \ell(X_{\text{---}j}) < \infty$ implies $\sup_j \mfl(X_{\text{---}j}) < \infty$. Let $\mfl \defeq \liminf_j \mfl(X_{\text{---}j})$. If $\mfl = 0$, then our claim follows from $\sup_j \mfl(X_{\text{---}j}) < \infty$. If $\mfl > 0$, then by passing to a subsequence if necessary, we obtain positive words $X_n$ and $Y_n$ such that $X_n$ is properly reduced, $X_{n+1} \equiv X_n Y_n$, $\mfl(X_n) = \mfl$ for all $n$, and $\eta = \chi_X = \lim_{n \to \infty} \chi_{X_n}$. Now $\mfl(X_{n+1}) = \mfl(X_n)$ for all $n$ implies that there must exist $w \in V$ such that $Y_n \in P_w$ for all $n$. This shows $\eta \in G. \Omega_{P_w}$, as desired.
\eproof
\pars

Now we turn to the following question: When do we have I. or II.?

In the following, we will assume without loss of generality that $P_v \neq \gekl{\bm{\epsilon}}$ for all $v \in V$, $P_e \neq \gekl{\bm{\epsilon}}$ for all $e \in A$ and $P_{\bar{e}}^{\bar{e}} \neq P_{o(e)}$ for all $e \in T$.

\blemma
\label{lem:I.Ex}
If there exists $e \in T$ with $P_e = \gekl{\bm{\epsilon}}$, then I. is satisfied.
\elemma
\nopar

\bproof
Let $x \in P_v \setminus \gekl{\bm{\epsilon}}$ or $x \in A$. In the latter case, set $v \defeq t(x)$. Let $\chi$ and $W$ be as in I.
\pari

First assume that there exists a strictly increasing sequence $(j_N)_N$ of positive integers such that, with $W_N \defeq W_{\text{---}j_N}$, $W_N$ can be represented by a properly reduced positive word with first letter in $P_v$ of first letter in $E$ with origin $v$, for all $N$. Assume that $[v,o(e)]$ does not contain $t(e)$, otherwise replace $e$ by $\bar{e}$. Take $Y \in P_{t(e)} \setminus \gekl{\bm{\epsilon}}$. Then $x Y W_N$ is reduced, and we can assume without loss of generality that $x Y W_N$ is properly reduced (when we replace $x$ and $W_N$ by suitable positive words representing them). Suppose that $x \in P_v$, the case $x \in A$ is similar. If $p_1 \dotsm p_{\mu} d_1 \dotsm$ is a properly reduced positive word representing $x Y W_N$, then we have $x = p_1 a$ or $x a = p_1$. In the first case, we are done. The second case leads to $a = \bm{\epsilon}$ using that $P_{\bar{e}}^{\bar{e}} = \gekl{\bm{\epsilon}}$.

Now assume that there exists a strictly increasing sequence $(j_N)_N$ of positive integers such that, with $W_N \defeq W_{\text{---}j_N}$, $W_N$ can be represented by a properly reduced positive word with first letter not in $P_v$ of first letter in $E$ with origin not equal to $v$, for all $N$. Assume that $[v,o(e)]$ does not contain $t(e)$, otherwise replace $e$ by $\bar{e}$. Take $y_t \in P_{t(e)} \setminus \gekl{\bm{\epsilon}}$ and $y_v \in P_v \setminus P_f^f$, where $[t(e),v]$ ends with $f \in T$. Define $Y \defeq y_t y_v$. Then $x Y W_N$ is reduced, and we can assume without loss of generality that $x Y W_N$ is properly reduced (when we replace $x$ and $W_N$ by suitable positive words representing them). The same argument as in the first case shows that I. holds.
\eproof
\pars

To get examples satisfying II., we now assume that $G_v \subseteq (\Rz,+)$ for all $v \in V$. We will still use multiplicative notation. 

\blemma
If $P_{\bar{e}}^{\bar{e}} \neq P_{o(e)}$ and $P_e^e \neq P_{t(e)}$ for all $e \in A$, then condition (LCM) implies that $P_{\bar{e}}^{\bar{e}}$ is not dense in $P_{o(e)}$, for all $e \in A \amalg T$, and that $P_e^e$ is not dense in $P_{t(e)}$, for all $e \in A \amalg T$. Thus $P_e \cong \Zz_{\geq 0}$ or $P_e = \gekl{\bm{\epsilon}}$ for all $e \in A \amalg T$.
\elemma
\nopar

\bproof
Otherwise, suppose we can find $p \in P_{o(e)} \setminus P_{\bar{e}}^{\bar{e}}$ and a sequence $(p_n)$ in $P_{\bar{e}}^{\bar{e}}$ such that $p \prec p_n$ and $\lim_{n \to \infty} p_n = p$. Then $p^{-1} p_n \in p^{-1}(P_{\bar{e}}^{\bar{e}})$. Hence $p^{-1} p_n \in p^{-1,e} P_{\bar{e}}^{\bar{e}}$, which implies that $p^{-1,e} \prec p^{-1} p_n$ for all $n$. This contradicts $\lim_{n \to \infty} p_n = p$.
\eproof
\pars

This motivates the following
\bdefin
\label{def:D}
Assume that $G_v \subseteq (\Rz,+)$ for all $v \in V$. Then we say that condition (D) is satisfied if $P_e \cong \Zz_{\geq 0}$ or $P_e = \gekl{\bm{\epsilon}}$ for all $e \in A \amalg T$.
\edefin

\blemma
\label{lem:II.Ex}
Assume that $G_v \subseteq (\Rz,+)$ for all $v \in V$, and that conditions (LCM) and (D) are satisfied. If $P_e \neq \gekl{\bm{\epsilon}}$ for all $e \in T$ and $\# V > 1$ or $\# A_+ > 0$, then II. is satisfied.
\elemma
\nopar

\bproof
First assume $\# A_+ > 0$. Take $e \in A_+$. Let $\bm{b}$ be the generator of $P_e^e \cong \Zz_{\geq 0}$. Take $x \in P_v \setminus \gekl{\bm{\epsilon}}$ or $x \in A$. Set $Y \defeq e$. Let $\chi$ and $W$ be as in II. Let $(j_N)_N$ be arbitrary and define $W_N$ as in II. Then $x Y W_N$ is reduced, and we can assume without loss of generality that $x Y W_N$ is properly reduced (when we replace $x$ and $W_N$ by suitable positive words representing them). Let us now treat the case that $x \in P_v$, the case $x \in A$ is similar.
\pars

Let $p_1 \dotsm p_\mu Y Z$ be a properly reduced positive word representing $x Y W_N$. If $x = p_1 a_1$, then we are done. Otherwise, we must have $p_1 = x a_1$, and then we obtain $a_1 p_2 = a_2, \dotsc, a_{\mu-1} p_\mu = a_\mu$. If $a_1 \neq \bm{\epsilon}$, then $a_2 \neq \bm{\epsilon}, \dotsc, a_\mu \neq \bm{\epsilon}$. So $a_\mu \in P_{\bar{e}}^{\bar{e}} = P_e^e$, and hence $a_\mu = \bm{b}^i$ for some $i > 0$. This implies
$$
 p_1 \dotso p_\mu = x a_1 p_2 \dotso p_\mu = x a_2 p_3 \dotso p_\mu = \dotso = x a_\mu = x \bm{b}^i.
$$
Thus $x \bm{b}^i \in p_1 P_T$ and $x Y W_N = p_1 \dotso p_\mu Y Z = x \bm{b}^i Y Z = x Y \bm{b}^{\bar{i}} Z$ for some $\bar{i} > 0$. We conclude that $W_N \in \bm{b} P$, as desired.

Now assume that $\# V > 1$. Take $e \in T$ and let $\bm{b}$ be the generator of $P_e^e \subseteq P_{t(e)}$. Then $\bm{b}$ is also the generator of $P_{\bar{e}}^{\bar{e}} \subseteq P_{o(e)}$. Take $x \in P_v \setminus \gekl{\bm{\epsilon}}$ or $x \in A$. In the latter case, set $v \defeq t(x)$. Let $\chi$ and $W$ be as in II.

Suppose that there exists a strictly increasing sequence $(j_N)_N$ of positive integers such that, with $W_N \defeq W_{\text{---}j_N}$, $W_N$ can be represented by a properly reduced positive word with first letter in $P_{v_N}$ of first letter in $E$ with origin $v_N$. Assume that $v$ and $v_N$ are on the same side of $e$ for all $N$, and that $[v,o(e)]$ does not contain $t(e)$ (the other case is similar). Take $Y \in P_{t(e)}$ such that $Y \prec z$ and $Y \neq z$ for all $z \in P_e^e \setminus \gekl{\bm{\epsilon}}$. Then $x Y W_N$ is reduced, and we can assume without loss of generality that $x Y W_N$ is properly reduced (when we replace $x$ and $W_N$ by suitable positive words representing them). Let us now treat the case that $x \in P_v$, the case $x \in A$ is similar.

Let $p_1 \dotsm p_\mu Y W'$ be a properly reduced positive word representing $x Y W_N$. If $x = p_1 a_1$, then we are done. If $x a_1 = p_1$, then we obtain $a_1 p_2 = a_2, \dotsc, a_{\mu-1} p_\mu = a_\mu$, with $a_\mu \in P_{\bar{e}}^{\bar{e}} = P_e^e$ and hence $a_\mu = \bm{b}^i$ for some $i$. Moreover, $a_\mu Y' = Y a_{\mu+1}$. As $Y \prec a_\mu$ but $Y \neq a_\mu$, we must have $a_{\mu+1} \neq \bm{\epsilon}$, so that $a_{\mu+1} = \bm{b}^{\bar{i}}$ for some $\bar{i} > 0$. This implies
$$
 p_1 \dotso p_\mu = x a_1 p_2 \dotso p_\mu = x a_2 p_3 \dotso p_\mu = \dotso = x a_\mu = x \bm{b}^i.
$$
Thus $x \bm{b}^i \in p_1 P_T$. Moreover, $x Y W_N = p_1 \dotso p_\mu Y' W' = x a_\mu Y' W' = x Y a_{\mu+1} W' = x Y \bm{b}^{\bar{i}} W'$. We conclude that $W_N = \bm{b}^{\bar{i}} W' \in \bm{b}^{\bar{i}} P \subseteq \bm{b} P$, as desired.

Now suppose that there exists a strictly increasing sequence $(j_N)_N$ of positive integers such that, with $W_N \defeq W_{\text{---}j_N}$, $W_N$ can be represented by a properly reduced positive word with first letter in $P_{v_N}$ of first letter in $E$ with origin $v_N$. Assume that $v$ and $v_N$ are on opposite sides of $e$ for all $N$, say $[v,v_N] = [v,o(e)]e[t(e),v_N]$. Take $y_t \in P_{t(e)}$ such that $y_t \prec z$ and $y_t \neq z$ for all $z \in P_e^e \setminus \gekl{\bm{\epsilon}}$. Furthermore, take $y_o \in P_{o(e)}$ such that $y_o \prec \bar{z}$ and $y_o \neq \bar{z}$ for all $\bar{z} \in P_{\bar{e}}^{\bar{e}} \setminus \gekl{\bm{\epsilon}}$. Define $Y \defeq y_t y_o$. Then $x Y W_N$ is reduced, and we can assume without loss of generality that $x Y W_N$ is properly reduced (when we replace $x$ and $W_N$ by suitable positive words representing them). Let us now treat the case that $x \in P_v$, the case $x \in A$ is similar.

Let $p_1 \dotsm p_\mu y'_t y'_o W'$ be a properly reduced positive word representing $x Y W_N$. If $x = p_1 a_1$, then we are done. If $x a_1 = p_1$, then we obtain $a_1 p_2 = a_2, \dotsc, a_{\mu-1} p_\mu = a_\mu$, with $a_\mu \in P_{\bar{e}}^{\bar{e}} = P_e^e$ and hence $a_\mu = \bm{b}^i$ for some $i$. Moreover, $a_\mu y'_t = y_t a_{\mu+1}$. As $y_t \prec a_\mu$ but $Y \neq a_\mu$, we must have $a_{\mu+1} \neq \bm{\epsilon}$. Similarly, we have $a_{\mu+1} y'_o = y_o a_{\mu+2}$. As $y_o \prec a_{\mu+1}$ but $y_o \neq a_{\mu+1}$, we must have $a_{\mu+2} \neq \bm{\epsilon}$. As $a_{\mu+2}$ lies in $P_e^e = \spkl{\bm{b}}^+$, we must have $a_{\mu+2} = \bm{b}^{\bar{i}}$ for some $\bar{i} > 0$. This implies
$$
 p_1 \dotsm p_\mu = x a_1 p_2 \dotso p_\mu = x a_2 p_3 \dotso p_\mu = \dotso = x a_\mu = x \bm{b}^i.
$$
Thus $x \bm{b}^i \in p_1 P_T$. Moreover, $x Y W_N = p_1 \dotso p_\mu y'_t y'_o W' = x a_\mu y'_t y'_o W' = x y'_t y'_o a_{\mu+2} W' = x y'_t y'_o \bm{b}^{\bar{i}} W'$. We conclude that $W_N = \bm{b}^{\bar{i}} W' \in \bm{b}^{\bar{i}} P \subseteq \bm{b} P$, as desired.
\eproof
\pars

Now we turn to closed invariant subspaces of $\Omega$.

\blemma
\label{lem:I.dense...}
Assume that I. holds. If there exists $v \in V$ and a sequence $x_n \in P_v \setminus \gekl{\bm{\epsilon}}$ with $x_{n+1} \prec x_n$ such that, for every $p \in P_v \setminus \gekl{\bm{\epsilon}}$, $x_n \prec p$ and $x_n \neq p$ for all sufficiently big $n$, then $\partial \Omega = \Omega$.
\elemma
\nopar

\bproof
Let $\chi \in \Omega$ be arbitrary and $Y_n$ and $W^{(n)}$ be as in I for $x = x_n$. We now claim that $\lim_{n \to \infty} \chi_{x_n Y W^{(n)}} = \chi_{\bm{\epsilon}}$. As in Lemma~\ref{lem:I.Ex}, we may assume without loss of generality that $x_n Y W^{(n)}_{N_n}$ is properly reduced for all $N_n$. If $\chi_{x_n Y W^{(n)}}(p) = 1$ then $x_n Y W^{(n)}_{N_n} \in pP$ for all sufficiently big $n$ and $N_n$. Assume that $p \neq \bm{\epsilon}$. Let $p_1 \dotsm p_{\mu} d_1 \dotsm$ be a properly reduced positive word representing $p$. We treat the case $p_1 \in P_{v_1} \setminus \gekl{\bm{\epsilon}}$. The case that our properly reduced word starts with $d_1 \in A$ is straightforward. $x_n Y W^{(n)}_{N_n} \in pP$ means that $x_n Y W^{(n)}_{N_n} = pz$ for some $z$. By Lemma~\ref{lem:XY_m}~(i), there is a properly reduced positive word with first letter $p_1 z'$ representing $pz$, and $z = z' z''$. Comparing properly reduced positive words, we must have $p_1 z' \in P_v$ by Lemma~\ref{lem:W=W'implies...}. If $v_1 = v$, then I. implies $p_1 \prec x_n$ for sufficiently big $n$, which contradicts our choice of $x_n$. If $v_1 \neq v$, let $[v,v_1]$ start with $d \in T$. Lemma~\ref{lem:px=q} implies that $z' \in p_1^{-1,d} P$. Hence I. implies that $x_n \in p_1 z' P \subseteq p_1 p_1^{-1,d} P$ and thus $x_n \in p_1 p_1^{-1,d} P_v$. In other words, $p_1 p_1^{-1,d} \prec x_n$ for all sufficiently big $n$. This again contradicts our choice of $x_n$.
\eproof
\pars

The following is an immediate consequence of Lemmas~\ref{lem:I.Ex}, \ref{lem:I.dense...}, \cite[Theorem~5.7.2 and Corollary~5.7.17]{CELY}.
\bcor
\label{cor:I,pis}
Assume that condition (LCM) is satisfied. If there exists $e \in T$ with $P_e = \gekl{\bm{\epsilon}}$, and if there exists $v \in V$ and a sequence $x_n \in P_v \setminus \gekl{\bm{\epsilon}}$ with $x_{n+1} \prec x_n$ such that, for every $p \in P_v \setminus \gekl{\bm{\epsilon}}$, $x_n \prec p$ and $x_n \neq p$ for all sufficiently big $n$, then $C^*_{\lambda}(P)$ is purely infinite simple.
\ecor

Now we assume that $G_v \subseteq (\Rz,+)$ for all $v \in V$. Our goal is to determine all closed invariant subspaces of $\Omega$ in the case where $\# V > 1$ or $\# V = 1$, say $V = \gekl{v}$, $G_v \subseteq (\Rz,+)$ dense and $\#A > 0$. The remaining case that $\# V = 1$, $V = \gekl{v}$, $G_v \cong \Zz$, $P_v \cong \Zz_{\geq 0}$ and $\#A > 0$ will be treated separately, and the case that $P \subseteq (\Rz,+)$ is treated in \cite{Li2015}.

The following notation will be convenient: If $\# V > 1$ and $P_e \neq \gekl{\beps}$ for all $e \in T$, then let $\bm{b}$ be the generator of $P_e^e \cong \Zzg$ for some $e \in T$, and if $\# V = 1$, say $V = \gekl{v}$, choose $\bm{b} \in P_v \setminus \gekl{\beps}$ arbitrary. Define
$$
 \Omega_{\bm{b},\infty} \defeq \menge{\chi \in \Omega}{(g.\chi)(\bm{b}^i) = 1 \ \forall \ g, \, i},
$$
where $i$ runs through all natural numbers, and we only consider those $g \in G$ such that $g.\chi$ is defined. Note that this is independent of the choice of $\bm{b}$, coincides with our previous definition (in those cases where there is an overlap), and we always have $\Omega_{\bm{b},\infty} \subseteq \Omega_{\infty}$.

Our goal is to prove the following:
\btheo
\label{thm:clinvsub_gen}
Suppose that $G_v \subseteq (\Rz,+)$ for all $v \in V$, and that $\# V > 1$ or $\# V = 1$, $V = \gekl{v}$, $G_v \subseteq (\Rz,+)$ dense and $\#A > 0$. Further assume that conditions (LCM) and (D) are satisfied.
\nopar

\begin{enumerate}
\item[(i)] Assume that there exists $e \in T$ with $P_e = \gekl{\beps}$. 
\begin{enumerate}
\item[(i$_1$)] If there exists $v \in V$ such that $G_v$ is dense in $\Rz$, then the following is the list of all closed invariant subspaces of $\Omega$: $\partial \Omega = \Omega$.
\item[(i$_2$)] If $P_v \cong \Zz_{\geq 0}$ for all $v \in V$, then the following is the list of all closed invariant subspaces of $\Omega$: $\partial \Omega = \overline{\Omega_{\infty}} \subseteq \Omega$.
\end{enumerate}
\item[(ii)] Assume that $P_e \neq \gekl{\beps}$ for all $e \in T$.
\begin{enumerate}
\item[(ii$_1$)] If there exists $v \in V$ such that $G_v$ is dense in $\Rz$ and $\# A \geq 1$, then the following is the list of all closed invariant subspaces of $\Omega$: $\Omega_{\bm{b},\infty} = \partial \Omega \subsetneq \Omega$.
\item[(ii$_2$)] If $\# A = 0$ (and $\# V > 1$), then the following is the list of all closed invariant subspaces of $\Omega$: $\gekl{\infty} = \partial \Omega \subsetneq \overline{\Omega_{\infty}} \subseteq \Omega$.
\item[(ii$_3$)] If $P_v \cong \Zz_{\geq 0}$ for all $v \in V$ (and $\# V > 1$), then the following is the list of all closed invariant subspaces of $\Omega$: $\Omega_{\bm{b},\infty} = \partial \Omega \subsetneq \overline{\Omega_{\infty}} \subseteq \Omega$.
\end{enumerate}
\end{enumerate}
\etheo

For the proof, we need a series of Lemmas.
\pars

\blemma
\label{lem:II.dense}
Suppose that $P_e \neq \gekl{\beps}$ for all $e \in T$, i.e., $P_e \cong \Zz_{\geq 0}$ for all $e \in T$. Then $P_T$ is Ore, and we write $\partial \Omega_{P_T} = \gekl{\infty}$. If $\# A \geq 1$ and there exists $v \in V$ such that $G_v$ is dense in $\Rz$, then $\partial \Omega = \Omega_{\bm{b},\infty}$. Moreover, for every $\chi \notin \Omega_{\bm{b},\infty}$, we have $\overline{G.\chi} = \Omega$. If $\# A = 0$, then for all $\chi \neq \infty$ and $\eta \in \Omega_{\infty}$, we have $\eta \in \overline{G.\chi}$.
\elemma
\nopar

\bproof
In the first case, take $e \in A$ and a strictly decreasing sequence $(x_n)$ in $P_v$ such that $\lim_{n \to \infty} x_n = \bm{\epsilon}$. Let $\chi \in \Omega$ be arbitrary and write $\chi = \chi_W$ for some infinite positive word $W$. By compactness, we can --- by passing to a subsequence if necessary --- assume that $\psi \defeq \lim_{n \to \infty} \chi_{x_n e W}$ exists. We claim that $\psi \in \Omega_{P_T}$. Indeed, if not, then we must have $\psi(p e) = 1$ for some $p \in P_T$. It follows that $p G_e^{\bar{e}} = x_n G_e^{\bar{e}}$ for all $n$. Hence $x_m G_e^{\bar{e}} = x_n G_e^{\bar{e}}$ for all $m$ and $n$. But this contradicts $\lim_{n \to \infty} x_n = \bm{\epsilon}$. So we obtain that $\Omega_{P_T} \cap \overline{G.\chi} \neq \emptyset$, so that $\infty \in \overline{G.\chi}$. We conclude that $\Omega_{\bm{b},\infty} \subseteq \partial \Omega$, as desired. Now we show $\overline{G.\chi} = \Omega$ for every $\chi \notin \Omega_{\bm{b},\infty}$. We may assume that $\chi(\bm{b}) = 0$. If $\# V > 1$ or $\# A_+ > 0$, then a similar argument as in Lemma~\ref{lem:II.gen} (or Lemma~\ref{lem:I.gen}) shows the following: If we take $e \in A$ and a sequence $(x_n)$ in $P_v$ such that $\lim_{n \to \infty} x_n = \bm{\epsilon}$ and write $\chi = \chi_W$ for some infinite positive word $W$, then $\lim_{n \to \infty} \chi_{x_n e W} = \chi_{\bm{\epsilon}}$. If $\# V = 1$ and $A = A_-$, and if we write $\chi = \chi_W$ for some positive word $W$, then $\chi(\bm{b}) = 0$ implies that no $e \in A_-$ can appear in $W$, so that $\chi \in \Omega_{P_v} \setminus \gekl{\infty}$. Now our claim follows because $G_v \curvearrowright \Omega_{P_v} \setminus \gekl{\infty}$ is minimal (see \cite{Li2015}).
\pars

Now we turn to the second case. If $\eta = \chi_X$ for some infinite positive word $X$ with $\lim_{l \to \infty} \ell(X_{\text{---}l}) = \infty$, then we already know that $\eta \in \overline{G.\chi}$. Otherwise Lemma~\ref{lem:linftyvsPv} implies that $\eta \in \Omega_{P_v}$ for some $v \in V$. If $P_v \cong \Zz_{\geq 0}$, then $\eta \in \Omega_{\infty}$ implies $\eta = \infty$, and our claim follows. If $P_v \not\cong \Zz_{\geq 0}$, then $G_v$ must be dense in $\Rz$. Let $(x_n)$ be a sequence in $P_v$ such that $\eta = \lim_{n \to \infty} \chi_{x_n}$. Without loss of generality we may assume $\chi(\bm{b}) = 0$. Let $Y$ and $W$, $W_N$ be as in II. for $x = x_n$. Note that in the proof of Lemma~\ref{lem:II.Ex}, $Y$ and $W_N$ were constructed so that they only depend on $v$, not on $x_n$. Moreover, as in the proof of Lemma~\ref{lem:II.Ex}, the first letter of $Y$ lies in $P_t$, and suppose that $[v,t]$ starts with $d \in T$. Without loss of generality we may assume that $x_n \prec z$ and $x_n \neq z$ for all $z \in P_d^{\bar{d}}$. This is because $G_v \curvearrowright \Omega_{P_v, \infty} \setminus \gekl{\infty}$ is minimal (see \cite{Li2015}). We claim that $\lim_{n \to \infty} \chi_{x_n Y W} = \eta$. Indeed, suppose that $\chi_{x_n Y W}(p) = 1$. Then $x_n Y W_N \in pP$. As before, $x_n Y W_N$ is reduced, and we can assume without loss of generality that $x_n Y W_N$ is properly reduced (when we replace $x_n$ and $W_N$ by suitable positive words representing them). Suppose that $W_p = p_0 p_1 \dotso p_m$ is a properly reduced word representing $p$, with $p_k \in P_{v_k}$. We proceed inductively on $\mfl(W_p)$ to show that $x_n \in pP$. $x_n Y W_N \in pP$ implies that $x_n Y W_N = pz$ for some $z$ in $P$. If $\mfl(W_p) = 0$, then $p = p_0$, and Lemma~\ref{lem:XY_m} implies that $pz$ can be represented by a properly reduced positive word with first letter of the form $p_0 z'$. Now II. implies that $x_n \in (p_0 z') P_T$ as otherwise, we would get $W_N \in \bm{b}P$, contradicting $\chi(\bm{b}) = 0$. Now suppose that $\mfl(W_p) \geq 1$. First let $W_z$ be a properly reduced positive word representing $z$. If $W_p W_z$ is properly reduced, then Lemma~\ref{lem:W=W'implies...} implies that $p_0 \in P_v$ and $[v_0,v_1]$ must start with $d$. As before, II. and $\chi(\bm{b}) = 0$ imply that $x_n = p_0 a$ for some $a \in P_d^{\bar{d}}$. But $x_n \prec z$ and $x_n \neq z$ for all $z \in P_d^{\bar{d}}$ implies $a = \bm{e}$, and we are done. If $W_p W_z$ is not properly reduced, then we can write $W_p W_z \equiv (W_p W'_z) W''_z$ such that $\mfl(W_p W'_z) < \mfl(W_p)$. By induction hypothesis, we obtain $x_n \in (W_p W'_z) P \subseteq W_p P = pP$, as desired.
\eproof
\pars

\blemma
\label{lem:InfIsIn}
Assume that $P_v \cong \Zz_{\geq 0}$ for all $v \in V$ and that $\# V > 1$. Then for all $v \in V$ and $\chi \in \Omega$, $\infty_{P_v} \in \overline{G.\chi}$.
\elemma
\nopar

\bproof
Suppose that there exists $e \in T$ with $P_e = \gekl{\beps}$. Then condition I. holds by Lemma~\ref{lem:I.Ex}. 
\pari

If $P_v *_{P_e} P_w \subseteq P_T$, write $P_v = \spkl{\alpha}^+$ and $P_w = \spkl{\beta}^+$. Let $X = \alpha \beta \alpha \beta \dotso$. Then $\alpha^{\infty} = \infty_{P_v *_{P_e} P_w} \in \overline{G. \chi_X}$, and Lemma~\ref{lem:linftyvsPv} implies that $\chi_X \in \overline{G. \chi}$.

If $P_v *_{\gekl{\bm{\epsilon}}} P_w \subseteq P_T$, write $P_v = \spkl{\alpha}^+$ and $P_w = \spkl{\beta}^+$ as before. We claim that $\lim_{n \to \infty} \chi_{\alpha^n \beta \alpha \beta \alpha \dotso} = \chi_{\alpha^{\infty}}$. Indeed, if $\chi_{\alpha^n \beta \alpha \beta \alpha \dotso}(p) = 1$ for all sufficiently big $n$, then we must have $\alpha^n \beta \alpha \dotso \beta \alpha \in pP$ for all sufficiently big $n$. As $p$ is fixed, this implies that $\alpha^n \in pP$ for all sufficiently big $n$. Thus $\chi_{\alpha^{\infty}}(p) = 1$. Therefore, we may take $X = \beta \alpha \beta \alpha \dotso$. Then $\alpha^{\infty} \in \overline{G. \chi_X}$, and Lemma~\ref{lem:linftyvsPv} implies that $\chi_X \in \overline{G. \chi}$.
\pars

Now suppose that $P_e \neq \gekl{\beps}$ for all $e \in T$. Then $P_T$ is Ore, and $\Omega_{\infty} \cap \Omega_{P_v} = \gekl{\infty}$. So it suffices to show that $\Omega_{P_T} \cap \overline{G. \chi} \neq \emptyset$. Take $w \in V$ with $w \neq v$ and $\alpha$, $\beta$ as above. Take $\chi \in \Omega$. If $\chi \in \Omega_{P_T}$, then there is nothing to show. If $\chi \notin \Omega_{P_T}$, then there exists $q \in P_T$ and $e \in A$ with $\chi(qe) = 1$. By compactness, we can find a sequence $n_i$ such that $(\alpha \beta)^{n_i}. \chi$ converges to $\eta$. We claim that $\eta \in \Omega_{P_T}$. If not, then there exists $p \in P_T$ such that $\eta(p e) = 1$. It follows that $(\alpha \beta)^{n_i} q G_e^{\bar{e}} = p G_e^{\bar{e}}$ and thus $(\alpha \beta)^{n_i} q G_e^{\bar{e}} = (\alpha \beta)^{n_j} q G_e^{\bar{e}}$ for all $i, j$. Hence, if we set $m_j \defeq n_j - n_1$, then $(\alpha \beta)^{m_j} q = q g_j$ for some $g_j \in G_e^{\bar{e}}$. Consider a reduced $\bm{o}$-word representing $q g_j$. Its path in $T$ is of finite length, and the length is independent of $j$. Now consider reduced $\bm{o}$-words representing $(\alpha \beta)^{m_j} q$. It is easy to see that the lengths of the paths in $T$ corresponding to these reduced $\bm{o}$-words tend to infinity as $j \to \infty$ because $m_j \to \infty$. So this is a contradiction, as desired.
\eproof
\pars

\blemma
\label{lem:binfNOTinf}
Suppose that $P_v \cong \Zz_{\geq 0}$ for all $v \in V$ and $P_e \neq \gekl{\bm{\epsilon}}$ for all $e \in T$. Then $\Omega_{\bm{b},\infty} = \Omega_{\infty}$ if and only if $\# V = 1$ and $A = A_-$. In particular, if $\# V > 1$, then $\Omega_{\bm{b},\infty} \subsetneq \overline{\Omega_{\infty}}$.
\elemma
\nopar

\bproof
\an{$\Leftarrow$} is clear. For \an{$\Rightarrow$}, first suppose that $\# V > 1$. Take $v, w \in V$ with $v \neq w$ and let $\alpha$ and $\beta$ be the generators of $P_v$ and $P_w$, and set $X \defeq \alpha \beta \alpha \beta \dotso$. Then we have $\chi_X \in \Omega_{\infty} \setminus \Omega_{\bm{b},\infty}$. If $A \neq A_-$, then take $e \in A_+$ and set $X \defeq e^{\infty} = e e e \dotso$. Then $\chi_X \in \Omega_{\infty} \setminus \Omega_{\bm{b},\infty}$. 
\eproof
\pars

\bproof[Proof of Theorem~\ref{thm:clinvsub_gen}]
(i$_1$) follows from Lemma~\ref{lem:I.dense...}.
\pari

(i$_2$) follows from Lemma~\ref{lem:linftyvsPv}, Lemma~\ref{lem:InfIsIn} and the fact that $\Omega_{P_v} = P_v \cup \gekl{\infty_{P_v}}$.

(ii$_1$) and (ii$_2$) follow from Lemma~\ref{lem:II.dense}.

Finally, we explain how to derive (ii$_3$): Clearly, $\Omega_{\bm{b},\infty} \supseteq \partial \Omega$. Now take $\chi \in \partial \Omega$ and $\eta \in \Omega_{\bm{b},\infty}$. Lemma~\ref{lem:linftyvsPv}~(ii) implies that $\eta \in \overline{G.\chi}$ or $\eta = \infty$. In the latter case, Lemma~\ref{lem:InfIsIn} implies $\eta \in \overline{G.\chi}$ as well. Now take $\chi \notin \Omega_{\bm{b},\infty}$. We may assume that $\chi(\bm{b}) = 0$. Take $\eta \in \Omega_{\infty}$ arbitrary. Then Lemma~\ref{lem:linftyvsPv}~(i) implies that $\eta \in \overline{G.\chi}$ or $\eta = \infty$. In the latter case, Lemma~\ref{lem:InfIsIn} implies $\eta \in \overline{G.\chi}$ as well. Moreover, Lemma~\ref{lem:binfNOTinf} implies that $\Omega_{\bm{b},\infty} \subsetneq \overline{\Omega_{\infty}}$. This completes the proof.
\eproof
\pars

Note that $\Omega_{\infty}$ is either closed or $\overline{\Omega_{\infty}} = \Omega$. For completeness, we characterize when the latter happens.
\blemma
\label{lem:Omegainf-vs-Omega}
Suppose that $G_v \subseteq (\Rz,+)$ for all $v \in V$. Then $\overline{\Omega_{\infty}} = \Omega$ if and only if one of the following is satisfied:
\nopar

\begin{enumerate}
\item[(a)] There exists $v\in V$ such that $G_{v}$ is dense in $\mathbb{R}$;
\item[(b)] $P_{v}\cong \mathbb{Z}_{\geq 0}$ for all $v\in V$ and $\# V=\infty$;
\item[(c)] $P_{v}\cong \mathbb{Z}_{\geq 0}$ for all $v\in V$ and $\# A_{+}=\infty$.
\end{enumerate}
\end{lemma}
\nopar

\bproof
If (a) holds, then $\overline{\Omega_{\infty}} = \Omega$ because $G_v \curvearrowright \Omega_{P_v} \setminus \gekl{\infty}$ is minimal (see \cite{Li2015}). Now assume $P_{v}\cong \mathbb{Z}_{\geq 0}$ for all $v\in V$. If $\# V=\infty$, set $V:=\{v_{i}\}_{i\in \mathbb{N}}$. Let $b_{v_{i}}$ be the generator of $P_{v_{i}}$ and define $\chi_{n}:=\chi_{b_{v_{n}} b_{v_{n+1}} b_{v_{n+2}} \dotsm},\ n\in \mathbb{N}$. It is easy to see that $\text{lim}_{n\rightarrow \infty}\chi_{n}=\chi_{\beps}$. If $\# A_{+}=\infty$, set $A_{+}:=\{e_{i}\}_{i\in \mathbb{N}}$. Define $\chi_{n}:=\chi_{e_{n} e_{n+1} e_{n+2} \dotsm},\ n\in \mathbb{N}$. It is easy to see that $\text{lim}_{n\rightarrow \infty}\chi_{n}=\chi_{\beps}$.
\pars

Conversely, assume $P_{v}\cong \mathbb{Z}_{\geq 0}$ for all $v\in V$, $\# V<\infty$ and $\# A_{+}<\infty$. Let $b_{v}$ be the generator of $P_{v}$. Then for all $\chi\in \Omega_{\infty}$, either $\chi(b_{v})=1$ for some $v\in V$ or $\chi(e)=1$ for some $e\in A_{+}$. Take a convergent sequence $(\chi_{n})_{n}\subseteq \Omega_{\infty}$, then either there exists $v\in V$ such that $\chi_{n}(b_{v})=1$ for all sufficiently big $n$ or there exists $e\in A_+$ such that $\chi_{n}(e)=1$ for all sufficiently big $n$, which implies $\text{lim}_{n\rightarrow \infty}\chi_{n}\neq \chi_{\beps}$.\\
\eproof
\pars

\subsection{The one vertex GBS case}

Let us now consider the remaining case of one vertex generalized Baumslag-Solitar (abbreviated by one vertex GBS) groups and their monoids, as in Example~\ref{ex:pi}~(ii). Suppose that we are in the same setting as introduced at the beginning of \S~\ref{s:clinvsub}, with $\# V = 1$, $V = \gekl{v}$, $P_v \cong \Zzg$, $\# A > 0$ and $P_e \cong \Zzg$ for all $e \in A$. Then our one vertex GBS monoid admits the presentation
$$
 P = \spkl{ \gekl{b} \cup A \ \vert \ b^{n_e} e = e b^{m_e} \ \forall \ e \in A_+, \, b^{n_e} e b^{m_e} = e \ \forall \ e \in A_-}^+.
$$
Here $m_e, n_e \in \Zz$ with $m_e, n_e \geq 1$ for all $e \in A$. The enveloping group $G$ of $P$ admits the same presentation.

For each $e \in A$, define a homomorphism $\theta_e: \: P \to \Zzg$ by $\theta_e(b) \defeq 0$, $\theta_e(e') \defeq \delta_{e,e'}$ for all $e' \in A$. We set $\theta_+ \defeq \sum_{e \in A_+} \theta_e$, $\theta_- \defeq \sum_{e \in A_-} \theta_e$ and $\theta \defeq \theta_+ + \theta_-$. We extend these maps to infinite words by allowing the value $\infty$.

We need the following standard forms, which follows by the results in \cite{Serre, Bass} (see also \cite[\S~2.2]{BMPST}).
\bprop
Each element of $G$ has unique representations in the two forms
\nopar

\begin{enumerate}
\item[(L)] $b^{j_{0}}e_1^{\iota_1}b^{j_{1}}e_2^{\iota_2} \dotsm e_p^{\iota_p}b^{j_{p}}$, where $\iota_{\bullet} \in \{\pm 1\}$, $0 \leq j_{\bullet} < n_{e_{\bullet +1}}$ if $\iota_{\bullet + 1} =1$, and $0 \leq j_{\bullet } < m_{e_{\bullet +1}}$ if $\iota_{\bullet + 1} =-1,\ j_{p}\in \mathbb{Z}$;
\item[(R)] $b^{j_{0}}e_1^{\iota_1}b^{j_{1}}e_2^{\iota_2} \dotsm e_p^{\iota_p}b^{j_{p}}$, where $\iota_{\bullet} \in \{\pm 1\}$, $0 \leq j_{\bullet} < m_{e_{\bullet}}$ if $\iota_{\bullet} =1$, and $0 \leq j_{\bullet} < n_{e_{\bullet}}$ if $\iota_{\bullet} =-1,\ j_{0}\in \mathbb{Z}$.
\end{enumerate}

Each element of $P$ has unique representations in the two forms
\begin{enumerate}
\item[(L)] $b^{j_{0}} e_1 b^{j_{1}} e_2 \dotsm e_p b^{j_{p}}$, $0 \leq j_{\bullet} < n_{e_{\bullet +1}},\ j_{p}\in \mathbb{Z}$;
\item[(R)] $b^{j_{0}} e_1 b^{j_{1}} e_2 \dotsm e_p b^{j_{p}}$, $0 \leq j_{\bullet} < m_{e_{\bullet}},\ j_{0}\in \mathbb{Z}$.
\end{enumerate}
\eprop
\pars

Our aim is to classify all closed invariant subspaces of $\Omega$. In order to state our main result, we introduce the following notation.
\bdefin
Define $\Omega_{A,\infty} \defeq \menge{\chi_w \in \Omega}{\theta(w) = \infty}$ and $\Omega_{b,\infty} \defeq \menge{\chi \in \Omega}{(g.\chi)(b^i) = 1 \ \forall \ g, \, i}$, where $i$ runs through all natural numbers, and we only consider those $g \in G$ such that $g.\chi$ is defined in the definition of $\Omega_{b,\infty}$.
\edefin
\nopar

Note that we always have $\Omega_{A,\infty}, \Omega_{b,\infty} \subseteq \Omega_{\infty}$, and that our definition of $\Omega_{b,\infty}$ is consistent with previous notation.
\pars

Our main result reads as follows:
\btheo
\label{thm:clinvsub_GBS}
Let $P$ be a one vertex GBS monoid. The closed invariant subspaces of $\Omega$ are given as follows (with inclusion relations precisely as indicated):
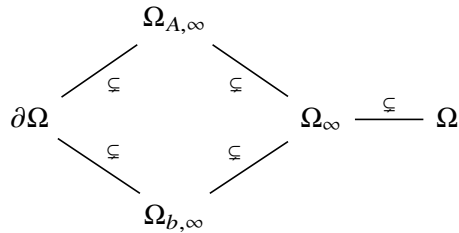
\nopar

\begin{enumerate}
\item[(i)]
$$
\begin{tikzcd}
 & \Omega_{A,\infty} \ar[dash, "\subsetneq"]{ld} & & \\
\partial \Omega \ar[dash, "\subsetneq"]{rd} & & \Omega_{\infty} \ar[dash, , "\subsetneq"]{lu} \ar[dash, "\subsetneq"]{r} & \Omega\\
 & \Omega_{b,\infty} \ar[dash, "\subsetneq"]{ru} & & 
\end{tikzcd}
$$
if $0 < \# A_+ < \infty$, $0 = \# A_-$;
\item[(ii)]
$
 \partial \Omega = \Omega_{A,\infty} \subsetneq \Omega_{\infty} = \Omega_{b,\infty} \subsetneq \Omega
$
if $0 = \# A_+$, $0 < \# A_- < \infty$;
\item[(iii)]
$
 \partial \Omega \subsetneq \Omega_{b,\infty} \subsetneq \Omega_{\infty} \subsetneq \Omega
$
if $0 < \# A_+ < \infty$, $0 < \# A_- < \infty$;
\item[(iv)]
$
 \partial \Omega = \Omega_{b,\infty} = \Omega_{\infty} \subsetneq \Omega
$
if $0 = \# A_+$, $\# A_- = \infty$;
\item[(v)]
$
 \partial \Omega = \Omega_{b,\infty} \subsetneq \Omega_{\infty} \subsetneq \Omega
$
if $0 < \# A_+ < \infty$, $\# A_- = \infty$;
\item[(vi)]
$
 \partial \Omega = \Omega_{b,\infty} \subsetneq \Omega
$
if $\# A_+ = \infty$.
\end{enumerate}
\etheo
\pars

Let us now prove Theorem~\ref{thm:clinvsub_GBS}. We start with the following description of $\Omega_{\max}$.
\bprop
\label{prop:Omegamax_GBS}
Let $\chi_w \in \Omega_{\infty}$. We have $\chi_w \in \Omega_{\max}$ if and only if $\theta_-(w) = \infty$ or $\theta_-(w) < \infty$ and the following hold:
\nopar

\begin{enumerate}
\item[(a)] $\theta_+(w) = \infty$, 
\item[(b)] there exists $j \in \Nz$ and $w \equiv w_{\text{---}j} w'$ such that $\theta_-(w') = 0$ and $\chi_{w'}(b^i) = 1$ for all $i \in \Nz$.
\end{enumerate}
\eprop
\nopar

Note that (b) is equivalent to saying that $\chi_w \in \Omega_{b,\infty}$.
\bproof
Suppose that $\theta_-(w) = \infty$. Take $x\in P$ with $\chi_{w}(x)=0$, and take $y\in P$ with $\theta(y)>\theta(x)$, $\theta_e(y) > \theta_e(x)$ for some $e \in A_-$ and $\chi_{w}(y)=1$. Such $y$ exists because $\theta_-(w) = \infty$. We claim $xP\cap yP=\emptyset$.
Let
$$x=b^{j_{0}} e_{1} b^{j_{1}} e_2 \dotsm b^{j_{k-1}} e_k b^{p}$$
be its standard L-form, and let $x'=b^{j_{0}} e_1 b^{j_{1}} e_2 \dotsm b^{j_{k-1}} e_k$. If $xP\cap yP\neq \emptyset$, there exist $r,\ s,\ t\in P$ such that $r=xs=yt$ and $xP\cap yP=rP$. $xs$ and $yt$ admit the same standard L-form, so $x'$ is a prefix of the standard L-form of $yt$ and hence of $y$. That is, $y=x'z$ for some $z\in P$ and $\chi_{w}(x')=1$. If $p\leq 0$, then $x'P\subseteq xP$ and thus $\chi_{w}(x)=1$, contradicting our choice of $x$. If $p>0$, then since $\theta_e(y)>\theta_e(x)$, we have $\theta_e(z)>0$. In this case, we have $z\in b^{p}P$ and thus $y\in xP$. This again leads to the conclusion $\chi_{w}(x)=1$, contradicting our choice of $x$. This proves our claim. Lemma~\ref{lem:Omegamax_genchar} implies that $\chi_{w}\in \Omega_{\text{max}}$.
\pars

Now suppose that $\theta_-(w) < \infty$. Since $\Omega_{\max}$ is invariant, it suffices to treat the case that $\theta_-(w) = 0$. We know that $b^i P \cap x P \neq \emptyset$ for all $i \in \Nz$ and $x \in P$. If $\chi_w \in \Omega_{\max}$, then Lemma~\ref{lem:Omegamax_genchar} implies that $\chi_w(b^i) = 1$ for all $i \in \Nz$ and $\theta_+(w) = \infty$ because $\chi_{b^{\infty}} \notin \Omega_{\max}$. Conversely, suppose that $w$ satisfies (a) and (b). For the sake of contradiction, assume that $\chi_w \notin \Omega_{\max}$. Then Lemma~\ref{lem:Omegamax_genchar} implies that there exists $x\in P$ with $\chi_{w}(x)=0$ such that for all $y\in P$ with $\chi_{w}(y)=1$, $xP\cap yP\neq \emptyset$. Let 
$$x=b^{j_{0}} e_1 b^{j_{1}} e_2 \dotsm b^{j_{k-1}} e_k b^{p}$$
be its standard L-form and let
$$x'=b^{j_{0}} e_1 b^{j_{1}} e_2 \dotsm b^{j_{k-1}} e_k.$$
Take $y\in P$ with $\theta(y) > \theta(x)$. Since $xP\cap yP\neq \emptyset$, there exist $r,\ s,\ t\in P$ such that $r=xs=yt$ and $xP\cap yP=rP$. $xs$ and $yt$ admit the same standard L-form, so $x'$ is a prefix of the standard L-form of $yt$ and hence of $y$. That is, $y=x'z$ for some $z\in P$ and $\chi_{w}(x')=1$. At the same time, $x'P\cap b^{i}P\neq \emptyset$ for all $i\in \mathbb{N}$. Actually, there exist $j\in \mathbb{N}$ and $\ x''\in P$ such that $x'b^{j}=b^{i}x''$ and $x'P\cap b^{i}P=x'b^{j}P$. Furthermore, when $i$ goes up to infinity, $j$ also tends to $\infty$. Take $i$ big enough such that $j>p$ and hence that $x'b^{j}P\subseteq xP$. Since $\chi_{w}(x')=\chi_{w}(b^{i})=1$, we have $\chi_{w}(x'b^{j})=1$ and hence $\chi_{w}(x)=1$, leading to a contradiction.
\eproof
\pars

\bcor
\label{cor:Omegab=..._GBS}
We always have $\Omega_{b,\infty} = G.\chi_{b^{\infty}} \cup \Omega_{\max}$.
\ecor

Now we describe $\partial \Omega$. 
\bprop
\label{prop:bOmega_GBS}
\begin{enumerate}
\item[(i)] If $\# A < \infty$, then $\chi_{b^{\infty}} \notin \partial \Omega$, and $\partial \Omega = \Omega_{\max}$.
\item[(ii)] If $\# A = \infty$, then $\chi_{b^{\infty}} \in \partial \Omega$, and $\partial \Omega = \Omega_{b,\infty}$.
\end{enumerate}
\eprop
\nopar

\bproof
(i) Assume $\chi_{b^{\infty}} \in \partial \Omega$. Then there exists a sequence $\{\chi_{w_{i}}\}_{i}\subseteq \Omega_{\text{max}}$ such that $\chi_{w_{i}}$ converges to $\chi_{b^{\infty}}$. For each $\chi_{w_{i}}$, there exist $e \in A$ and $0 \leq j < n_e$ such that $\chi_{w_{i}}(b^j e)=1$. Here we use that $\chi_{w_i} \in \Omega_{\max}$. Since there are only finitely many possible values for the pair $(e,j)$, there must be some common $e \in A$ and $0 \leq j < n_e$ such that $\chi_{w_{i}}(b^j e)=1$ for infinitely many $i$. Taking the limit, we get $\chi_{b^{\infty}}(b^j e)=1$, which contradicts the fact $\chi_{b^{\infty}}^{-1}(1) = \{P,\ bP,\ b^{2}P,\ \cdots\}$.

(ii) Let $e_n$, $n \in \Nz$, be pairwise distinct elements of $A$. Let $k_n \in \Nz$ grow sufficiently fast such that, with $w=b^{k_{0}} e_1 b^{k_{1}} e_2 b^{k_{2}} e_3 \dotsm$, we have $\chi_{w}\in \Omega_{\text{max}}$. Set $$w_m:=b^m b^{k_{m}} e_{m+1} b^{k_{m+1}} e_{m+2} b^{k_{m+2}} e_{m+3} \dotsm.$$
Then $\lim_{m \to \infty}\ \chi_{w_m}=\chi_{b^{\infty}}$. Indeed, for any $x\in P$ with $\theta(x)>0$, $\chi_{w_m}(x)=0$ for $m$ big enough. For all $i\in \mathbb{N}$, $\chi_{w_m}(b^{i})=1$ for $m$ big enough. Therefore, $\chi_{b^{\infty}} \in \partial \Omega$.
\eproof
\pars

Now let us study $\Omega_{A,\infty}$.
\bprop
\label{prop:OmegaA_GBS}
\begin{enumerate}
\item[(i)] $\Omega_{A,\infty} = \Omega_{\max}$ if and only if $\# A_+ = 0$.
\item[(ii)] If $\# A_+ > 0$, then for all $\chi \in \Omega_{A,\infty} \setminus \partial \Omega$, we have $\Omega_{A,\infty} \subseteq \overline{G.\chi}$.
\item[(iii)] $\Omega_{A,\infty}$ is closed if and only if $\# A < \infty$ and either $\# A_+ = 0$ or $\# A_- = 0$.
\end{enumerate}
\eprop
\nopar

\bproof
(i) follows from Proposition~\ref{prop:Omegamax_GBS} and the observation that for $e \in A_+$, $\chi_{e^{\infty}} \notin \Omega_{\max}$. 
\pars

(ii) We may assume without loss of generality that $\chi = \chi_w$ with $\theta_-(w) = 0$ and $\chi_w(b) = 0$. Let $w'$ be another infinite word with $\theta_-(w') = 0$. Then we claim that $w'_{\text{---}i}.\chi_{w}$ converges to $\chi_{w'}$. Indeed, take $x\in P$. If $\chi_{w'}(x)=1$, then there exists $M_{2}\in \mathbb{N}$ such that $w'_{\text{---}i}\in xP$ for all $i\geq M_{2}$. For these $i$, $w'_{\text{---}i}\chi_{w}(x)=1$ and hence $\text{lim}\ w'_{\text{---}i}.\chi_{w}(x)=1$. If $\chi_{w'}(x)=0$, we have also $\text{lim}\ w'_{\text{---}i}.\chi_{w}(x)=0$. Otherwise, take $i$ big enough with $w'_{\text{---}i}.\chi_{w}(x)=1$. Since $\chi_{w'}(x)=0$, $w'_{\text{---}i} \notin xP$ and thus $w'_{\text{---}i} w_{\text{---}j} \in xP$ for some $j$. That is, $w'_{\text{---}i} w_{\text{---}j}=xy$ for some $y\in P$. Let
$$x=b^{k_{0}} e_1 b^{k_{1}} e_2 \dotsm b^{k_{M-1}} e_M b^{p}$$
be its standard L-form and let $x'=b^{k_{0}} e_1 b^{k_{1}} e_2 \dotsm b^{k_{M-1}} e_M$. By the uniqueness of the standard L-form, we have $p\geq 0$ and there exists $z\in P$ such that $x'z=w'_{\text{---}i}$ and that $z w_{\text{---}j}=b^{p}y$. Since $\chi_{w} (b)=0$, $w_{\text{---}j} \notin bP$ and thus $z\in b^{p}P$. This means $w'_{\text{---}i}=x'z\in x'b^{p}P=xP$, contradicting the assumption $\chi_{w'}(x)=0$.

(iii) If $\# A < \infty$ and $\# A_+ = 0$, then (i) implies that $\Omega_{A,\infty} = \Omega_{\max}$ and hence $\Omega_{A,\infty} = \partial \Omega$ by Proposition~\ref{prop:Omegamax_GBS}. If $\# A < \infty$ and $\# A_- = 0$, then take a sequence $\{\chi_{i}\}$ in $\Omega_{A,\ \infty}$ converging to some character $\chi \in \Omega$. For $M\in \mathbb{N}$, there exist unique elements $e_{1}, e_{2},\ \cdots,\ e_{M} \in A = A_+$ and $(k_{0},\ k_{1},\ \cdots,\ k_{M-1})$ with $0 \leq k_{\mu} < n_{e_{\mu+1}}$ such that
$$\chi_{i}(b^{k_{0}} e_1 b^{k_{1}} e_2 \dotsm b^{k_{M-1}} e_M)=1$$
for all $i$ big enough. As a result, $\chi(b^{k_{0}} e_1 b^{k_{1}} e_2 \dotsm b^{k_{M-1}} e_M)=1$. Since $M$ is arbitrary, $\chi\in \Omega_{A,\infty}$. Conversely, suppose that $\# A < \infty$ and that $0 < \# A_+, \# A_-$. Take $e \in A_+$ and $f \in A_-$ and set $w_{k}:=b^{k} f e e e \dotsm,\ k\in \mathbb{N}$. Then $\chi_{w_{k}}\in \Omega_{A,\ \infty}$. We claim that $\chi_{w_{k}}$ converges to $\chi_{b^{\infty}}$. Firstly, $\chi_{w_{k}}(b^{i})=1$ for all $i\in \mathbb{N}$ and thus $\lim_{k}\ \chi_{w_{k}}(b^{i})=1$ for all $i\in \mathbb{N}$. Furthermore, if $\limsup_{k}\ \chi_{w_{k}}(x)=1$ for some $x\in P$ with $\theta(x)>0$, then there exists $l\in \mathbb{N}$ such that $xP \subseteq b^{l} f P$ and that $\limsup_{k}\ \chi_{w_{k}}(b^{l} f)=1$. At the same time, $b^{k} f e^{n} \notin b^{l} f P$ for all $k>l$ and all $n\in \mathbb{N}$, contradicting $\limsup_{k}\ \chi_{w_{k}}(b^{l} f)=1$. So $\lim_{k}\ \chi_{w_{k}}(x)=0$ for all $x\in P$ with $\theta(x)>0$. This proves our claim, which implies that $\Omega_{A,\ \infty}$ is not closed.
\eproof
\pars

\bproof[Proof of Theorem~\ref{thm:clinvsub_GBS}]
Using $\Omega_{\infty} = \Omega_{b,\infty} \cup \Omega_{A,\infty}$, Theorem~\ref{thm:clinvsub_GBS} follows from Lemma~\ref{lem:Omegainf-vs-Omega}, Corollary~\ref{cor:Omegab=..._GBS} and Propositions~\ref{prop:bOmega_GBS}, \ref{prop:OmegaA_GBS}.
\eproof

\section{Topological freeness}

Assume that $G_v \subseteq (\Rz,+)$ for all $v \in V$, and $\# V > 1$ or $\# A > 0$, and that conditions (LCM) from Definition~\ref{def:LCM} and (D) from Definition~\ref{def:D} are satisfied. We set out to determine, for every closed invariant subspace $X \subseteq \Omega$, when $G \curvearrowright X$ is topologically free. To state our main result, we need the following terminology. Assume that for all $e \in T$, $P_e \neq \gekl{\beps}$. Given $e \in A$, let $v = o(e)$ and $w = t(e)$. Let $b_v$ be the generator of $P_v$ and $b_w$ the generator of $P_w$. Let $m_e, n_e \in \Zz_{\geq 1}$ be such that $(\cdot)^{\bar{e}}: \: P_{\bar{e}} \into P_v$ is given by $z \ma n_e z$ and $(\cdot)^e: \: P_{\bar{e}} \into P_w$ is given by $z \ma \pm m_e z$. Then we have $b_v^{n_e} e = e b_w^{\pm m_e}$ in $G$. Moreover, as $P_e \neq \gekl{\beps}$ for all $e \in T$, we have $\langle b_v^{n_e} \rangle \cap \langle b_w^{m_e} \rangle = \langle b_v^{l_e n_e} \rangle = \langle b_w^{k_e m_e} \rangle$ for some $k_e, l_e \in \Zz_{\geq 1}$.

\btheo
\label{thm:topfree}
\begin{enumerate}
\item[(i)] $G \curvearrowright \Omega$ is always topologically free.
\item[(ii)] $G \curvearrowright \partial \Omega$ is topologically free if and only if one of the following holds:
\begin{enumerate}
\item[(ii$_0$)] There exists $e \in T$ with $P_e = \gekl{\beps}$,
\item[(ii$_1$)] For all $e \in T$, $P_e \neq \gekl{\beps}$, $\# A > 0$ and there exists $e \in A$ with $k_e \not\mid l_e$,
\item[(ii$_2$)] For all $e \in T$, $P_e \neq \gekl{\beps}$, $\# A > 0$, for all $e \in A$, $k_e \mid l_e$ and $\big( \bigcap_{e \in A} \langle b_{o(e)}^{k_e n_e} \rangle \big) \cap \big( \bigcap_{v \in V} G_v \big) = \gekl{\beps}$.
\end{enumerate}
\item[(iii)] If $1 < \# V < \infty$, $P_v \cong \Zzg$ and $\# A_+ < \infty$, then $G \curvearrowright \Omega_{\infty}$ is topologically free if and only if one of the following holds:
\begin{enumerate}
\item[(a)] There exists $e \in T$ with $P_e = \gekl{\beps}$,
\item[(b)] $\# A_+ > 0$,
\item[(c)] $\# V > 2$,
\item[(d)] $\# A_+ = 0$, $\# V = 2$, and for the unique $e \in T$ and the embeddings $(\cdot)^{\bar{e}}: \: P_e \cong \Zzg \into \Zzg \cong P_{o(e)}, \ z \ma kz$ and $(\cdot)^e: \: P_e \cong \Zzg \into \Zzg \cong P_{t(e)}, \ z \ma lz$, we have $(k,l) \neq (2,2)$.
\end{enumerate}
\item[(iv)] If $\# V = 1$, $\# A < \infty$, $G \curvearrowright \Omega_{b,\infty}$ and $G \curvearrowright \Omega_{\infty}$ are not topologically free.
\item[(v)] If $\# V = 1$, $0 < \# A_+ < \infty$, $\# A_- = 0$, $G \curvearrowright \Omega_{A,\infty}$ is topologically free if and only if $\# A_+ \geq 2$, or $\# A_+ = 1$, $A_+ = \gekl{e}$ and $m_e \neq 1$.
\item[(vi)] If $\# V = 1$, $0 < \# A_+ < \infty$, $\# A_- = \infty$, $G \curvearrowright \Omega_{\infty}$ is topologically free if and only if $\# A_+ \geq 2$, or $\# A_+ = 1$, $A_+ = \gekl{e}$ and $m_e \neq 1$.
\end{enumerate}
\etheo
\nopar

Note that for one vertex GBS monoids given by the presentation
$$
 P = \spkl{\gekl{b} \cup A \ \vert \ b^{n_e} e = e b^{m_e} \ \forall \ e \in A_+, \, b^{n_e} e b^{m_e} = e \ \forall \ e \in A_-}^+,
$$
Theorem~\ref{thm:topfree}~(ii) says that $G \curvearrowright \partial \Omega$ is topologically free if and only if either there exists $e \in A$ with $n_e \not\mid m_e$, or $n_e \mid m_e$ for all $e \in A$ and $ \gcd(\menge{n_e}{e \in A}) = \infty$.
\pars

In order to prove the theorem, we first determine when $G \curvearrowright \partial \Omega$ is topologically free. It follows immediately from Theorem~\ref{thm:clinvsub_gen} that $G \curvearrowright \partial \Omega$ is topologically free if there exists $e \in T$ with $P_e = \gekl{\beps}$, and that $G \curvearrowright \partial \Omega$ is not topologically free if $G = G_T$ (i.e., $\# A = 0$). We can thus focus on the case where $P_e \neq \gekl{\beps}$ for all $e \in T$. Recall that $G^c = \menge{g \in G}{g pP \cap pP \neq \emptyset \ \forall \ p \in P}$.
\bprop
\label{prop:topfree_bOmega}
Assume that $P_e \neq \gekl{\beps}$ for all $e \in T$ and that $\# A > 0$. Then $G^c = \gekl{\beps}$ if there exists  $e \in A$ with $k_e \not\mid l_e$. If for all $e \in A$, $k_e \mid l_e$, then $G^c = \big( \bigcap_{e \in A} \langle b_{t(e)}^{k_e n_e} \rangle \big) \cap \big( \bigcap_{v \in V} G_v \big)$.
\eprop
\nopar

\bproof
We first treat the one vertex GBS case, i.e., $P = \spkl{\gekl{b} \cup A \ \vert \ b^{n_e} e = e b^{m_e} \ \forall \ e \in A_+, \, b^{n_e} e b^{m_e} = e \ \forall \ e \in A_-}^+$. It is clear that $G^c \subseteq \spkl{b}$. Given $g \in G^c$ and $e \in A$, $g eeP \cap eeP \neq \emptyset$ implies that $g \in \spkl{b^{n_e}}$, say $g = b^{\lambda n_e}$, and that there exist $x, y \in P$ with $g eex = eey$, which implies $e b^{\pm \lambda m_e} ex = eey$ and thus $b^{\pm \lambda m_e} \in \spkl{b^{n_e}}$, i.e., $\lambda m_e = \lambda' n_e$. Continuing this way, we obtain $\lambda' m_e = \lambda'' n_e, \dotsc, \lambda^{(i-1)} m_e = \lambda^{(i)} n_e$ for some $\lambda^{(\bullet)} \in \Zz$. Hence, for all $i \in \Nz$, $\lambda^{(i)} = \frac{m_e}{n_e} \lambda^{(i-1)} = \dotso = (\frac{m_e}{n_e})^i \lambda$. Thus $(\frac{m_e}{n_e})^i \lambda \in \Zz$ for all $i \in \Nz$, which implies $n_e \mid m_e$ unless $\lambda = 0$. This shows that $G^c = \gekl{\beps}$ if there exists $e \in A$ with $n_e \not\mid m_e$.
\pars

In the general case, the same argument shows that $G^c = \gekl{\beps}$ if there exists  $e \in A$ with $k_e \not\mid l_e$. Now assume that for all $e \in A$, $k_e \mid l_e$. Take $g \in G^c$ and $e \in A$. $g eeP \cap eeP \neq \emptyset$ implies that $g \in \langle b_{o(e)}^{k_e n_e} \rangle$, say $g = b_{o(e)}^{\kappa}$. Take $v \notin \gekl{o(e),t(e)}$, suppose that $[v,o(e)]$ starts with $d_o$, $[v,t(e)]$ starts with $d_t$, and choose $\gamma \in P_v$ with $\gamma \notin P_{d_o}^{\bar{d}_o}$, $\gamma \notin P_{d_t}^{\bar{d}_t}$. Then $g \gamma eP \cap \gamma eP \neq \emptyset$ implies that $\gamma^{-1} g \gamma = b_{o(e)}^{\lambda}$, so that $\gamma^{-1} b_{o(e)}^{\kappa} \gamma = b_{o(e)}^{\kappa}$. Comparing normal forms, we must have $g = b_{o(e)}^{\kappa} \in G_v$. Hence $g \in \bigcap_{v \in V} G_v$. This shows that $G^c \subseteq \big( \bigcap_{e \in A} \langle b_{t(e)}^{k_e n_e} \rangle \big) \cap \big( \bigcap_{v \in V} G_v \big)$. The reverse inclusion is straightforward.
\eproof
\pars

In the following, given $\chi \in \Omega$, we write $\Stab(\chi) \defeq \menge{g \in G}{g.\chi = \chi}$.
\bproof[Proof of Theorem~\ref{thm:topfree}]
(i) follows from $P^* = \gekl{\beps}$ (see \cite[Theorem~5.7.2]{CELY}). (ii) follows from Proposition~\ref{prop:topfree_bOmega} and Theorem~\ref{thm:Gc}. 
\pari

Let us prove (iii). In case (a), suppose that $e \in T$ satisfies $P_e = \gekl{\beps}$. Let $v=o(e)$, $w=t(e)$, $P_v = \spkl{\alpha}^+$ and $P_w = \spkl{\beta}^+$. Set $X \defeq \alpha^{k_1} \beta^{k_2} \alpha^{k_3} \beta^{k_4} \dotsm$, where $(k_i)_i$ is aperiodic. Then $\Stab(\chi_X) = \gekl{\beps}$. In case (b), take $e\in A_{+}$, $v \in V$, $v \neq o(e)$ and let $\alpha\in P_v$ be the generator, set $X:=\alpha^{k_{1}} e \alpha^{k_{2}} e \cdots$, where $(k_{i})_i$ is an aperiodic sequence in $\{0,\ 1\}$. Then $\Stab(\chi_X) = \gekl{\beps}$. For (c), take $u,\ v,\ w\in V$ and let $\alpha\in P_{u},\ \beta\in P_{v},\ \gamma \in P_{w}$ be the generators, set $X:=\alpha \beta \gamma^{k_{1}}\alpha \beta \gamma^{k_{2}}\cdots$, where $(k_{i})_i$ is an aperiodic sequence in $\{0,\ 1\}$. Then $\Stab(\chi_X) = \gekl{\beps}$. In case (d), let $\alpha\in P_{o(e)},\ \beta\in P_{t(e)}$ be the generators. If $k>2$, set $X:=\alpha^{k_{1}} \beta \alpha^{k_{2}} \beta \cdots$, where $(k_{i})_i$ is an aperiodic sequence in $\{1,\ 2\}$. Then $\Stab(\chi_X) = \gekl{\beps}$. The case $l > 2$ is analogous. If $k=l=2$, $\Omega_{\infty}\setminus \Omega_{\textbf{b},\ \infty}$ is a single orbit containing $\chi_{Y}$ with $Y=\alpha \beta \alpha \beta \cdots$, and $\Stab(\chi_{Y})\neq \{\beps\}$, so $G\curvearrowright \Omega_{\infty}$ is not topologically free.

(iv) holds because $\Omega_{b,\infty} = G. \chi_{b^{\infty}} \cup \partial \Omega$ and $\Stab(\chi_{b^{\infty}}) \neq \gekl{\beps}$. 

Let us prove (v). If $\# A_+ \geq 2$, take $e_1,\ e_2 \in A_+$ with $e_1 \neq e_2$ and let $\chi_{w}\in \Omega_{A,\ \infty}$ with $w= e_{j_{1}}e_{j_{2}}e_{j_{3}}\cdots$, where $j_{\mu}$ is an aperiodic sequence in $\{1,\ 2\}$. We then have $\chi_{w}\notin \partial\Omega$ and $\Stab(\chi_{w})=\{\beps\}$. If $\# A_+ = 1$, say $A_+ = \gekl{e}$, let $w=b^{i_{0}} e b^{i_{1}} e b^{i_{2}} e \cdots$ be such that $i_{\bullet}$ is an aperiodic sequence in $\{0,\ 1\}$. Then $\Stab(\chi_w) = \gekl{\beps}$ if $m_e \neq 1$. If $m_e = 1$, then $\Omega_{A,\ \infty}\setminus \partial \Omega = G.\chi_{a^{\infty}}$ and $\Stab(\chi_{a^{\infty}}) \neq\{\beps\}$. Hence $G\curvearrowright \Omega_{A,\ \infty}$ is not topologically free. The proof of (vi) is analogous.
\eproof
\pars

We deduce the following immediate consequences.
\bcor
\label{cor:tfAll}
$G \curvearrowright X$ is topologically free for every closed invariant subspace $X$ of $\Omega$ if and only if one of the following holds:
\nopar

\begin{enumerate}
\item[(i)] There exists $e \in T$ with $P_e = \gekl{\beps}$,
\item[(ii)] For all $e \in T$, $P_e \neq \gekl{\beps}$, $\# V > 1$, $\# A > 0$, and one of the following holds:
\begin{enumerate}
\item[(ii$_1$)] There exists $e \in A$ with $k_e \not\mid l_e$,
\item[(ii$_2$)] For all $e \in A$, $k_e \mid l_e$ and $\big( \bigcap_{e \in A} \langle b_{o(e)}^{k_e n_e} \rangle \big) \cap \big( \bigcap_{v \in V} G_v \big) = \gekl{\beps}$.
\end{enumerate}
\item[(iii)] $\# V = 1$, $\# A = \infty$, $\# A_+ \in \gekl{0, \infty}$, and (ii$_1$) or (ii$_2$) holds,
\item[(iv)] $\# V = 1$, $\# A = \infty$, $\# A_+ < \infty$, and (ii$_1$) or (ii$_2$) holds, and either $\# A_+ \geq 2$ or $\# A_+ = 1$, $A_+ = \gekl{e}$ and $m_e \neq 1$.
\end{enumerate}
\ecor
\pars

\bcor
\label{cor:ideals=clinvsub}
If one of (i) -- (iv) in Corollary~\ref{cor:tfAll} is satisfied, then the assignment $X \ma C^*_r(G \ltimes (\Omega \setminus X))$ is a one-to-one correspondence between closed invariant subspaces of $X$ and ideals of $C^*_{\lambda}(P) \cong C^*_r(G \ltimes \Omega)$.
\ecor
\nopar

\bproof
It follows from \cite{Gue} that $G$ is exact. Hence the claim follows from Corollar~\ref{cor:tfAll} and \cite[Theorem~A]{BL20}.
\eproof
\pars

\bremark
For one vertex GBS monoids $P$ for which Corollary~\ref{cor:ideals=clinvsub} does not apply, the primitive ideal space of $C^*_{\lambda}(P)$ has been completely described in \cite{Chen_PhD}.
\eremark

\section{Amenability and nuclearity}

Let $P$ be as in \S~\ref{ss:presentation} and assume that condition (LCM) is satisfied. For the nuclearity of the reduced semigroup $C^{*}$-algebra $C^{*}_{\lambda}(P)$, we have the following theorem, whose proof is based on the idea of \an{controlled maps} as in \cite{LR96,CL,CHR,HRT,HNSY}.
\btheo
\label{thm:NucP-PT}
$C^*_{\lambda}(P)$ is nuclear if $C^{*}_{\lambda}(P_{T})$ is nuclear.
\etheo
\nopar

\bproof
By Proposition~\ref{prop:LCM}, $P$ is right LCM. We have the following expression:
$$C^{*}_{\lambda}(P)=\overline{\text{span}\{\lambda_{p}\lambda^{*}_{q},\ p,\ q\in P\}}.$$
\pars

Let $\theta:\ P\rightarrow \mathbb{N}$ be a semigroup homomorphism such that $\theta(e)=1$ for all $e\in A$ and that $\theta(x)=0$ for all $x\in P_{T}$. Define a unitary $u_{z}, z\in \mathbb{T}$, on $\ell_{2}(P)$ by
$$u_{z}(\delta_{x})=z^{\theta(x)}\delta_{x},\ x\in P.$$
Then $\text{Ad}(u_{z})$ is a $*$-isomorphism of $C^{*}_{\lambda}(P)$. Furthermore, we have $$\text{Ad}(u_{z})(\lambda_{p}\lambda^{*}_{q})=z^{-\theta(p)+\theta(q)}\lambda_{p}\lambda^{*}_{q}.$$
Define an action $\alpha$ of $\mathbb{T}$ on $C^{*}_{\lambda}(P)$ by $\alpha(z):=\text{Ad}(u_{\Bar{z}})$, then the $k$th spectral subspace for $\alpha$ is given by:
$$B_{k}=\overline{\text{span}\{\lambda_{p}\lambda^{*}_{q},\ \theta(p)-\theta(q)=k,\ p,\ q\in P\}},\ k\in \mathbb{Z}.$$
It is easy to see that $B_{k}=B^{k}_{1},\ k\in \mathbb{Z}^{*}$, which implies, by \cite[Proposition (4.8)]{Exe94}, that the action $\alpha$ is semi-saturated. If $\alpha$ is regular, by \cite[Theorem 4.21]{Exe94}, $C^{*}_{\lambda}(P)$ is isomorphic to a partial crossed product of $B_{0}$ by a partial automorphism. In this case, $C^{*}_{\lambda}(P)$ is nuclear if and only if $B_{0}$ is nuclear.

If $\alpha$ is not regular, then tensoring it by the trivial circle action on $\mathcal{K}$, we get a stable action $\alpha'$. Furthermore, $\alpha'$ is still semi-saturated. This implies that $\alpha'$ is regular by \cite[Corollary 4.5]{Exe94}. Again by \cite[Theorem 4.21]{Exe94}, $C^{*}_{\lambda}(P)\otimes \mathcal{K}$ is isomorphic to a partial crossed product of $B_{0}\otimes \mathcal{K}$ by a partial automorphism. In this case, $C^{*}_{\lambda}(P)$ is nuclear if and only if $B_{0}\otimes \mathcal{K}$ is nuclear. And the latter holds if and only if $B_{0}$ is nuclear. Therefore, $C^{*}_{\lambda}(P)$ is nuclear if and only if $B_{0}$ is nuclear.

For $p,\ q\in P$, let
\begin{eqnarray*}
 p &=& h_{0}e_{1}h_{1}e_{2}\cdots h_{k-1}e_{k}h_{k}, h_{i-1}\in P_{T},\ e_{i}\in A,\ 1\leq i \leq k,\ h_{k}\in G^{e_{k}}_{e_{k}}P_{T},\\
 q &=& h'_{0}e'_{1}h'_{1}e'_{2}\cdots h'_{l-1}e'_{l}h'_{l}, h'_{j-1}\in P_{T},\ e'_{j}\in A,\ 1\leq j \leq l,\ h'_{l}\in G^{e'_{l}}_{e'_{l}}P_{T}
\end{eqnarray*}
be the compact forms. We say $p\sim q$ if
$$h_{0}e_{1}h_{1}e_{2}\cdots h_{k-1}e_{k}G^{e_{k}}_{e_{k}}=h'_{0}e'_{1}h'_{1}e'_{2}\cdots h'_{l-1}e'_{l}G^{e'_{l}}_{e'_{l}}.$$
Alternatively, $p\sim q$ if $k=l$, $e_{i}=e'_{i}$ for all $1\leq i \leq k$ and there exists $x\in G^{e_{k}}_{_{k}}$ such that
$$h_{0}e_{1}h_{1}e_{2}\cdots h_{k-1}e_{k}=h'_{0}e'_{1}h'_{1}e'_{2}\cdots h'_{l-1}e'_{l}x.$$
It is easy to check that $\sim$ is a well-defined equivalence relation in $P$.

For $p\in P$ with a compact form as above, define $\bar{p} \defeq h_{0}e_{1}h_{1}e_{2}\cdots h_{k-1}e_{k}$. Then $\bar{p}$ is unique up to the equivalence relation $\sim$. Moreover, for all $p,\ q\in P$, $p\sim q$ if and only if $\Bar{p}\sim \Bar{q}$.

Let $P_{l}:=\{p\in P,\ \theta(p)=l\},\ l\in \mathbb{N}$ and let $B_{0,\ l}:=\overline{\text{span}\{\lambda_{p}\lambda^{*}_{q},\ \ p,\ q\in P_{l}\}}$. Then $B_{0,\ l}$, restricted on $\ell_{2}(\cup_{k<l}P_{k})$, is 0. Therefore, we can regard $B_{0,\ l}$ as a $C^{*}$-algebra on the Hilbert space $\ell_{2}(\cup_{k\geq l}P_{k})$.

When $A_{-}=\emptyset$, $\lambda_{p}\lambda^{*}_{q}$ is of the form $\lambda_{\Bar{p}}\lambda_{h}\lambda_{h'}^{*}\lambda^{*}_{\Bar{q}},\ h,\ h'\in P_{T}$. Furthermore, we have in $B_{0,\ l}$,
\begin{equation*}
   \lambda_{\Bar{p}_{1}}\lambda_{h_{1}}\lambda_{h_{1}'}^{*}\lambda^{*}_{\Bar{q}_{1}}\cdot\lambda_{\Bar{p}_{2}}\lambda_{h_{2}}\lambda^{*}_{h'_{2}}\lambda^{*}_{\Bar{q}_{2}}=
    \begin{cases}
    \lambda_{\Bar{p}_{1}}\lambda_{h_{1}}\lambda_{h_{1}'}^{*}\lambda^{*}_{x}\lambda_{h_{2}}\lambda_{h_{2}'}^{*}\lambda^{*}_{\Bar{q}_{2}},&  \Bar{q}_{1}=\Bar{p}_{2}x,\ x\in P_{T},\\
    \lambda_{\Bar{p}_{1}}\lambda_{h_{1}}\lambda_{h_{1}'}^{*}\lambda_{x}\lambda_{h_{2}}\lambda_{h_{2}'}^{*}\lambda^{*}_{\Bar{q}_{2}},&  \Bar{q}_{1}x=\Bar{p}_{2},\ x\in P_{T},\\
    0,& \text{otherwise}.\\
    \end{cases}
\end{equation*}
\newline
Let $H_{l} \defeq \ell_{2}( \{\Bar{p}: \: \theta(p)=l\} )$ and define a linear map
$$V:\ H_{l}\otimes \ell_{2}(P)\ \rightarrow \ell_{2}(\cup_{k\geq l}P_{k})$$
by sending $\delta_{\Bar{p}}\otimes \delta_{x}$ to $\delta_{\Bar{p}x}$, then $V$ is a unitary. Let $\cK_{l} \defeq \cK(H_l)$. Then the map
$$\varphi:\ B_{0,\ l}\rightarrow \cK_{l}\otimes \mathcal{L}(\ell_{2}(P)),\ T\mapsto V^{*}TV$$
is an injective $*$-homomorphism. Furthermore, it maps
$\lambda_{\Bar{p}}\lambda_{h}\lambda_{h'}^{*}\lambda^{*}_{\Bar{q}}$ to $\varepsilon_{\Bar{p},\ \Bar{q}}\otimes \lambda_{h}\lambda_{h'}^{*}$ (where $\varepsilon_{\bar{p},\bar{q}}$ denotes the standard matrix units) and hence $\varphi(B_{0,\ l})= \cK_{l}\otimes C^{*}(\lambda(P_{T}))$. It is straightforward to see that $C^{*}(\lambda(P_{T}))\cong C^{*}_{\lambda}(P_{T})$. Since $C^{*}_{\lambda}(P_{T})$ is nuclear, so is $B_{0,\ l}$.

When $A_{-}\neq \emptyset$, for $p\in P$ with compact form $p=h_{0}e_{1}h_{1}e_{2}\cdots h_{k-1}e_{k}h_{k}$, $h_{i-1}\in P_{T}$, $e_{i}\in A$, $1\leq i \leq k$, $h_{k}\in G^{e_{k}}_{e_{k}}P_{T}$, define
$X_{\bar{p}} \defeq \{x\in P^{e_{k}}_{e_{k}},\ \bar{p}x^{-1}\in P\}$. If $X_{p}\neq \{\beps\}$, then there must exist a sequence $(x^{(n)}_{p})_{n\in \mathbb{N}}\subseteq X_{p}$ with $x^{(n)}_{p}\prec x^{(n+1)}_{p}$ such that for all $x\in X_{p}$, $x\prec x^{(n)}_{p}$ for some $n\in \mathbb{N}$ since every group $G_{v},\ v\in V$ is countable and totally ordered. For each $n\in \mathbb{N}$, let
\begin{equation*}
    \bar{p}^{(n)}:=
    \begin{cases}
    \Bar{p},&  X_{p}=\{\epsilon\},\\
    \Bar{p}(x^{(n)}_{p})^{-1},& X_{p}\neq \{\epsilon\}.\\
    \end{cases}
\end{equation*}
Define
$$B^{(n)}_{0,\ l}:=\overline{\text{span}\{\lambda_{\bar{p}^{(n)}}\lambda_{h}\lambda^{*}_{h'}\lambda^{*}_{\bar{q}^{(n)}},\ p,\ q\in P_{l},\ h,\ h'\in P_{T}\}}$$
and define $\cK_l$ as before. Similarly as in the case when $A_{-}=\emptyset$, we obtain $B^{(n)}_{0,\ l}\cong \cK_{l}\otimes C^{*}(\lambda(P_{T}))$, which means that $B^{(n)}_{0,\ l}$ is nuclear. Noting that $B^{(n)}_{0,\ l} \subseteq B^{(n+1)}_{0,\ l}$, we conclude that $B_{0,\ l}=\overline{\cup_{n\in \mathbb{N}}B^{(n)}_{0,\ l}}$ is nuclear as an inductive limit of nuclear $C^{*}$-algebras. Define $B_{0,\ \leq l} \defeq \sum_{0\leq k\leq l} B_{0,\ k}$. Then $B_{0,\ l},\ l\geq 1$ is an ideal in $B_{0,\ \leq l}$ and the corresponding quotient is a quotient of $B_{0,\ \leq l-1}$. Since quotients and extensions of $C^{*}$-algebras preserve nuclearity, we get, by induction, that $B_{0,\ \leq l}$ is nuclear. As an inductive limit of nuclear $C^{*}$-algebras, $B_{0}=\overline{\cup_{l\geq 0}B_{0,\ \leq l}}$ is nuclear. Therefore, $C^{*}_{\lambda}(P)$ is nuclear.
\eproof
\pars

Let us prove the converse.
\bprop
If $C^*_{\lambda}(P)$ is nuclear, then $C^{*}_{\lambda}(P_{T})$ is nuclear.
\eprop
\nopar

\bproof
Let $\Omega_T \defeq \Omega_{P_T}$. There is a canonical map $\Omega_T \to \Omega$ sending $\chi_w$ to $\chi_w$, where $w$ is a word in $\menge{P_v}{v \in V}$. This map is continuous and a homeomorphism onto its image because for all $p, q_i \in P_T$, we have $\Omega(p;q_i) \cap \Omega_T = \Omega_T(p;q_i)$ by Proposition~\ref{prop:LCMP,PT}. This map induces a groupoid embedding $G_T \ltimes \Omega_T \into G \ltimes \Omega$ which is a homeomorphism onto its image, which is closed. If $C^*_{\lambda}(P)$ is nuclear, then $G \ltimes \Omega$ is amenable, so that $G_T \ltimes \Omega_T$ is also amenable as a closed subgroupoid (see \cite[Proposition~5.1.1]{AR}).
\eproof
\pars

Now assume that $G_v \subseteq (\Rz,+)$ for all $v \in V$, and $\# V > 1$ or $\# A > 0$.
\btheo
\label{thm:CPTnuc}
The following are equivalent:
\nopar

\begin{enumerate}
\item[(i)] $C^*_{\lambda}(P)$ is nuclear,
\item[(ii)] $C^*_{\lambda}(P_T)$ is nuclear,
\item[(iii)] For all $T' \subseteq T$ with $P_e \neq \gekl{\beps}$ for all $e \in T'$, either $T'$ consists of a single vertex or $T'$ consists of exactly two vertices $v, w$ and one pair of edges $e, \bar{e}$ with $o(e) = v$, $t(e) = w$, such that $P_v \cong \Zzg$, $P_w \cong \Zzg$, and the embeddings $(\cdot)^{\bar{e}}$, $(\cdot)^e$ are both given by $\Zzg \to \Zzg, \, z \ma 2z$.
\end{enumerate}
\etheo
\pars

For the proof, we need some preparations. We start with the following generalization of a construction from \cite{Hoc}.
\blemma
\label{lem:submonoid-amenable}
Let $P$ be a submonoid of a group $G$ and $Q$ a submonoid of a group $H$. Assume that $Q^* = \gekl{\beps}$. Then $P * Q$ embeds into $(\bigoplus_H G) \rtimes H$.
\elemma
\nopar

\bproof
Given $h \in H$ and $g \in G$, let $\delta(h,g) \in \bigoplus_H G$ be given by $\delta(h,g)(x) = \beps$ if $x \neq h$ and $\delta(h,g)(h) = g$. By definition of the $H$-action on $\bigoplus_H G$, we have $h.\delta(b,a) = \delta(hb,a)$. Now define 
$$
 \varphi: \: P * Q \into (\bigoplus_H G) \rtimes H, \, P \ni p \ma (\delta(\beps,p),\beps), \, Q \ni q \ma (\beps,q).
$$
A straightforward proof by induction shows that
$$
 \varphi(p_1 q_1 \dotsm p_n q_n) = (\delta(\beps,p_1) \delta(q_1,p_2) \delta(q_1 q_2,p_3) \dotsm \delta(q_1 \dotsm q_{n-1},p_n), q_1 \dotsm q_n)
$$
for all $p_i \in P$, $q_i \in Q$. Now we show that $\varphi$ is injective, i.e., $\varphi(p_1 q_1 \dotsm p_m q_m) = \varphi(r_1 s_1 \dotsm r_n s_n)$ implies $p_1 q_1 \dotsm p_m q_m = r_1 s_1 \dotsm r_n s_n$ for all $p_i, r_j \in P \setminus \gekl{\beps}$ and $q_i, s_j \in Q \setminus \gekl{\beps}$. We proceed inductively on $\max(m,n)$. The induction start $\max(m,n) = 1$ is clear. For the induction step, suppose that $\varphi(p_1 q_1 \dotsm p_m q_m) = \varphi(r_1 s_1 \dotsm r_n s_n) = (f,h)$. Then $p_1 = f(\beps) = r_1$ because $Q^* = \gekl{\beps}$. Moreover, both $q_1$ and $s_1$ can be characterized as the minimal element (with respect to $\prec$) $x$ of $Q \subseteq H$ such that $x \neq \beps$ and $f(x) \neq \beps$. Hence it follows that $q_1 = s_1$. We deduce that $\varphi(p_1 q_1 \dotsm p_m q_m) = \varphi(r_1 s_1 \dotsm r_n s_n) = \varphi(r_1 s_1) \varphi(r_2 s_2 \dotsm r_n s_n) = \varphi(p_1 q_1) \varphi(r_2 s_2 \dotsm r_n s_n)$ and thus $\varphi(p_2 q_2 \dotsm p_m q_m) = \varphi(r_2 s_2 \dotsm r_n s_n)$. Induction hypothesis implies $p_2 q_2 \dotsm p_m q_m = r_2 s_2 \dotsm r_n s_n$.
\eproof
\pars

\bcor
\label{cor:submonoid-amenable}
Free products of monoids with no non-trivial invertible elements and which embed into amenable groups again embed into amenable groups.
\ecor

\bproof[Proof of Theorem~\ref{thm:CPTnuc}] 
First of all, for a tree $T'$, $G_{T'}$ is amenable if and only if either $T'$ consists of a single vertex or $T'$ consists of exactly two vertices $v, w$ and one pair of edges $e, \bar{e}$ with $o(e) = v$, $t(e) = w$, such that $G_v \cong \Zz$, $G_w \cong \Zz$, and the embeddings $(\cdot)^{\bar{e}}$, $(\cdot)^e$ are both given by $\Zz \to \Zz, \, z \ma 2z$. 

Now assume that $C^*_{\lambda}(P_T)$ is nuclear. Then for all subtrees $T' \subseteq T$ with $P_e \neq \gekl{\beps}$ for all $e \in T'$, $G_{T'}$ must be amenable, because they are the stabilizer groups of $\infty_{P_{T'}}$. This shows \an{$\Rarr$}. Conversely, if all $T' \subseteq T$ with $P_e \neq \gekl{\beps}$ for all $e \in T'$ are as in the statement of Theorem~\ref{thm:CPTnuc}, then it follows that
$P_T \cong \Conv_i P_{T_i}$,
where for each $i$, $T_i$ consists of a single vertex or $T'$ consists of exactly two vertices $v, w$ and one pair of edges $e, \bar{e}$ with $o(e) = v$, $t(e) = w$, such that $P_v \cong \Zzg$, $P_w \cong \Zzg$, and the embeddings $(\cdot)^{\bar{e}}$, $(\cdot)^e$ are both given by $\Zzg \to \Zzg, \, z \ma 2z$. Hence for all $i$, $P_{T_i}$ is a submonoid of the amenable group $G_{T_i}$. Corollary~\ref{cor:submonoid-amenable} implies that $P_T$ embeds into an amenable group, and thus $C^*_{\lambda}(P_T)$ is nuclear by \cite[Corollary 3.16]{Li17}.
\eproof

\bremark
The proof of Theorem~\ref{thm:CPTnuc} shows that $C^*_{\lambda}(P_T)$ is nuclear if and only if $P_T$ embeds into an amenable group.
\eremark

\section{K-theory}

We now assume that $G_v \subseteq (\Rz,+)$ for all $v \in V$, and $\# V > 1$ or $\# A > 0$, and that conditions (LCM) and (D) are satisfied. Our goal is to compute K-theory for $C^*_r(G \ltimes X)$ for all closed invariant subspaces $X \subseteq \Omega$. In the following, given a unital C*-algebra $C$, we set $K_*(C) \defeq (K_0(C),[1_C]_0,K_1(C))$. For $e \in A$, define $\sgn(e) \defeq +1$ if $e \in A_+$ and $\sgn(e) \defeq -1$ if $e \in A_-$. In the one vertex GBS case, i.e., $\# V = 1$ and $P_v \cong \Zzg$ for all $v \in V$, recall that we have the presentation $P = \spkl{\gekl{b} \cup A \ \vert \ b^{n_e} e = e b^{m_e} \ \forall \ e \in A_+, \, b^{n_e} e b^{m_e} = e \ \forall \ e \in A_-}^+$.
\btheo
\label{thm:K}
\begin{enumerate}
\item[(i)] The canonical inclusion $\Cz \cdot 1 \into C^*_r(G \ltimes \Omega)$ induces $K_*(C^*_r(G \ltimes \Omega)) \cong K_*(\Cz)$.
\item[(ii)] Assume that $G = G_T$. Then $\partial \Omega = \gekl{\infty}$, and $C^*_r(G \ltimes \gekl{\infty}) \cong C^*_{\lambda}(G)$, so that $K_*(C^*_r(G \ltimes \gekl{\infty})) \cong K_*(C^*_{\lambda}(G))$.
\item[(iii)] If $P_e \neq \gekl{\beps}$ for all $e \in T$, $\# A \geq 1$ and $X = \Omega_{b,\infty}$ if $\# V = 1$, $P_v \cong \Zzg$ for all $v \in V$, $X = \Omega_{\bm{b},\infty}$ if $\# V > 1$, then the map $C^*_{\lambda}(G_T) \to C^*_r(G \ltimes X), \, \lambda_g \ma 1_{\gekl{g} \times X}$ induces $K_*(C^*_r(G \ltimes X)) \cong K_*(C^*_{\lambda}(G_T))$.
\item[(iv)] Suppose that $\# A_+ < \infty$, $\# V < \infty$ and $P_v \cong \Zzg$ for all $v \in V$.
\begin{enumerate}
\item[(iv$_1$)] If $P_e \neq \gekl{\beps}$ for all $e \in T$, then $K_*(C^*_r(G \ltimes \Omega_{\infty})) \cong (\Zz,1, \Zz)$.
\item[(iv$_2$)] If $P_e = \gekl{\beps}$ for some $e \in T$, then $K_*(C^*_r(G \ltimes \Omega_{\infty})) \cong (\Zz/N,1, \gekl{0})$, where $N = \frac{1}{2} \# \menge{e \in T}{P_e = \gekl{\beps}}$.
\end{enumerate}
\item[(v)] Suppose that $\# V = 1$, $P_v \cong \Zzg$ for all $v \in V$, $\# A < \infty$, $\# A_+ = 0$ or $\# A_- = 0$.
\begin{enumerate}
\item[(v$_1$)] If $\sum_{e \in A} n_e \neq 1$, then $K_*(C^*_r(G \ltimes \Omega_{A,\infty})) \cong (\Zz / (1 - \sum_{e \in A} n_e), 1, \Zz / (1 + \sum_{e \in A_-} m_e))$.
\item[(v$_2$)] If $\sum_{e \in A} n_e = 1$, i.e., $\# A = 1$, $A = \gekl{e}$ and $n_e = 1$, then $K_*(C^*_r(G \ltimes \Omega_{A,\infty})) \cong (\Zz, 1, \Zz \oplus \Zz / (1 + \sum_{e \in A_-} m_e))$.
\end{enumerate}
\item[(vi)] Suppose that $\# V = 1$, $P_v \cong \Zzg$ for all $v \in V$, $\# A < \infty$ and $\# A_+ > 0$.
\begin{enumerate}
\item[(vi$_1$)] Assume $\sum_{e \in A} n_e \neq 1$.
\begin{enumerate}
\item[(vi$_{1a}$)] If $\sum_{e \in A} \sgn(e) m_e) \neq 1$, then $K_*(C^*_r(G \ltimes \partial \Omega)) \cong (\Zz / (1 - \sum_{e \in A} n_e), 1, \Zz / (1 - \sum_{e \in A} \sgn(e) m_e))$.
\item[(vi$_{1b}$)] If $\sum_{e \in A} \sgn(e) m_e) = 1$, then $K_*(C^*_r(G \ltimes \partial \Omega)) \cong (\Zz / (1 - \sum_{e \in A} n_e) \oplus \Zz, (1,0), \Zz)$.
\end{enumerate}
\item[(vi$_2$)] Assume $\sum_{e \in A} n_e = 1$, i.e., $\# A = 1$, $A = A_+ = \gekl{e}$ and $n_e = 1$.
\begin{enumerate}
\item[(vi$_{2a}$)] If $m_e \neq 1$, then $K_*(C^*_r(G \ltimes \partial \Omega)) \cong (\Zz, 1, \Zz \oplus \Zz / (1 - m_e)$.
\item[(vi$_{2b}$)] If $m_e = 1$, then $K_*(C^*_r(G \ltimes \partial \Omega)) \cong (\Zz \oplus \Zz, (1,0), \Zz \oplus \Zz)$.
\end{enumerate}
\end{enumerate}
\end{enumerate}
\etheo

Let us now prove our main result about K-theory. We start with the following result, which is an immediate consequence of \cite[Corollary~1.3]{Li2021} because $G$ satisfies the Baum-Connes conjecture with coefficients by \cite{OO}.
\bprop
\label{prop:K:C*P}
The canonical inclusion $\Cz \cdot 1 \into C^*_r(G \ltimes \Omega)$ induces $K_*(C^*_r(G \ltimes \Omega)) \cong K_*(\Cz)$.
\eprop

\bprop
\label{prop:K:Omegab}
If $P_e \neq \gekl{\beps}$ for all $e \in T$, $\# A \geq 1$ and $X = \Omega_{b,\infty}$ if $\# V = 1$, $P_v \cong \Zzg$ for all $v \in V$, $X = \Omega_{\bm{b},\infty}$ if $\# V > 1$, then the map $C^*_{\lambda}(G_T) \to C^*_r(G \ltimes X), \, \lambda_g \ma 1_{\gekl{g} \times X}$ induces $K_*(C^*_r(G \ltimes X)) \cong K_*(C^*_{\lambda}(G_T))$.
\eprop
\nopar

\bproof
We claim that $\{g X\}_{g\in G}$ is a $G$-invariant regular basis for the compact open subsets of $X$, in the sense of \cite[Definition~2.12]{Li2021}. It is easy to see that $g X$ is a compact open subset of $X$ for all $g\in G$ and that $\{gX\}_{g\in G}$ is $G$-invariant. Therefore, it remains to show that $\{g X\}_{g\in G}$ is a regular basis. First of all, if $\bigcap_{1 \leq i \leq n} p_{i} X \neq \emptyset$ with $p_{i}\in P$, $1 \leq i \leq n$ and $n\in \mathbb{N}$, then we must have $\bigcap_{1 \leq i \leq n} p_{i} P \neq \emptyset$ and thus $\bigcap_{1 \leq i \leq n} p_{i}P=rP$ for some $r\in P$ because $P$ is right LCM. Therefore, $\bigcap_{1 \leq i \leq n} p_{i}X = r X$. Secondly, for every basic compact open subset $\mathcal{O}$ in $X$, there exist $p,\ p_{i},\ 1\leq i\leq n\in P$ such that $\mathcal{O}=\{\chi \in X,\ \chi(p)=1,\ \chi(p_{i})=0\}$. In this case, we have $\mathcal{O}=pX \setminus (\bigcup_{1\leq i\leq n}p_{i}X)$. Thirdly, if $pX=\bigcup_{1 \leq i \leq n} p_{i}X$ for some $p,\ p_{i},\ 1\leq i\leq n\in P$, then we must have $pP=\bigcup_{1 \leq i \leq n} p_{i}P$ and thus $pP=p_{i}P$ for some $i$ as $\bm{\epsilon} \in P$ (in other words, because $P$ satisfies independence). In this case, $pX=p_{i}X$. These observations, together with the fact that for all $g\in G$ there exists $p\in P$ such that $g X = p X$, yields our claim that $\{g X\}_{g\in G}$ is a $G$-invariant regular basis for the compact open subsets of $X$. Since $G$ satisfies the Baum-Connes conjecture with coefficients by \cite{OO}, Proposition~\ref{prop:K:Omegab} follows from \cite[Theorem~1.2]{Li2021}.
\eproof
\pars

\bprop
\label{prop:K:OmegaInf}
Suppose that $\# A_+ < \infty$, $\# V < \infty$ and $P_v \cong \Zzg$ for all $v \in V$.
\nopar

\begin{enumerate}
\item If $P_e \neq \gekl{\beps}$ for all $e \in T$, then $K_*(C^*_r(G \ltimes \Omega_{\infty})) \cong (\Zz,1, \Zz)$.
\item If $P_e = \gekl{\beps}$ for some $e \in T$, then $K_*(C^*_r(G \ltimes \Omega_{\infty})) \cong (\Zz/N,1, \gekl{0})$, where $N = \frac{1}{2} \# \menge{e \in T}{P_e = \gekl{\beps}}$.
\end{enumerate}
\eprop
\bproof
Because of the following short exact sequence of $C^{*}$-algebras
$$0\rightarrow C^{*}_{r}(G\ltimes (\Omega\setminus \Omega_{\infty}))\rightarrow C^{*}_{r}(G\ltimes \Omega)\rightarrow C^{*}_{r}(G\ltimes \Omega_{\infty})\rightarrow 0,$$
we obtain the six term exact sequence of their K-theories as follows
\begin{equation}
\label{e:K:OmegaInf}
    \begin{tikzcd}
        K_{0}\big(C^{*}_{r}\big(G\ltimes (\Omega\setminus \Omega_{\infty})\big)\big)\arrow{r}\ &K_{0}\big(C^{*}_{r}(G\ltimes \Omega)\big)\arrow{r}\ &K_{0}\big(C^{*}_{r}(G\ltimes \Omega_{\infty})\big)\arrow{d}\\
        K_{1}\big(C^{*}_{r}(G\ltimes \Omega_{\infty})\big)\arrow{u}\ &K_{1}\big(C^{*}_{r}(G\ltimes \Omega)\big) \arrow{l}\ &K_{1}\big(C^{*}_{r}\big(G\ltimes (\Omega\setminus \Omega_{\infty})\big)\big)\arrow{l}.\\
    \end{tikzcd}
\end{equation}
Since $C^{*}_{r}(G\ltimes (\Omega\setminus \Omega_{\infty}))\cong \mathcal{K}(\ell^2 P)$ and $C^*_r(G \ltimes \Omega) \cong C^*_{\lambda}(P)$, we have
$$K_{0}(C^{*}_{r}(G\ltimes (\Omega\setminus \Omega_{\infty})))\cong \mathbb{Z}\ \text{and}\ K_{1}(C^{*}_{r}(G\ltimes (\Omega\setminus \Omega_{\infty})))\cong \gekl{0}.$$
Plugging this into \eqref{e:K:OmegaInf} and using Proposition~\ref{prop:K:C*P}, we obtain
$$
    \begin{tikzcd}
        \mathbb{Z}\arrow{r}{K_{0}(\iota)}\ &\mathbb{Z}\arrow{r}\ &K_{0}\big(C^{*}_{r}(G\ltimes \Omega_{\infty})\big)\arrow{d}\\
        K_{1}\big(C^{*}_{r}(G\ltimes \Omega_{\infty})\big)\arrow{u}\ &0 \arrow{l}\ &0\arrow{l},\\
    \end{tikzcd}
$$
where $\iota$ is the canonical inclusion. To calculate the K-theory of $C^{*}_{r}(G\ltimes \Omega_{\infty})$, we need to determine the map $K_{0}(\iota)$ from $\mathbb{Z}$ to $\mathbb{Z}$. It suffices to compute $K_{0}(\iota)([\varepsilon_{\beps,\beps}]_{0})$ (as before, $\varepsilon_{\bullet,\bullet}$ are the standard matrix units). Here we are making use of the isomorphism $C^{*}_{r}(G\ltimes (\Omega\setminus \Omega_{\infty}))\cong \mathcal{K}(\ell^2 P)$.
\pars

If $P_e \neq \gekl{\beps}$ for all $e \in T$, let $b_{v}$ be the generator of $P_{v}$. Then we have relations $b_{v}^{m_{v,w}}=b_{w}^{m_{w,\ v}}$ and $b_{v}^{m_{v,\ e}}e=eb_{w}^{m_{e,\ w}}$ for all $e \in A_+$. It is easy to see that
$$(\lambda_e \lambda_e^*) \cap (\lambda_{b_{o(e)}} \lambda_{b_{o(e)}}^*)=  \lambda_{b_{o(e)}^{m_{o(e),\ e}}e} \lambda_{b_{o(e)}^{m_{o(e),\ e}}e}^* \ \text{if}\ e\in A_{+}$$
and that
$$ \lambda_e \lambda_e^* \subseteq \lambda_{b_{o(e)}} \lambda_{b_{o(e)}}^* \ \text{if}\ e\in A_{-}.$$
For $v\in V$, denote by $v(w)$ the vertex connected to $w$ in the geodesic path $[v,\ w]\subseteq T$ for all $v\neq w\in V$. Since $\# V<\infty$ and $\# A_{+}<\infty$, we always have
$$ \varepsilon_{\beps,\beps} = 1- \big( \lambda_{b_{v}}\lambda^{*}_{b_{v}}+\sum_{v\neq w\in V}(\lambda_{b_{w}}\lambda^{*}_{b_{w}}-\lambda_{b_{v(w)}^{m_{v(w),\ w}}}\lambda^{*}_{b_{v(w)}^{m_{v(w),\ w}}})+\sum_{e\in A_{+}}(\lambda_{e}\lambda^{*}_{e}-\lambda_{b_{o(e)}^{m_{o(e),\ e}}e}\lambda^{*}_{b_{o(e)}^{m_{o(e),\ e}}e}) \big).$$
It follows that $[\iota(\varepsilon_{\beps,\beps})]_0 = 0$ in $K_*(C^{*}_{\lambda}(P))$. From the six term exact sequence, it follows that $K_*(C^*_r(G \ltimes \Omega_{\infty})) \cong (\Zz,1, \Zz)$.

If there exists $e\in T$ such that $P_{e}=\{\beps\}$, we get similarly as above that $[\iota(\varepsilon_{\beps,\beps})]_{0}=-N[1]_0$. From the six term exact sequence, we get $K_*(C^*_r(G \ltimes \Omega_{\infty})) \cong (\Zz/N,1, \gekl{0})$.
\eproof
\pars

Let us now turn to the one vertex GBS case, i.e., $P = \spkl{\gekl{b} \cup A \ \vert \ b^{n_e} e = e b^{m_e} \ \forall \ e \in A_+, \, b^{n_e} e b^{m_e} = e \ \forall \ e \in A_-}^+$.
\bprop
\label{prop:K:OmegaAInf}
Suppose that $\# V = 1$, $P_v \cong \Zzg$ for all $v \in V$, $\# A < \infty$, $\# A_+ = 0$ or $\# A_- = 0$.
\nopar

\begin{enumerate}
\item If $\sum_{e \in A} n_e \neq 1$, then $K_*(C^*_r(G \ltimes \Omega_{A,\infty})) \cong (\Zz / (1 - \sum_{e \in A} n_e), 1, \Zz / (1 + \sum_{e \in A_-} m_e))$.
\item If $\sum_{e \in A} n_e = 1$, i.e., $\# A = 1$, $A = \gekl{e}$ and $n_e = 1$, then $K_*(C^*_r(G \ltimes \Omega_{A,\infty})) \cong (\Zz, 1, \Zz \oplus \Zz / (1 + \sum_{e \in A_-} m_e))$.
\end{enumerate}
\eprop
\bproof
We have the following exact sequence of $C^{*}$-algebras,
$$0\rightarrow C^{*}_{r}\big(G\ltimes (\Omega_{\infty} \setminus \Omega_{A,\ \infty})\big)\rightarrow C^{*}_{r}\big(G\ltimes \Omega_{\infty})\rightarrow C^{*}_{r}(G\ltimes \Omega_{A,\ \infty})\rightarrow 0,$$
and the corresponding six term exact sequence of their K-theories,
\begin{equation}
\label{e:K:OmegaAInf}
    \begin{tikzcd}
        K_{0}\big(C^{*}_{r}\big(G\ltimes (\Omega_{\infty}\setminus \Omega_{A,\ \infty})\big)\big)\arrow{r}\ &K_{0}\big(C^{*}_{r}(G\ltimes \Omega_{\infty})\big)\arrow{r}\ &K_{0}\big(C^{*}_{r}(G\ltimes \Omega_{A,\ \infty})\big)\arrow{d}\\
        K_{1}\big(C^{*}_{r}(G\ltimes \Omega_{A,\ \infty})\big)\arrow{u}\ &K_{1}\big(C^{*}_{r}(G\ltimes \Omega_{\infty})\big) \arrow{l}\ &K_{1}\big(C^{*}_{r}\big(G\ltimes (\Omega_{\infty}\setminus \Omega_{A,\ \infty})\big)\big)\arrow{l}.\\
    \end{tikzcd}
\end{equation}
Since $\Omega_{\infty}\setminus \Omega_{A,\ \infty} = G. \chi_{b^{\infty}}$ and $\Stab(\chi_{b^{\infty}}) = \spkl{b} \cong \Zz$, we have
$$C^{*}_{r}\big(G\ltimes (\Omega_{\infty}\setminus \Omega_{A,\ \infty})\big)\cong \mathcal{K}\otimes C(\mathbb{T})$$
and thus
$$K_{0}\big(C^{*}_{r}\big(G\ltimes (\Omega_{\infty}\setminus \Omega_{A,\ \infty})\big)\big)\cong \mathbb{Z}\ \text{and}\ K_{1}\big(C^{*}_{r}\big(G\ltimes (\Omega_{\infty}\setminus \Omega_{A,\ \infty})\big)\big)\cong \mathbb{Z}.$$
By Proposition~\ref{prop:K:OmegaInf}, we have
$$K_{0}\big(C^{*}_{r}(G\ltimes \Omega_{\infty}) \big)\cong \mathbb{Z}\ \text{and}\ K_{1}\big(C^{*}_{r}(G\ltimes \Omega_{\infty}) \big)\cong \mathbb{Z}.$$
\pars

Plugging this into \eqref{e:K:OmegaAInf}, we obtain the six term exact sequence
$$
    \begin{tikzcd}
        \mathbb{Z}\arrow{r}{K_{0}(\iota)}\ &\mathbb{Z}\arrow{r}\ &K_{0}\big(C^{*}_{r}(G\ltimes \Omega_{A,\ \infty})\big)\arrow{d}\\
        K_{1}\big(C^{*}_{r}(G\ltimes \Omega_{A,\ \infty})\big)\arrow{u}\ &\mathbb{Z} \arrow{l}\ &\mathbb{Z}\arrow{l}{K_{1}(\iota)},\\
    \end{tikzcd}
$$
where $\iota$ is the inclusion map from $C^{*}_{r}\big(G\ltimes (\Omega_{\infty} \setminus \Omega_{A,\ \infty})\big)$ into $C^{*}_{r}\big(G\ltimes \Omega_{\infty})$.

Since $C^{*}_{r}\big(G\ltimes (\Omega_{\infty}\setminus \Omega_{A,\ \infty})\big)\cong \mathcal{K}\otimes C(\mathbb{T})$, $K_0(C^{*}_{r}\big(G\ltimes (\Omega_{\infty}\setminus \Omega_{A,\ \infty})\big))$ is generated by $[1_{(\beps,\ \chi_{b^{\infty}})}]_{0}$. Let $\pi$ be the quotient map $\pi:\ C^{*}_{r}(G\ltimes \Omega)\rightarrow C^{*}_{r}(G\ltimes \Omega_{\infty})$. Then $\pi(1_{\{\beps\}\times X}) = 1_{(\beps,\ \chi_{b^{\infty}})}$, where $X$ is determined by
$$1_{\{\beps\}\times X}=1-\sum_{e \in A,\ 0\leq j\leq n_e-1}1_{\{b^{j}e\}\times \Omega}1^{*}_{\{b^{j}e\}\times \Omega}\in C^*_r(G\ltimes \Omega).$$
In case $\# A_- = 0$, then $X=\{\chi_{b^{k}},\ k\in\ \mathbb{N}\}\cup\{\chi_{b^{\infty}}\}$, and in case $\# A_+ = 0$, then $X=\{\chi_{b^{\infty}}\}\cup\{\chi_{b^{k}},\ k\in\ \mathbb{N}\}\cup \big(\cup_{e \in A}\{\chi_{b^{k} e b^{j}},\ k\geq n_e,\ 0\leq j\leq m_e-1\}\big)$. 
In $K_{0}(C^{*}_{r}(G\ltimes \Omega))$,\ $[1_{\{\beps\}\times X}]_{0}= (1-\sum_{e \in A} n_e)[1]_0$. As we have seen in the proof of Proposition~\ref{prop:K:OmegaInf}, $\pi$ induces an isomorphism in $K_0$. Thus $[1_{(\beps,\ \chi_{b^{\infty}})}]_{0}=(1-\sum_{e \in A} n_e)[1]_0$ in $K_{0}(C^{*}_{r}(G\ltimes \Omega_{\infty}))$.

Using $C^{*}_{r}\big(G\ltimes (\Omega_{\infty}\setminus \Omega_{A,\ \infty})\big)\cong \mathcal{K}\otimes C(\mathbb{T})$, it is easy to see that $K_1(C^{*}_{r}\big(G\ltimes (\Omega_{\infty}\setminus \Omega_{A,\ \infty})\big))$ is generated by $[u]_1$, where $u \defeq 1_{(b,\ \chi_{b^{\infty}})}+1-1_{(\beps,\ \chi_{b^{\infty}})}\in C^{*}_{r}(G\ltimes \Omega_{\infty})$. Let
$$v \defeq \begin{pmatrix} 1_{\{b\}\times X} +1-1_{\{\beps\}\times X}\ &\ 1_{\{\beps\}\times (X\setminus bX)}\\ 0\ &\ 1_{\{b^{-1}\}\times bX} +1-1_{\{\beps\}\times X}\ \end{pmatrix}.$$
We have $\pi(v)=\begin{pmatrix} u\ &\ 0\\ 0\ &\ u^{*}\ \end{pmatrix}$ and $p \defeq v\begin{pmatrix} 1\ &\ 0\\ 0\ &\ 0\ \end{pmatrix}v^{*}=\begin{pmatrix} 1-1_{\{\beps\}\times (X\setminus bX)}\ &\ 0\\ 0\ &\ 0\ \end{pmatrix}$. Therefore, the index map $\delta_{1}:\ K_{1}\big(C^{*}_{r}(G\ltimes \Omega_{\infty})\big)\rightarrow K_{0}\big(C^{*}_{r}(G\ltimes (\Omega\setminus \Omega_{\infty}))\big)$ sends $[u]_1$ to
$-[1_{\{\beps\}\times (X\setminus bX)}]_{0}$. If $\# A_- = 0$, then $[1_{\{\beps\}\times (X\setminus bX)}]_{0}=1$. If $\# A_+ = 0$, then $[1_{\{\beps\}\times (X\setminus bX)}]_{0}=1+\sum_{e \in A_-}m_e$.

It is now straightforward to deduce the desired result about $K_*(C^*_r(G \ltimes \Omega_{A,\infty}))$.
\eproof
\pars

\bprop
Suppose that $\# V = 1$, $P_v \cong \Zzg$ for all $v \in V$, $\# A < \infty$ and $\# A_+ > 0$.
\nopar

\begin{enumerate}
\item[(1)] Assume $\sum_{e \in A} n_e \neq 1$.
\begin{enumerate}
\item[(1a)] If $\sum_{e \in A} \sgn(e) m_e \neq 1$, then $K_*(C^*_r(G \ltimes \partial \Omega)) \cong (\Zz / (1 - \sum_{e \in A} n_e), 1, \Zz / (1 - \sum_{e \in A} \sgn(e) m_e))$.
\item[(1b)] If $\sum_{e \in A} \sgn(e) m_e = 1$, then $K_*(C^*_r(G \ltimes \partial \Omega)) \cong (\Zz / (1 - \sum_{e \in A} n_e) \oplus \Zz, (1,0), \Zz)$.
\end{enumerate}
\item[(2)] Assume $\sum_{e \in A} n_e = 1$, i.e., $\# A = 1$, $A = A_+ = \gekl{e}$ and $n_e = 1$.
\begin{enumerate}
\item[(2a)] If $m_e \neq 1$, then $K_*(C^*_r(G \ltimes \partial \Omega)) \cong (\Zz, 1, \Zz \oplus \Zz / (1 - m_e))$.
\item[(2b)] If $m_e = 1$, then $K_*(C^*_r(G \ltimes \partial \Omega)) \cong (\Zz \oplus \Zz, (1,0), \Zz \oplus \Zz)$.
\end{enumerate}
\end{enumerate}
\eprop
\bproof
We have the following exact sequence of $C^{*}$-algebras,
$$0\rightarrow C^{*}_{r}\big(G\ltimes (\Omega_{b, \infty} \setminus \partial \Omega)\big)\rightarrow C^{*}_{r}\big(G\ltimes \Omega_{b, \infty})\rightarrow C^{*}_{r}(G\ltimes \partial \Omega)\rightarrow 0,$$
and the corresponding six term exact sequence of their K-theories,
\begin{equation}
\label{e:K:OmegaBInf}
    \begin{tikzcd}
        K_{0}\big(C^{*}_{r}\big(G\ltimes (\Omega_{b, \infty} \setminus \partial \Omega)\big)\big)\arrow{r}\ &K_{0}\big(C^{*}_{r}(G\ltimes \Omega_{b, \infty})\big)\arrow{r}\ &K_{0}\big(C^{*}_{r}(G\ltimes \partial \Omega)\big)\arrow{d}\\
        K_{1}\big(C^{*}_{r}(G\ltimes \partial \Omega)\big)\arrow{u}\ &K_{1}\big(C^{*}_{r}(G\ltimes \Omega_{b, \infty})\big) \arrow{l}\ &K_{1}\big(C^{*}_{r}\big(G\ltimes (\Omega_{b, \infty} \setminus \partial \Omega)\big)\big)\arrow{l}.\\
    \end{tikzcd}
\end{equation}
\pars

We have $\Omega_{b, \infty}\setminus \partial \Omega=\Omega_{\infty}\setminus \Omega_{A,\ \infty}$, and thus
$$C^{*}_{r}\big(G\ltimes (\Omega_{b, \infty} \setminus \partial \Omega)\big)=C^{*}_{r}\big(G\ltimes (\Omega_{\infty}\setminus \Omega_{A,\ \infty})\big)\cong \mathcal{K}\otimes C(\mathbb{T})$$
as in the proof of Proposition~\ref{prop:K:OmegaAInf}. This implies
$$K_{0}\big(C^{*}_{r}\big(G\ltimes (\Omega_{b, \infty} \setminus \partial \Omega)\big)\big) = \mathbb{Z} [1_{(\beps,\ \chi_{b^{\infty}})}]_0 \ \text{and}\ K_{1}\big(C^{*}_{r}\big(G\ltimes (\Omega_{b, \infty} \setminus \partial \Omega)\big)\big) = \mathbb{Z} [u]_1,$$
where $u=1_{(b,\ \chi_{b^{\infty}})} + 1 - 1_{(\beps,\ \chi_{b^{\infty}})}$. By Proposition~\ref{prop:K:Omegab}, we have
$$K_{0}\big(C^{*}_{r}(G\ltimes \Omega_{b, \infty}) \big) = \mathbb{Z} [1]_0 \ \text{and}\ K_{1}\big(C^{*}_{r}(G\ltimes \Omega_{b, \infty}) \big)  = \mathbb{Z} [1_{\gekl{b} \times \Omega_{b, \infty}}]_1.$$
Plugging this into \eqref{e:K:OmegaBInf}, we obtain
$$
    \begin{tikzcd}
        \mathbb{Z}\arrow{r}{K_{0}(\iota)}\ &\mathbb{Z}\arrow{r}\ &K_{0}\big(C^{*}_{r}(G\ltimes \partial \Omega)\big)\arrow{d}\\
        K_{1}\big(C^{*}_{r}(G\ltimes \partial \Omega)\big)\arrow{u}\ &\mathbb{Z} \arrow{l}\ &\mathbb{Z}\arrow{l}[swap]{K_{1}(\iota)},\\
    \end{tikzcd}
$$
\nopar

where $\iota$ is the inclusion map from $C^{*}_{r}\big(G\ltimes (\Omega_{b, \infty} \setminus \partial \Omega)\big)$ into $C^{*}_{r}\big(G\ltimes \Omega_{b, \infty})$. Since
$$1_{(\beps,\ \chi_{b^{\infty}})}=1-\sum_{e \in A,\ 0\leq j\leq n_e-1}1_{\{b^{j}e\}\times \Omega_{b, \infty}}1^{*}_{\{b^{j}e\}\times \Omega_{b, \infty}}\in C^*_r(G\ltimes \Omega_{b, \infty}),$$
we have $K_{0}(\iota)([1_{(\beps,\ \chi_{b^{\infty}})}]_{0})=(1-\sum_{e \in A} n_e)[1]_0$.
\pars

Define, for each $e \in A$,
$$u_e \defeq 1+\sum_{\ 0\leq j\leq n_e-1}(1_{\{b\}\times \Omega_{b, \infty}}-1)1_{\{b^{j}e\}\times \Omega_{b, \infty}}1^{*}_{\{b^{j}e\}\times \Omega_{b, \infty}}\in C^*_r(G\ltimes \Omega_{b, \infty}).$$
Then we have $u\cdot\prod_{e \in A}u_e=1_{\{b\}\times \Omega_{b, \infty}}$. Define
$$\Pi_{e,j}:=1_{\{b^{j}a_{i}\}\times \Omega_{b, \infty}}1^{*}_{\{b^{j}a_{i}\}\times \Omega_{b, \infty}},$$
and let $u'_e = u_e \sum_{j} \Pi_{e,j}$. Then $u_e = u'_e + (1 - \sum_{j} \Pi_{e, j})$. With respect to the pairwise orthogonal projections $\Pi_{e, j}$, $0 \leq j < n_e$, we have
\[u'_e=
\left( \begin{array}{c|c}
0 & 1_{\{b^{n_e}\}\times \Omega_{b, \infty}}\\
  \hline
\begin{matrix}
    1 & & \\
    & \ddots & \\
    & & 1
  \end{matrix} & 0\\
\end{array}
\right).
\]
Multiplying $u'_{e}$ by the permutation matrix
\[
\left( \begin{array}{c|c}
0 & \begin{matrix}
    1 & & \\
    & \ddots & \\
    & & 1
  \end{matrix}\\
  \hline
1 & 0\\
\end{array}
\right)
\]
from the  right hand side, we get the following diagonal matrix
\[u''_{e}=
\left( \begin{array}{c|c}
 1_{\{b^{n_e}\}\times \Omega_{b, \infty}} & 0\\
  \hline
 0 &
\begin{matrix}
    1 & & \\
    & \ddots & \\
    & & 1
  \end{matrix}\\
\end{array}
\right).
\]
Therefore, $u'_e$ is homotopic to $u''_e$, and hence $u_e$ is homotopic to
\begin{align*}
 & 1+(1_{\{b^{n_e}\}\times \Omega_{b, \infty}}-1)1_{\{e\}\times \Omega_{b, \infty}}1^{*}_{\{e\}\times \Omega_{b, \infty}}\\
 = \ & (1-1_{\{e\}\times \Omega_{b, \infty}}1^{*}_{\{e\}\times \Omega_{b, \infty}})+1_{\{e\}\times \Omega_{b, \infty}}1_{\{b^{\sgn(e)m_e}\}\times \Omega_{b, \infty}}1^{*}_{\{e\}\times \Omega_{b, \infty}}.
\end{align*}

It follows from Lemma~\cite[Lemma 4.6.2]{HR} that $[u_e]_1 = [1_{\{b^{\sgn(e)m_e}\}\times \Omega_{b, \infty}}]_1$. That is, $[u_e]_{1}=\sgn(e)m_e [1_{\{b\}\times \Omega_{b, \infty}}]_1$ and $[u]_{1}= (1-\sum_{e \in A} \sgn(e) m_e) [1_{\{b\}\times \Omega_{b, \infty}}]_1$. 

It is now straightforward to deduce the desired result about $K_*(C^*_r(G \ltimes \partial \Omega))$.
\eproof
\pars

\section{Classification of boundary quotients, and families of Cartan subalgebras in UCT Kirchberg algebras}

\subsection{Classification of boundary quotients}

Assume that $G_v \subseteq (\Rz,+)$ for all $v \in V$, and $\# V > 1$ or $\# A > 0$, and that conditions (LCM) from Definition~\ref{def:LCM} and (D) from Definition~\ref{def:D} are satisfied. The following is an immediate consequence of Theorem~\ref{thm:topfree}~(ii) and Theorems~\ref{thm:CPTnuc} and \ref{thm:K}.
\btheo
\label{thm:bC}
$\partial C^*_{\lambda}(P) = C^*_r(G \ltimes \partial \Omega)$ is a UCT Kirchberg algebra if the following two conditions are satisfied.
\nopar

\begin{enumerate}
\item[(TF)] One of the following holds:
\begin{enumerate}
\item[(a)] There exists $e \in T$ with $P_e = \gekl{\beps}$,
\item[(b)] For all $e \in T$, $P_e \neq \gekl{\beps}$, $\# A > 0$ and there exists $e \in A$ with $k_e \not\mid l_e$,
\item[(c)] For all $e \in T$, $P_e \neq \gekl{\beps}$, $\# A > 0$, for all $e \in A$, $k_e \mid l_e$ and $\big( \bigcap_{e \in A} \langle b_{t(e)}^{k_e n_e} \rangle \big) \cap \big( \bigcap_{v \in V} G_v \big) = \gekl{\beps}$.
\end{enumerate}
\item[(N)] For all $T' \subseteq T$ with $P_e \neq \gekl{\beps}$ for all $e \in T'$, either $T'$ consists of a single vertex or $T'$ consists of exactly two vertices $v, w$ and one pair of edges $e, \bar{e}$ with $o(e) = v$, $t(e) = w$, such that $P_v \cong \Zzg$, $P_w \cong \Zzg$, and the embeddings $(\cdot)^{\bar{e}}$, $(\cdot)^e$ are both given by $\Zzg \to \Zzg, \, z \ma 2z$.
\end{enumerate}

In that case, the K-theory of $\partial C^*_{\lambda}(P) = C^*_r(G \ltimes \partial \Omega)$ is given as follows:
\begin{enumerate}
\item[(i)] If there exists $e \in T$ with $P_e = \gekl{\beps}$ and one of the following holds:
\begin{enumerate}
\item[(i$_1$)] There exists $v \in V$ such that $G_v$ is dense in $\Rz$,
\item[(i$_2$)] $\# V = \infty$,
\item[(i$_3$)] $\# A_+ = \infty$,
\end{enumerate}
then $K_*(\partial C^*_{\lambda}(P)) \cong K_*(\Cz)$.
\item[(ii)] If there exists $e \in T$ with $P_e = \gekl{\beps}$, $P_v \cong \Zzg$ for all $v \in V$, $\# V < \infty$ and $\# A_+ < \infty$, 

then
$K_*(\partial C^*_{\lambda}(P)) \cong (\Zz/N,1, \gekl{0})$, where $N = \frac{1}{2} \# \menge{e \in T}{P_e = \gekl{\beps}}$.
\item[(iii)] If $P_e \neq \gekl{\beps}$ for all $e \in T$ and one of the following holds:
\begin{enumerate}
\item[(iii$_1$)] $\# V > 1$,
\item[(iii$_2$)] $\# V = 1$, $V = \gekl{v}$, and $G_v$ is dense in $\Rz$,
\item[(iii$_3$)] $\# V = 1$, $V = \gekl{v}$, $P_v \cong \Zzg$, $\# A = \infty$,
\end{enumerate}
then $K_*(\partial C^*_{\lambda}(P)) \cong K_*(C^*_{\lambda}(G_T))$.
\item[(iv)] If $\# V = 1$, $V = \gekl{v}$, $P_v \cong \Zzg$, $\# A < \infty$ and $\# A_+ > 0$, then $\sum_{e \in A} n_e \neq 1$, and we have the following:
\begin{enumerate}
\item[(iv)$_1$] If $\sum_{e \in A} \sgn(e) m_e \neq 1$, then $K_*(\partial C^*_{\lambda}(P)) \cong (\Zz / (1 - \sum_{e \in A} n_e), 1, \Zz / (1 - \sum_{e \in A} \sgn(e) m_e))$.
\item[(iv)$_2$] If $\sum_{e \in A} \sgn(e) m_e = 1$, then $K_*(\partial C^*_{\lambda}(P)) \cong (\Zz / (1 - \sum_{e \in A} n_e) \oplus \Zz, (1,0), \Zz)$.
\end{enumerate}
\end{enumerate}
\etheo

By the Kirchberg-Phillips classification theorem \cite{KP,Phi}, if conditions (TF) and (N) are satisfied, then $\partial C^*_{\lambda}(P)$ is completely classified by K-theory. In particular, we obtain the following
\bcor
\label{cor:bC}
In case (i) of Theorem~\ref{thm:bC}, we have $\partial C^*_{\lambda}(P) \cong \cO_{\infty}$, and in case (ii) of Theorem~\ref{thm:bC}, we have $\partial C^*_{\lambda}(P) \cong \cO_{N+1}$.
\ecor
\pars

\subsection{Families of Cartan subalgebras in UCT Kirchberg algebras}

In the following, we call the up to homeomorphism unique totally disconnected, second countable, locally compact non-compact Hausdorff space without isolated points the non-compact locally compact Cantor space. Theorem~\ref{thm:bC} has the following application.
\btheo
\label{thm:CartanInKirchberg}
Let $A$ be a UCT Kirchberg algebra. For every abelian, torsion-free, finite rank group $\Gamma$ which is not free abelian, there exists a Cartan subalgebra $B_{\Gamma}$ of $A$ such that $\Spec B_{\Gamma}$ is homeomorphic to the Cantor space if $A$ is unital and to the non-compact locally compact Cantor space if $A$ is not unital, and, for all such groups $\Gamma$ and $\Lambda$, $(A,B_{\Gamma}) \cong (A,B_{\Lambda})$ implies $\Gamma \cong \Lambda$.
\etheo
\nopar

\bproof
Given $\Gamma$ as in the theorem, we view $\Gamma$ as a subgroup of $(\Rz,+)$ by choosing an embedding $\Gamma \into \Rz$ and let $\Gamma^+ \defeq \Gamma \cap [0,\infty)$, $P \defeq \Zzg * \Gamma^+$, and $\cG_{\Gamma} \defeq (\Zz * \Gamma) \ltimes \Omega_P$. Note that $\cG_{\Gamma}$ depends on the choice of the embedding $\Gamma \into \Rz$. Elements of $\Omega_P$ are in one-to-one correspondence to finite or infinite reduced words in $\Zzg \cup \Gamma^+$. It is now straightforward to see that $\Stab(\cG_{\Gamma}) = \gekl{\gekl{\beps}, \Zz, \Gamma}$. Here and in the sequel, given a groupoid $\cG$, $\Stab(\cG)$ denotes the set of isotropy groups $\cG_x^x$, $x \in \cG^{(0)}$, up to isomorphism.
\pari

Now let $A$ be a UCT Kirchberg algebra. Let $\cG_A$ be the groupoid model for $A$ as in \cite[\S~5]{LR}, i.e., $\cG_A$ is an {\'e}tale locally compact Hausdorff groupoid such that $A \cong C^*_r(\cG_A)$. Moreover, the construction in \cite[\S~5]{LR} yields that $\Stab(\cG_A) \subseteq \gekl{\gekl{\beps},\Zz,\Zz^2}$. We then obtain, using Corollary~\ref{cor:bC}, that $C^*_r(\cG_A \times \cG_{\Gamma}) \cong A \otimes \cO_{\infty} \cong A$. Let $B_{\Gamma} \subseteq A$ be the image of $C^*_r(\cG_A^{(0)} \times \cG_{\Gamma}^{(0)})$ under this isomorphism. It is an immediate consequence of the results in \cite[\S~5]{LR} that $\Spec B_{\Gamma}$ is of the required homeomorphism type. Moreover, we have (up to isomorphism) $\gekl{\gekl{\beps},\Zz,\Gamma} \subseteq \Stab(\cG_A \times \cG_{\Gamma}) \subseteq \gekl{\gekl{\beps},\Zz,\Zz^2,\Zz^3,\Gamma,\Zz \times \Gamma, \Zz^2 \times \Gamma}$. It follows that $\Gamma$ can be characterized, up to isomorphism, as the group in $\Stab(\cG_A \times \cG_{\Gamma})$ which is not free abelian and of minimal rank. Now if $(A,B_{\Gamma}) \cong (A,B_{\Lambda})$, then $\Stab(\cG_A \times \cG_{\Gamma}) = \Stab(\cG_A \times \cG_{\Lambda})$ and thus $\Gamma \cong \Lambda$.
\eproof
\pars

\end{document}